\documentclass[12pt]{amsart} 

\usepackage{amsmath,amssymb,color,amsthm,epsfig,extarrows,tikz,enumitem,
hyperref,url}

\theoremstyle{definition}
\newtheorem{prop}{Proposition}[section]
\newtheorem{algor}[prop]{Algorithm} 
\newtheorem{thm}[prop]{Theorem}
\newtheorem{rem}[prop]{Remark}
\newtheorem{lem}[prop]{Lemma}
\newtheorem{dfn}[prop]{Definition}

\newtheorem{mahler}{Mahler's Theorem}

\newtheorem{hensel}{Hensel's Lemma}

\newtheorem{ex}[prop]{Example}
\newtheorem{cor}[prop]{Corollary}

\newtheorem{not*}[prop]{Notation}

\newcommand{\F}{\mathbb{F}}
\newcommand{\bF}{\overline{\mathbb{F}}}
\newcommand{\N}{\mathbb{N}}
\newcommand{\Q}{\mathbb{Q}}
\newcommand{\R}{\mathbb{R}}
\newcommand{\C}{\mathbb{C}}

\newcommand{\Z}{\mathbb{Z}}

\newcommand{\res}{\mathrm{Res}}
\newcommand{\ba}{\bar{a}}

\newcommand{\barf}{\bar{f}}

\newcommand{\cD}{\mathcal{D}}

\newcommand{\cT}{\mathcal{T}}

\newcommand{\np}{\mathbf{NP}}
\newcommand{\zpp}{\mathbf{ZPP}}

\newcommand{\norm}[1]{\left|#1\right|_p}
\newcommand{\ord}{\operatorname{ord}}
\newcommand{\newt}{\operatorname{Newt}}
\newcommand{\anewt}{\operatorname{Newt}_\infty}

\newcommand{\tf}{\tilde{f}}
\newcommand{\tg}{\tilde{g}}
\newcommand{\tq}{\tilde{q}}
\renewcommand{\qed}{$\blacksquare$}
\newcommand{\dia}{$\diamond$}
\newcommand{\eps}{\varepsilon}

\newcommand{\floor}[1]{\left\lfloor #1 \right\rfloor}
\newcommand{\ceil}[1]{\left\lceil #1 \right\rceil}

\title[Root Repulsion and Faster Solving]{\mbox{}\\ 
\vspace{-1.3in}Root Repulsion and Faster Solving for Very Sparse Polynomials  
Over $p$-adic Fields}

\author{J.\ Maurice Rojas}
\thanks{Partially supported by NSF grants CCF-1900881 and 
CCF-1409020. A much shorter, preliminary version of this work appeared 
in the proceedings of the conference ISSAC 2021 (July 19--23, virtual event) 
\cite{issacversion}.}
\email{jmauricerojas@gmail.com}
\author{Yuyu Zhu}
\email{yuyu.zhu1213@gmail.com}
\address{Texas A\&{}M University, TAMU 3368, College Station, Texas \ 
77843-3368 }

\begin{document} 

\begin{abstract}
For any {\em fixed} field $K\!\in\!\{\Q_2,\Q_3,\Q_5,\ldots\}$, we prove that 
all polynomials $f\!\in\!\Z[x]$ with exactly $3$ (resp.\ $2$) monomial terms, 
degree $d$, and all coefficients having absolute value at most $H$, 
can be solved over $K$ within deterministic time $\log^{4+o(1)}(dH)\log^3(d)$ 
(resp.\ $\log^{2+o(1)}(dH)$) in the classical Turing model: Our 
underlying algorithm correctly counts the
number of roots of $f$ in $K$, and for 
each such root generates an approximation in $\Q$ with logarithmic height 
$O(\log^2(dH)\log(d))$ that converges at a rate of 
$O\!\left((1/p)^{2^i}\right)$ 
after $i$ steps of Newton iteration. We also prove significant speed-ups in 
certain settings, a minimal spacing bound of $p^{-O(p\log^2_p(dH)\log d)}$ for 
distinct roots in $\C_p$, and even stronger repulsion when 
there are nonzero degenerate roots in $\C_p$: $p$-adic distance 
$p^{-O(\log_p(dH))}$. On the other hand, we prove that there is an explicit 
family of tetranomials with distinct nonzero roots in $\Z_p$ indistinguishable 
in their first $\Omega(d\log_p H)$ most significant base-$p$ digits.   
\end{abstract}

\keywords{p-adic, Hensel, Newton, iteration, trinomial, approximate, root 
counting}

\maketitle

\tableofcontents

\vspace{-.4in} 
\section{Introduction} 
Solving polynomial equations over the $p$-adic rational numbers $\Q_p$ 
underlies many important computational questions in number theory 
(see, e.g., \cite{padicfastmult,balakrishnan,kedlaya,dzb})  
and is close to\linebreak 
\scalebox{.955}[1]{applications in coding theory (see, e.g., \cite{blq13}). 
Furthermore, the complexity of 
solving {\em structured}}\linebreak 
\scalebox{.947}[1]{equations --- such as those with a fixed number of 
monomial terms or invariance with respect to}\linebreak 
\scalebox{.965}[1]{a group action --- arises naturally in many computational 
geometric applications and is closely}\linebreak 
related to a deeper understanding of circuit 
complexity (see, e.g., \cite{koiranrealtau}). 
So we will study how sparsity affects the complexity of separating and 
approximating roots in $\Q_p$. Unless stated otherwise, all 
$O$-constants and $\Omega$-constants are effective and absolute. 

Recall that thanks to 17th century work of Descartes, and 20th century 
work of Lenstra \cite{len99} and Poonen \cite{poonen},    
it is known that univariate polynomials with exactly $t$ monomial 
terms have at most $t^{O(1)}$ roots in a fixed field $K$ {\em only} 
when $K$ is $\R$ or a finite algebraic extension of $\Q_p$ for some prime 
$p\!\in\!\N$. 
We'll use $|\cdot|_p$ (resp.\ $|\cdot|$) for the absolute value on the 
$p$-adic complex numbers $\C_p$ normalized so that $|p|_p\!=\!\frac{1}{p}$ 
(resp.\ the standard absolute value on $\C$). Recall also that for any function 
$f$ analytic on $K$, the 
corresponding {\em Newton endomorphism} is $N_f(z):=z-\frac{f(z)}{f'(z)}$, 
and the corresponding sequence of {\em Newton iterates} of a {\em start-point} 
$z_0\!\in\!K$ is the sequence $(z_i)^\infty_{i=0}$ where 
$z_{i+1}\!:=\!N_f(z_i)$ for all $i\!\geq\!0$. 

Our first main result is that we can efficiently count the roots of 
univariate trinomials in $\Q_p$, {\em and} find succinct start-points in $\Q$ 
under which Newton iteration converges quickly to all the roots in $\Q_p$. 
We use $\#S$ for the cardinality of a set $S$.  

\begin{thm} 
\label{thm:big} {\em 
\vbox{For any prime $p$ and a trinomial $f\!\in\!\Z[x]$ with 
degree $d$ and all its coefficients having (Archimedean) absolute value
$\leq\!H$, we can find in deterministic time\\ 
\mbox{}\hfill $O\!\left(p^3\log^4(dH)\log^3_p(d)
\log(p\log(dH))\right)$\hfill\mbox{}\\ 
a set $\{\frac{\alpha_1}{\beta_1},\ldots,\frac{\alpha_m}{\beta_m}\}\subset\!\Q$ of cardinality 
$m\!=\!m(p,f)$ such that:\\ 
\mbox{}\hspace{.4cm} 
(1) For all $j$ we have $\alpha_j\!\neq\!0 
\Longrightarrow \log|\alpha_j|,\log|\beta_j|=O\!\left(p
\log^2_p(dH)\log(d)\right)$.\\  
\mbox{}\hspace{.4cm} 
(2) $z_0\!:=\!\alpha_j/\beta_j \Longrightarrow 
f$ has a root $\zeta_j\!\in\!\Q_p$ with sequence of Newton iterates
satisfying\\ 
\mbox{}\hspace{1.3cm}$|z_i-\zeta_j|_p\!\leq\!(1/p)^{2^i}|z_0-\zeta_j|_p$ 
for all $i,j\!\geq\!1$.\\  
\mbox{}\hspace{.4cm} 
(3) $m\!=\!\#\{\zeta_1,\ldots,\zeta_m\}$ is exactly 
the number of roots of $f$ in $\Q_p$.}}
\end{thm} 

\noindent 
We prove Theorem \ref{thm:big} in Section \ref{sub:mainalgor} 
via Algorithm \ref{algor:trinosolqp} there. The dependence on $p$ can 
be lowered significantly in certain natural settings, e.g., restricting to 
roots of the form $p^j+O(p^{j+1})$, making mild assumptions on the gcd 
of the exponents, or assuming 
the presence of degenerate roots in $\C^*_p$: See Corollaries 
\ref{cor:faster}, \ref{cor:degen}, and \ref{cor:final} 
below. An analogue of Theorem \ref{thm:big} also holds for $K\!=\!\R$ 
and will be presented in a sequel to this paper. We call a $z_0\!\in\!\Q_p$ 
satisfying the convergence condition from Theorem \ref{thm:big} {\em an 
approximate root of $f$ (in the sense of Smale\footnote{This 
terminology has only been applied over $\C$ so far \cite{smale}, so we take 
the opportunity here to extend it to the $p$-adic rationals. Note that 
we do not restrict $\zeta$ to be non-degenerate.}{\em )}, 
with associated true root $\zeta$}. This type 
of convergence provides an efficient encoding of an approximation that 
can be quickly tuned to any desired accuracy. 

\begin{rem} {\em Defining the {\em input size} of a univariate polynomial 
$f(x)\!:=\!\sum^t_{i=1} c_i x^{a_i}\!\in\!\Z[x]$ as 
$\sum^t_{i=1}\log((|c_i|+2)(|a_i|+2))$ we see that Theorem \ref{thm:big} 
implies that one can solve univariate\linebreak 
\scalebox{.91}[1]{trinomial equations, over any {\em fixed} 
$p$-adic field, in deterministic time polynomial in the input size. \dia}} 
\end{rem} 
\begin{rem} {\em 
Efficiently solving univariate $t$-nomial equations over $K$ 
in the sense of Theorem \ref{thm:big} is easier for $t\!\leq\!2$: 
The case $t\!=\!1$ is clearly trivial (with $0$ the only possible root) 
while the case $(K,t)\!=\!(\R,2)$ is implicit in work 
on computer arithmetic from the 1970s (see, e.g., \cite{borwein}). 
We review the case $(K,t)\!=\!(\Q_p,2)$ with $p$ prime in Corollary 
\ref{cor:binomod} and Theorem \ref{thm:binoqp} of Section 2 below. \dia}  
\end{rem} 

Despite much work on factoring univariate polynomials over $\Q_p$ 
(see, e.g., \cite{CG00,GNP12,bns13,blq13}), all known general algorithms for 
solving (or even just counting the solutions of) arbitrary degree $d$ 
polynomial equations over $\Q_p$ have complexity exponential in $\log d$.  
So Theorem \ref{thm:big} presents a significant new speed-up, and greatly 
improves an earlier complexity bound (membership in $\np$, for detecting roots 
in $\Q_p$) from \cite{airr}. We'll see in Sections \ref{sec:central} and 
\ref{sec:trinosolqp} how our speed-up depends on $p$-adic Diophantine 
approximation \cite{yu94,yu07}. Another key new ingredient in proving Theorem 
\ref{thm:big} is an efficient encoding of roots in $\Z/(p^k)$ from 
\cite{DMS19,krrz19}, with an important precursor in \cite{blq13}.   

\subsection{Dependence on $p$} 
While there are certainly number-theoretic algorithms with deterministic 
complexity having dependence $(\log p)^{O(1)}$ on an input prime $p$, solving 
sparse polynomial equations in just one variable 
over $\Q_p$ appears to have much larger complexity\linebreak 
\scalebox{.95}[1]{as a function of $p$. There is a naive reason (R1 below), 
and a subtle reason (R2 below), for this:}  

\smallskip
\noindent 
{\bf R1.} {\em Whereas a binomial has at most $3$ roots in $\R$ (e.g., 
$x^3-x$), a binomial can have as many as $\max\{p,3\}$ roots in $\Q_p$ (e.g., 
$x^{\max\{p,3\}}-x$). Furthermore, trinomials have at most $5$, $7$, $9$, or 
$3p-2$ roots in $K$, according as $K$ is $\R$, $\Q_2$ \cite{len99}, 
$\Q_3$ \cite{zhuthesis}, or $\Q_p$ with $p\!\geq\!5$ \cite{ak11,barba}, and 
each bound is sharp. \dia} 

\smallskip
The most natural $p$-adic analogue of a positive real number is a 
$p$-adic rational number {\em with most significant digit $1$}, i.e., a number 
of the form $p^j+O(p^{j+1})$. Restricting to such roots then 
cuts the aforementioned root cardinality bounds to $2$, $6$, $4$, and $3$ 
(respectively over $\R$, $\Q_2$, $\Q_3$, or $\Q_p$ with $p\!\geq\!5$), 
and yields a significant speed-up for solving that we detail  
in Corollary \ref{cor:faster} below. 
Alternatively, rather than restricting digits of roots, one can observe 
that trinomials over $\Z$ with many roots in $\Q_p$ are (arguably) rare. 
This enables another significant speed-up to our main algorithm for  
``most'' choices of exponents. 
\begin{cor} {\em
\label{cor:faster}
Following the notation of Theorem \ref{thm:big}, we can 
lower the deterministic time complexity bound to 
$O\!\left(p^2\log^4(dH)\log^3_p(d)\log(p\log(dH))\right)$, in either 
of following settings: (1)    
we only seek roots of the form $p^j+O(p^{j+1})$, or (2) 
we assume that the exponents are $\{0,a_2,a_3\}$ with 
$\gcd(a_2a_3(a_3-a_2),(p-1)p)\!\leq\!2$. 
In either case, the stated worst-case height bounds for the approximate roots 
remain the same.} 
\end{cor} 

\noindent 
We prove Corollary \ref{cor:faster} in Section \ref{sub:faster}, and leave 
average-case speed-ups, where one averages over {\em coefficients}, for future 
work. 
It follows from our framework that the speed-ups from Corollary 
\ref{cor:faster} continue to hold 
(modulo a multiple of $C^{O(1)}$) under softer 
assumptions like (a) restricting to roots with most significant digit 
in some cardinality $C$ subset of $\{1,\ldots,p-1\}$ or (b) assuming 
$\gcd(a_2a_3(a_3-a_2),(p-1)p)\!\leq\!C$. So our assumptions  
above are more restrictive merely for the sake of simplifying our exposition.  

\smallskip
\noindent 
{\bf R2.} {\em Approximating square-roots of $p$-adic integers not divisible 
by $p$, within accuracy $1$, is equivalent to finding square-roots in 
the finite field $\F_p$. The latter problem is {\em still} not known to 
be doable in deterministic time polynomial in $\log p$, even though 
the decision version is doable in deterministic time $\log^{2+o(1)} p$ 
(see, e.g., \cite{shoup,bs,poonenzeta}). Furthermore, it remains unknown how 
to find just a {\em single} $d$th root of a $d$th power in $\F^*_p$ in time 
$(\log(p)+\log d)^{O(1)}$, even if randomness is allowed (see, e.g., 
\cite{amm,caofan,chokwonlee}). \dia}  

\smallskip
Parallel to R2, even if one only wants to approximate a single 
root in $\Q_p$ of a trinomial, the minimal {\em currently} provable initial 
accuracy needed to make Newton iteration converge quickly appears to 
have {\em quasi-linear} dependence on $p$. This is because our key valuation  
bounds (see Section \ref{sec:central}) currently hinge on estimates for 
{\em linear forms in $p$-adic logarithms} \cite{bakerabc,yu94,yu07}, and 
further improvements to the latter estimates appear quite difficult. 

\subsection{Repulsion, and the Separation Chasm at Four Terms} 
\label{sub:tetrasep} 
\scalebox{.9}[1]{The $p$-adic rational roots of}\linebreak 
\scalebox{.97}[1]{sparse polynomials can range from well-separated 
to tightly spaced, already with just $4$ terms.}   
\begin{thm} 
{\em \label{thm:tetra}
Consider the family of tetranomials \\ 
\mbox{}\hfill $\displaystyle{f_{d,\eps}(x):=x^d - \eps^{-2h}x^2 
+ 2\eps^{-(h+1)}x - \eps^{-2}}$\hfill\mbox{}\\ 
with $h\!\in\!\N$, $h\!\geq\!3$, and $d\!\in\!\{4,\ldots,\lfloor e^h 
\rfloor\}$ even. Let $H\!:=\!\max\{\eps^{\pm 2h}\}$.
Then $f_{d,\eps}$ has distinct nonzero roots $\zeta_1,\zeta_2$ in the 
open unit disk of $K$ (centered at $0$) with 
$|\log|\zeta_1-\zeta_2|_p|\!=\!\Omega(d\log H)$ or
$|\log|\zeta_1-\zeta_2||\!=\!\Omega(d\log H)$, according as 
$(K,\eps)\!=\!(\Q_p,p)$ or $(K,\eps)\!=\!(\R,1/2)$. In 
particular, the coefficients of $p^{2h}f_{d,p}$  
all lie in $\Z$ and have $O(\log_p H)$ base-$p$ digits, 
and we need $\Omega(d\log_p H)$ many base-$p$ digits to 
distinguish the roots of $f$ in $\Z_p$. }  
\end{thm}

\vspace{-.1cm} 
\noindent
We prove Theorem \ref{thm:tetra} in Section \ref{sec:tetra}, where we will 
also see in Remark \ref{rem:tetra} that the basin of attraction for a root of 
$f_{d,p}$ in $\Q_p$ (under the Newton endomorphism $N_{f_{d,p}}$) can be 
exponentially small in $\log d$ as well. The special case $K\!=\!\R$ of Theorem 
\ref{thm:tetra} was derived earlier (in different notation) by Mignotte 
\cite{mig95}. (See also \cite{sag14}.) The cases $K\!=\!\Q_p$ with $p$ 
prime appear to be new, and our proof unifies the Archimedean and 
non-Archimedean cases via tropical geometry \cite{aknr}. Approximating roots 
in $\Q_p$ in average-case time sub-linear in $d$ for tetranomials 
(where one averages over the coefficients but fixes the exponents) is 
thus an intriguing open problem.  

Mignotte used the tetranomial $f_{d,1/2}$ in \cite{mig95} 
to show that an earlier root separation bound of Mahler 
\cite{mah64}, for {\em arbitrary} degree $d$ polynomials in $\Z[x]$, 
is asymptotically near-optimal. We recall the following paraphrased version: 

\begin{mahler} {\em 
Suppose $f\!\in\!\Z[x]$ has degree $d\!\geq\!2$, all coefficients of 
(Archimedean) absolute value at most $H$, and is irreducible in $\Z[x]$. Let 
$\zeta_1,\zeta_2\!\in\!\C$ be distinct roots of $f$. Then 
$|\zeta_1-\zeta_2|\!>\!\frac{\sqrt{3}}{(d+1)^{d+\frac{1}{2}}H^{d-1}}$. 
In particular, $|\log|\zeta_1-\zeta_2||\!=\!O(d\log(dH))$. \qed} 
\end{mahler} 

\noindent 
The very last statement is actually a small addendum, making use of the 
following classic fact: The complex roots of an $f$ as above lie in an open 
disk, centered at the origin, of radius $2H$ (see, e.g., \cite[Ch.\ 8]{rs} or 
Theorem \ref{thm:newt} in Section \ref{sub:newt} below). It is straightforward 
to prove\linebreak 
\scalebox{.96}[1]{an analogue of Mahler's bound, of the same asymptotic order 
for $|\log|\zeta_1-\zeta|_p|$, for roots in $\C_p$.}   

Our new algorithmic results are enabled by our third and final main
result: Mahler's bound can be dramatically improved for the roots of 
{\em tri}nomials in $\C_p$.   
\begin{thm} {\em 
\label{thm:tri} 
Suppose $p$ is prime and $f\!\in\!\Z[x]$ has exactly $3$ monomial terms, 
degree $d$, and all its coefficients have  
(Archimedean) absolute value at most $H$. Let $\zeta_1,\zeta_2\!\in\!\C_p$ be 
distinct roots of $f$. Then $\log H\!\geq\!\log|\zeta_1-\zeta_2|_p\!\geq\!
-O\!\left(p\log^2(dH)\log_p d\right)$. Furthermore, if $f$ has a 
degenerate root in $\C^*_p$, then the last lower bound can be sharpened to 
$-O(\log(dH))$. }    
\end{thm} 

\vspace{-.1cm} 
\noindent 
We prove Theorem \ref{thm:tri} in Section \ref{sec:trisepqp}. 
Theorem \ref{thm:tri} provides a $p$-adic analogue of a separation 
bound of Koiran for complex roots of trinomials \cite{koiransep}. 
As to whether our lower bound is optimal, there are recent examples from
\cite{fgr} showing that $\log|\zeta_1-\zeta_2|_p\!=\!-\Omega(\log\max\{d,H\})$ 
can occur. However, we are unaware of any examples exhibiting  
$\log|\zeta_1-\zeta_2|_p\!=\!-\Omega(p^\eps)$ for some $\eps\!>\!0$.   
Asymptotically optimal separation bounds,
over both $\C_p$ and $\C$, are already known for binomials and 
we review these bounds in Section \ref{sec:bisep}.  

The presence of degenerate roots appears to not only increase the repulsion 
of roots for trinomials but also speed up their approximation: 
\begin{cor} \label{cor:degen} {\em 
Following the notation of Theorem \ref{thm:big}, if $f$ has a degenerate 
root in $\C^*_p$, then we can find, in deterministic time 
$O\!\left(p[p^{1/2}\log^2(p)+\log^2(dH)\log(dp)\log\log(dpH)]\right)$, 
or Las Vegas randomized time 
$O\!\left(p\left[\log^{2+o(1)}(p)+\log^2(dH)\log(dp)\log\log(dpH)\right]
\right)$, a set of approximate roots in the sense of Smale, each in $\Q$ 
and with logarithmic height $O(\log(dH))$, with distinct associated true 
roots having union the zero set of $f$ in $\Q_p$.}  
\end{cor} 

\noindent 
We prove Corollary \ref{cor:degen} in Remark \ref{rem:main} of 
Section \ref{sub:mainalgor} below. 
It is not yet clear whether significantly better bounds for root spacing and 
root approximation can hold in complete generality: The apparent 
improvements implied by the presence of degenerate roots could   
just be a side-effect of our underlying techniques. Curiously, a similar 
``repulsion from degeneracy'' phenomenon also occurs in the (Archimedean)  
setting of roots in $\C$: See 
\cite[Proof of Thm.\ 18]{koiransep}. 

\subsection{Previous Complexity and Sparsity Results} 
Deciding the existence of roots 
over $\Q_p$ for univariate polynomials with an {\em arbitrary} number 
of monomial terms is already $\np$-hard with respect 
to randomized ($\zpp$, a.k.a.\ Las Vegas) reductions \cite{airr}. 
On the other hand, detecting roots over $\Q_p$ for $n$-variate $(n+1)$-nomials 
is known to be doable in $\np$ \cite{airr}. Speeding this up to 
polynomial-time, even for $n\!=\!2$ and fixed $p$, hinges upon detecting roots 
in $(\Z/(p^k))^2$ for bivariate trinomials of degree $d$ 
in time $(k+\log d)^{O(1)}$. The latter problem remains open, but some 
progress has been made in author Zhu's Ph.D.\ thesis \cite{zhuthesis}. 

On a related note, counting points on trinomial curves over the prime 
fields $\F_p$ in time $(\log(pd))^{O(1)}$ remains a challenging open question. 
Useful quantitative estimates in this direction were derived in  
\cite{huavandiver} and revisited via real quadratic optimization in 
\cite{avendanomorales}.  

\section{Background} \label{sec:back} 
Recall that the famous {\em Ultrametric Inequality} states that 
for any $\alpha,\beta\!\in\!\C_p$ we have\linebreak 
$\ord_p(\alpha\pm \beta)\!\geq\!\min
\{\ord_p\alpha,\ord_p\beta\}$. (Equivalently: $|\alpha\pm \beta|_p
\!\leq\!\max\{|\alpha|_p,|\beta|_p\}$.) We will frequently use (without further 
mention) this inequality, along with its natural implication 
$\ord_p \alpha\!<\ord_p \beta \Longrightarrow 
\ord_p(\alpha\pm \beta)\!=\!\ord_p \alpha$. We also recall that 
the metrics $|\cdot|$ and $|\cdot|_p$ are respectively called 
{\em Archimedean} and {\em non-Archimedean} because as $n\longrightarrow 
\infty$ we have $|n|\longrightarrow \infty$, while the sequence 
$|n|_p$ remains inside the bounded set $\{1,1/p,1/p^2,\ldots\}$.  

Let us also recall that a polynomial-time {\em Las Vegas randomized} algorithm 
is a polynomia-time algorithm that uses polynomially random bits in the 
input size, errs with probability at worst $1/2$, but correctly reports 
if it errs. Such an algorithm can be run $k$ times to boost the success 
probability to $O(1/2^k)$, and this 
type of randomization is standard in many number-theoretic algorithms,
such as the fastest current algorithms for factoring polynomials over finite
fields or primality checking (see, e.g., \cite{ku08,cheng}). In our 
setting, errors (for a Las Vegas speed-up) consist of reporting too few roots 
in $\Q_p$, but such errors can be detected and reported at no extra cost. 
 
\subsection{Newton Polygons and Newton Iteration: Archimedean and 
Non-Archimedean}  
\label{sub:newt} 
Definitive sources for $p$-adic arithmetic and analysis include 
\cite{serre,schikhof,Rob00}. 
We denote the standard $p$-adic valuation on $\C_p$ (normalized so that 
$\ord_p p\!=\!1$) by  $\ord_p : \C_p \longrightarrow \Q$. 
The {\em most significant ($p$-adic)} {\em digit} of 
$\sum^\infty_{j=s}a_j p^j\!\in\!\Q_p$ is $a_s$, assuming the 
$a_j\!\in\!\{0,\ldots,p-1\}$ and $a_s\!\neq\!0$.  

The notion of Newton polygon goes back to 17th century work of Newton 
on Puiseux series solutions to polynomial equations 
\cite[pp.\ 126--127]{newtonbook}. 
We will need variants 
of this notion over $\C_p$ and $\C$. (See, e.g., \cite{wei63} for the 
$p$-adic case and \cite{ostrowskiarch,aknr} for the complex case.) 
\begin{dfn} {\em \label{dfn:newt} Suppose 
$f(x)\!:=\!\sum_{i=1}^{t}c_ix^{a_i}\!\in\!\Z[x]$ with 
$c_i\!\neq 0$ for all $i$ and $a_1\!<\!\cdots\!<\!a_t$. We then define 
the {\em $p$-adic Newton polygon}, $\newt_p(f)$ (resp.\ {\em 
Archimedean Newton polygon,} $\anewt(f)$) to be the convex hull of the set of 
points $\{(a_i,\ord_pc_i)\; | \;  i\!\in\!\{1,\ldots,t\}\}$ (resp.\ the convex 
hull of $\{(a_i,-\log|c_i|)\; | \;  i\!\in\!\{1,\ldots,t\}\}$). We call an edge 
$E$ of a polygon in $\R^2$ {\em lower} if and only if $E$ has an inner 
normal with positive last coordinate. We also define the 
{\em horizontal length} of a line segment $E$ connecting $(r,s)$ and $(u,v)$ 
to be $\lambda(E)\!:=\!|u-r|$. \dia } 
\end{dfn}
\begin{ex} {\em \label{ex:newts} 
Following the notation of Theorem \ref{thm:tetra}, we set 
$h\!=\!3$ and illustrate $\newt_p\left(f_{5,p}\right)$ (for $p$ odd) and 
$\newt_\infty(f_{5,1/2})$ below: \\  
\mbox{}\hfill\scalebox{1}[1]{\epsfig{file=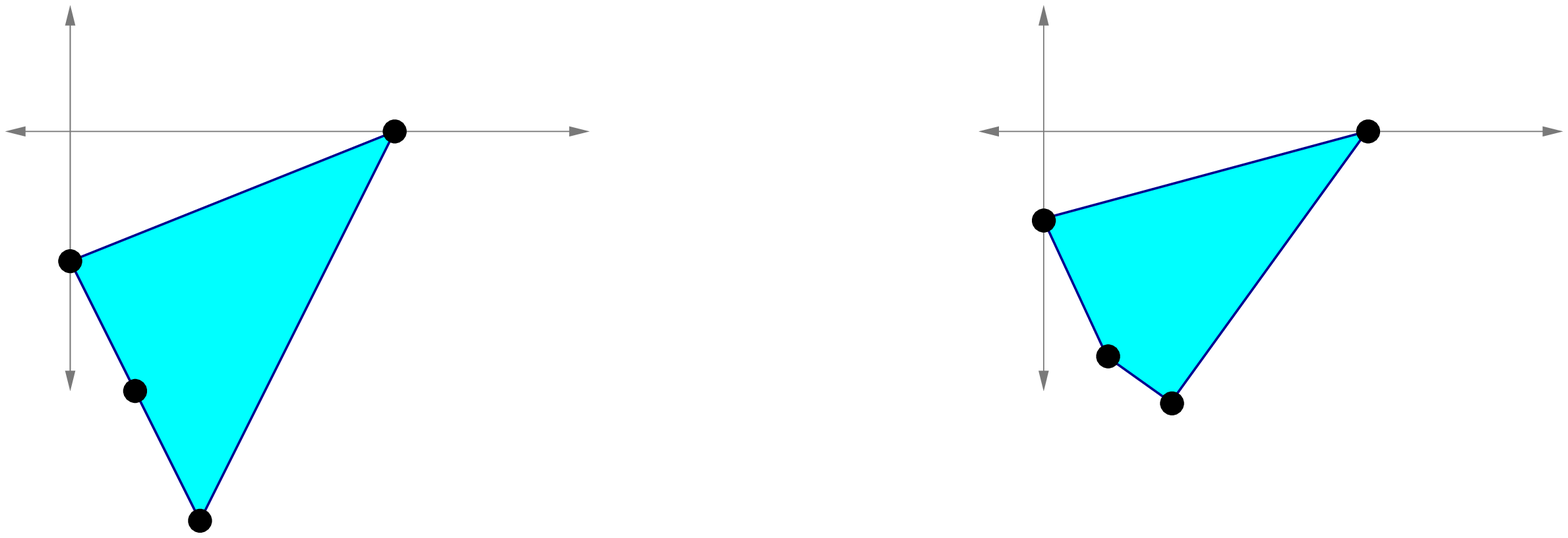,height=1.5in,
clip=}}\hfill\mbox{}\\ 
Note that the $p$-adic Newton polygon on the left has exactly $2$ lower 
edges (with horizontal lengths $2$ and $3$), while the Archimedean Newton 
polygon on the right has exactly $3$ lower edges (with horizontal lengths 
$1$, $1$, and $3$). \dia}  
\end{ex}  

\begin{thm} \label{thm:newt} {\em 
Following the notation above, the number\footnote{counting multiplicity} 
of roots of $f$ in $\C_p$ of valuation $v$ is exactly the horizontal length of 
the face of $\newt_p(f)$ with inner normal $(v,1)$. Furthermore, if $\anewt(f)$ 
has a lower edge $E$ with slope $v$, and no other lower edges with slope 
in the open interval $(v-\log 3,v+\log 3)$, then the number$^2$ of roots 
$\zeta\!\in\!\C$ of $f$ with $\log|\zeta|\!\in\!(v-\log 3,v+\log 3)$ 
is exactly $\lambda(E)$. \qed }  
\end{thm} 

\noindent 
The first portion of Theorem \ref{thm:newt} goes back to early 
20th century work of Hensel, while the second portion is an 
immediate consequence of \cite[Thm.\ 1.5]{aknr} (with an  
important precursor in \cite{ostrowskiarch}). The set of slopes 
of the lower edges of $\newt_p(f)$ (or of $\anewt(f)$) is an example 
of a {\em tropical variety} \cite{aknr}. 

We will also use the following version of Hensel's famous 
criterion for the rapid convergence of Newton's method over $\C_p$: 
\begin{hensel} 
{\em (See, e.g., \cite[Thm.\ 4.1 \& Inequality (5.7)]{conrad}.)  
Suppose $p$ is prime, $f\in \Z[x]$, $j\!\geq\!1$,   
$\zeta\!\in\!\Z_p$, $\ell\!=\!\ord_p f'(\zeta)\!<\!\infty$, and 
$f(\zeta)\equiv 0\mod p^{2\ell+j}$. 
Let $\zeta'\!:=\!\zeta-\frac{f(\zeta)}{f'(\zeta)}$. 
Then $f(\zeta')\!=\!0$ mod $p^{2\ell+2j}$, $\ord_p f'(\zeta')\!=\!\ell$, 
and $\zeta\!=\!\zeta'$ mod $p^{\ell+2j}$. \qed }  
\end{hensel}

\subsection{Separating Roots of Binomials} 
\label{sec:bisep}
When $f\!\in\!\Z[x]$ is a binomial, all of its roots in $\C$ are 
multiples of roots of unity that are evenly spaced on a circle. 
The same turns out to be true over $\C_p$, but the root spacing 
then depends more subtly on $p$ and much less on the degree. For convenience, 
we will sometimes write $|\cdot|_\infty$ instead of $|\cdot|$ 
for the standard norm on $\C$. Rather than stating  
lower bounds on $|\zeta_1-\zeta_2|_p$ (which always tend to $0$ as 
$H\longrightarrow\infty$ in our setting), we will instead state {\em upper} 
bounds on $|\log|\zeta_1-\zeta_2|_p|$: the latter clearly  
includes {\em both} a lower and upper bound on $|\zeta_1-\zeta_2|_p$. 
In summary, we have the following: 
\begin{prop}{\em \label{prop:bi} 
Suppose $f(x)\!:=\!c_1+c_2x^d\!\in\!\Z[x]$, $d\!\geq\!2$, $c_1c_2\!\neq\!0$, 
and $|c_1|,|c_2|\!\leq\!H$. Then for any distinct roots 
$\zeta_1,\zeta_2\!\in\!\C$ of $f$, we have $|\log|\zeta_1-\zeta_2||\!\leq\!
\log(d)+\frac{1}{d}\log H$. Also, for any distinct roots 
$\zeta_1,\zeta_2\!\in\!\C_p$ of $f$, we have that  
$|\log|\zeta_1-\zeta_2|_p|$ is at most  
$\frac{1}{d}\log H$ or $\frac{\log p}{p-1}+\frac{1}{d}\log H$,  
according as $d\!>\!p^{\ord_p d}$ or $d\!=\!p^{\ord_p d}\!\geq\!p$. }
\end{prop} 

\noindent
Put another way, if one fixes $p$ and $H$, and lets $d\longrightarrow \infty$, 
then the minimal root distance tends to $0$ at a rate of $\Theta(1/d)$ for the 
Archimedean case. However, in the non-Archimedean case, the minimal root 
distance {\em is never less than} $\frac{1}{Hp^{1/(p-1)}}$.  

\smallskip 
\noindent 
{\bf Proof of Proposition \ref{prop:bi}:} 
The case $p\!=\!\infty$ follows from an estimate for the
distance between the vertices of a regular $d$-gon. In particular,
the minimal spacing between distinct complex roots can easily be expressed
explicitly as $|c_1/c_2|^{1/d}\sqrt{2(1-\cos\frac{2\pi}{d})}$, which is
clearly bounded from below by $H^{-1/d} \sqrt{2(1-\cos\frac{2\pi}{d})}$.
From the elementary
inequality $1-\cos x\!\geq\!x^2\left(\frac{1}{2!}-\frac{\pi^2}{48}\right)$
we easily get $\left|\frac{1}{2}\log\left(1-\cos\frac{2\pi}{d}\right)
\right|\!\leq\!\log(d)-\frac{1}{2}
\log\left(4\pi^2-\frac{\pi^2}{6}\right)$ for all $d\!\geq\!6$. Observing that
$|\frac{1}{2}\log(1-\cos\frac{2\pi}{d})|\!\leq\!\log 2$ for
$d\!\in\!\{2,\ldots,5\}$ we get our stated bound via the Triangle Inequality
applied to
$\left|\log\left(H^{-1/d}\sqrt{2(1-\cos\frac{2\pi}{d})} \right)\right|$.

The case of prime $p$ follows easily from the Ultrametric Inequality and
classical facts on the spacing of $p$-adic roots of unity (see, e.g.,
\cite[Cor.\ 1, Pg.\ 105, Sec.\ 4.3 \& Thm.\ Pg.\ 107,
Sec.\ 4.4]{Rob00}). In particular, when $\gcd(d,p-1)$, the $d$th roots of
unity in $\C_p$ are all at unit distance.
At the opposite extreme of
$d\!=\!p^j$ for $j\!\geq\!1$, the set of distances between distinct $d$th roots
is exactly $\left\{\frac{1}{p^{1/(p-1)}}, \frac{1}{p^{1/
(p^1(p-1))}}, \ldots,\frac{1}{p^{1/(p^{j-1}(p-1))}}\right\}$. So the 
minimum distance is $1/p^{1/(p-1)}$ for $d$ a non-trivial $p$th power.
In complete generality, we see that there are distinct $d$th roots of unity at
distance $1$ if and only if $d$ is divisible by a prime other than $p$.
Observing that $\ord_p\!\left(H^{-1/d}\right)\!=\!-\frac{1}{d}\ord_p H\!
\geq\!-\frac{\log H}{d\log p}$ and $|x|_p\!=\!p^{-\ord_p x}$, we then see that
$\log|H^{-1/d}|_p\!\geq\!-\frac{1}{d}\log H$ and our bound follows again from
the Triangle Inequality. \qed

\subsection{Characterizing Roots of Binomials Over $\Q^*_p$} 
\label{sub:binocharqp}  
For any ring $R$ we let $R^*$ denote the multiplicatively invertible elements 
of $R$.  Counting roots of binomials over $\Q_p$ is more involved than 
counting their roots over $\R$, but is still quite efficiently doable. 
The first step is reducing the problem to $\Z/(p^k)$ for $k$ linear 
in the bit-size of the degree of the binomial. 
\begin{lem} {\em \label{lem:binoqp} 
Suppose $p$ is an odd prime and $f(x)\!:=\!c_1+c_2x^d\!\in\!\Z[x]$ with 
$|c_1|,|c_2|\!\leq\!H$, $c_1c_2\!\neq\!0$, and $\ell\!:=\!\ord_p d$. 
Then the number of roots of $f$ in $\Q_p$ is either $0$ or 
$\gcd(d,p-1)$. In particular, $f$ has roots in $\Q_p$ if and only if 
{\em both} of the following conditions hold:\\ 
\mbox{}\hfill 
(1) $d|\ord_p(c_1/c_2)$ \ and \ (2) $\left(-\frac{c_1}{c_2}p^{\ord_p(c_2/c_1)}
\right)^{p^{\ell}(p-1)/\gcd(d,p-1)}\!=\!1$ {\em mod} 
$p^{2\ell+1}$. \hfill \qed }  
\end{lem} 

\noindent 
Lemma \ref{lem:binoqp} is classical and follows from basic group theory 
(the fact that the multiplicative group $(\Z/(p^k))^*$ is cyclic, of order 
$p^{k-1}(p-1)$, for $p$ odd) and Hensel's Lemma. 

Recall that the only roots of unity in $\Q_2$ are $\{\pm 1\}$ (see, e.g.,
\cite{Rob00}). The following lemma is then a simple consequence of the
multiplicative group $(\Z/(2^k))^*$ being exactly the product
$\{\pm 1\}\times \left\{1,5,\ldots,5^{2^{k-3}} \text{ mod } 2^k 
\right\}$ (having cardinality $2^{k-1}$) when $k\!\geq\!3$
(see, e.g., \cite[Thm.\ 5.6.2, pg.\ 109]{bs}), and Hensel's Lemma.
\begin{lem} {\em \label{lem:binoq2}
Suppose $f(x)\!:=\!c_1+c_2x^d\!\in\!\Z[x]$ with
$|c_1|,|c_2|\!\leq\!H$, and $c_1c_2\!\neq\!0$. Then the number of roots of
the binomial $f$ in $\Q_2$ is either $0$ or $\gcd(d,2)$. In particular,
if $\ell\!:=\!\ord_2 d$ and $u\!:=\!\ord_2(c_2/c_1)$, then $f$ has roots in
$\Q_2$ if and only if {\em both} of the following conditions\linebreak 
\scalebox{.9}[1]{hold: (1) $d|u$ and (2) either (i) $d$ is odd or (ii) both
$\frac{c_1}{c_2}2^u\!=\!-1$ mod $8$ and $\left(-\frac{c_1}{c_2}2^u 
\right)^{2^{\ell-1}}\!\!\!\!\!\!=\!1$ mod $2^{2\ell+1}$. \qed}}
\end{lem}

\subsection{Bit Complexity Basics and Counting Roots of Binomials} 
\label{sub:binocountqp} The following bit-complexity estimates for finite ring 
arithmetic will be fundamental for our main algorithmic results, and follow 
directly from the development of \cite[Ch.\ 4 \& 11]{vzgbook} 
(particularly \cite[Cor.\ 11.13, pg.\ 327]{vzgbook}) assuming 
one uses the recent fast integer multiplication algorithm of Harvey and 
van der Hoeven \cite{harveymult}. See also 
\cite{sutherland} for an excellent exposition on most of the bounds below.   
We use $\log^* x$ to denote the minimal $k$ such that $k$ compositions of 
$\log$ applied to $x$ yield a real number $\leq\!1$.  
\begin{thm} \label{thm:cxity} {\em For any prime $p\!\in\!\N$ and 
$j,m,n\!\in\!\N$, we have the following bit-complexity bounds (in the Turing 
model) involving $A,a,b,c\!\in\!\N$ with $A,a,b\!\leq\!2^n-1$, 
$A\!\geq\!2^{n-1}$, $c\!\leq\!2^m-1$ with $m\!=\!O(\log n)$, 
$r,s\!\in\!\{0,\ldots,p^j-1\}$ with $p\nmid r$, and $f,g\!\in\!\F_p[x]$ 
both having degree $\leq\!d$: \\ 
\mbox{}\hfill 
\begin{tabular}{r|l}   
Operation          & Best Current $O$-bound (as of July 2021) \\ \hline  
$a+b$              & $O(n)$ \\ 
$a\cdot b$         & $O(n\log n)$ \\ 
$a$ mod $b$        & $O(n\log n)$ \\ 
$A$ mod $c$        & $O(nm)$ \\ 
$r\cdot s$ mod $p^j$ & $O(j\log(p)\log(j\log p))$ \\ 
$1/r$ mod $p^j$    & $O(j\log(p)\log^2(j\log p))$ \\ 
$r^s$ mod $p^j$    & $O(j^2\log^2(p)\log(j\log p))$ \\ 
$f\cdot g$         & $O\!\left(d\log(p)\log(d\log(p))4^{\log^*(d\log p)}
                      \right)$ \\ 
$\gcd(f,g)$        & $O(d\log(p)\log^2(d)\log(\log d)\log \log p)$ \\ 
\end{tabular} 
\hfill\qed } 
\end{thm} 

\noindent 
We note that the penultimate bound comes directly from \cite{harveymultFq}. 
The very last bound is actually a simple combination of the Half-gcd 
algorithm from \cite[Thm.\ 11.1, Ch.\ 11]{vzgbook} with the fast 
polynomial multiplication algorithm from \cite{cantorkaltofen}, and 
can likely be improved slightly via the techniques of \cite{harveymultFq}.  
\begin{cor} {\em \label{cor:binomod} 
Following the notation of Lemmata \ref{lem:binoqp} and \ref{lem:binoq2}, 
one can count exactly the number of roots of $f$ in $\Q_p$ in 
time $O\!\left(\log^2(dpH)\log\log(dpH)\right)$.  
Furthermore, for any root $\zeta\!\in\!\Q^*_p$ there is an 
$x_0\!\in\!\Z\left/\left( p^{2\ell+1}\right) 
\right.$ that is a root of the mod $p^{2\ell+1}$ reduction of\linebreak  
$\frac{c_1}{p^{\ord_pc_1}}+ \frac{c_2}{p^{\ord_pc_2}}x^d$,
and with $z_0\!:=\!p^{\ord_p(c_2/c_1)/d}x_0\!\in\!\Q$  
an approximate root of $f$ with associated true root $\zeta$. In 
particular, the logarithmic height\footnote{The logarithmic height of 
a rational number $a/b$ with $\gcd(a,b)\!=\!1$ is simply 
$\log\max\{|a|,|b|\}$ (and 
we declare the logarithmic height of $0$ to be $0$).} 
of $z_0$ is $O\!\left(\log\left(pH^{1/d} \right)\right)$. }  
\end{cor}  

\noindent 
{\bf Proof: (Case of odd $p$)}  
First note that $(\Z/p^{2\ell+1})^*$ is cyclic and
Lemma \ref{lem:binoqp} tells us that we can reduce deciding the feasibility of
$c_1+c_2x^d\!=\!0$ over $\Q^*_p$ to checking $d\stackrel{?}{|}\ord_p(c_1/c_2)$
and $(-c_1/c_2)^r\!\stackrel{?}{=}\!1$ mod $p^{2\ell+1}$ with $r\!=\!
p^\ell(p-1)/\gcd(d,p-1)$.

The $p$-adic valuation can be computed easily by bisection, ultimately
resulting in $O(\log H)$ divisions involving integers with
$O(\max\{\log p,\log H\})\!=\!O(\log(pH))$ bits. Checking divisibility by $d$ 
involves dividing an integer with $O(\log\log H)$ bits by an integer with 
$O(\log d)$ bits. By Theorem \ref{thm:cxity} these initial steps take 
time $O(\log(H)\log(pH)\log\log(pH))+O(m\log m)$, where 
$m\!=\!\max\{\log\log H,\log d\}$. 
By Theorem \ref{thm:cxity}, the $r$th power can be computed 
in time $O(\ell^2\log^2(p)\log(\ell \log p))$. So our 
overall complexity bound is \\ 
\mbox{}\hfill $O(\ell^2\log^2(p)\log(\ell \log p)+\log(H)\log(pH)\log\log(pH)
 +\log(d)\log\log d)$.\hfill\mbox{}\\  
Since $\ell\!\!\leq\!\log_p d$ our final bound becomes \\ 
\mbox{}\hfill $O(\log^2(d)\log(\log d)+\log(H)\log(pH)\log 
\log(pH))$.\hfill\mbox{}\\  
A simple over-estimate then yields our stated complexity bound.
The remainder of the lemma then follows easily from
Hensel's Lemma and Proposition \ref{prop:bi}. \qed 

\smallskip 
\noindent 
{\bf (Case of $p\!=\!2$)} 
The proof is almost identical to the odd $p$ case, 
save that we use Lemma \ref{lem:binoq2} in place of Lemma \ref{lem:binoqp}.
In particular, the case $\ell\!=\!0$ remains unchanged.

As for the case $\ell\!\geq\!1$, the only change is an extra congruence
condition (mod $8$) to check whether $\frac{c_1}{c_2}2^u$ is a square 
mod $2^{2\ell+1}$ (see, e.g., \cite[Ex.\ 38, pg.\ 192]{bs}). However, this 
additional complexity is negligible compared to the other steps, so we are 
done. \qed

\subsection{Trees and Roots in $\Z/(p^k)$ and $\Z_p$}  \label{sub:trees} 
Recall that for any field $K$, a root $\zeta\!\in\!K$ of $f$ is 
{\em degenerate} if and only if $f(\zeta)\!=\!f'(\zeta)\!=\!0$. 
The $p$-adic analogue of bisecting an isolating interval containing a real 
root is to approximate the next base-$p$ digit of an approximate root in 
$\Q_p$. Shifting from bisecting intervals to extracting 
digits is crucial since $\Q_p$ is not an\linebreak 
\scalebox{.96}[1]{ordered field. We will write $f'$ for 
the derivative of $f$ and $f^{(i)}$ for the $i$th order derivative of $f$.}  
\begin{dfn} {\em \label{dfn:crazytree} \cite{krrz19} 
For any $f\in \Z[x]$ let $\tf$ denote the mod $p$ reduction of $f$. Assume 
$\tf$ is not identically $0$. Then, for any 
degenerate root $\zeta_0\!\in\!\{0,\ldots,p-1\}$ of $\tilde{f}$, 
we then define $s(f,\zeta_0):=\min_{i\geq 0}\left\{ i
+\ord_p\frac{f^{(i)}(\zeta_0)}{i!}\right\}$.   
Fixing $k\in \N$, for $i\geq 1$, let us inductively define a set 
$\cT_{p,k}(f)$ of pairs $(f_{i-1,\mu}, k_{i-1,\mu})$   
$\in \Z[x]\times \N$: We set $(f_{0,0}, k_{0,0}) := (f,k)$. Then for any 
$i\geq 1$ with $(f_{i-1,\mu}, k_{i-1,\mu})\!\in\!
\cT_{p,k}(f)$, and any degenerate root $\zeta_{i-1}\!\in\!\F_p$   
of $\tilde{f}_{i-1,\mu}$ with
$s_{i-1}:= s(f_{i-1,\mu},\zeta_{i-1})\in \{2,\ldots,k_{i-1,\mu}-1\}$, 
we define $\zeta:= \mu + \zeta_{i-1}p^{i-1},k_{i,\zeta}:=k_{i-1,\mu} 
- s_{i-1}$,\linebreak 
\scalebox{.93}[1]{$f_{i,\zeta}(x) := p^{-s(f_{i-1,\mu},\zeta_{i-1})} 
f_{i-1,\mu}(\zeta_{i-1} + px)  \mod p^{k_{i,\zeta}}$, and 
then include append $(f_{i,\zeta},k_{i,\zeta})$ to $\cT_{p,k}(f)$. \dia}}  
\end{dfn}

\begin{ex} \label{ex:tri} 
{\em If $f(x)\!=\!x^{10}-10x+738$ and $p\!=\!3$ then $\tf(x)\!=\!x(x-1)^9$ mod 
$3$, $1$ is a degenerate root of $\tf$ in $\F_3$, and one can check that 
$s(f,1)\!=\!4$ (no greater than the multiplicity of the factor $x-1$ in 
$\tf$). In particular, $f_{1,1}$ has degree $10$ (and $10$ monomial terms) but 
$\tf_{1,1}(x)\!=\!x^3+2x^2$. \dia } 
\end{ex} 

The collection of pairs $(f_{i,\zeta},k_{i,\zeta})$ admits a tree structure 
that will give us a way to extend Hensel lifting to degenerate roots. 
\begin{dfn}
\label{dfn:tree} {\em \cite{krrz19} 
The set $\cT_{p,k}(f)$ naturally admits the structure of a labelled, 
rooted, directed tree as follows\footnote{This 
definition differs slightly from the original in \cite{krrz19}: the edges 
are unlabelled here.}  
\begin{itemize}
\item[(i)]{We set $f_{0,0}\!:=\!f$, $k_{0,0}\!:=\!k$, and let
$(f_{0,0},k_{0,0})$ be the label of the root node of
$\cT_{p,k}(f)$.}  
\item[(ii)]{\scalebox{.94}[1]{The non-root nodes of $\cT_{p,k}(f)$ are 
labelled by the $(f_{i,\zeta},k_{i,\zeta})\!\in\!\cT_{p,k}(f)$ 
with $i\!\geq\!1$.}}
\item[(iii)]{There is an edge from node $(f_{i-1,\mu},k_{i-1,\mu})$ to
node $(f_{i,\zeta},k_{i,\zeta})$ if and only if 
there is a degenerate root $\zeta_{i-1}\!\in\!\F_p$ of $\tf_{i-1,\mu}$ with 
$s(f_{i-1,\mu},\zeta_{i-1})
\!\in\!\{2,\ldots,k_{i-1,\mu}-1\}$ and 
$\zeta\!=\!\mu+\zeta_{i-1}p^{i-1}\!\in\!\Z/(p^i)$. \dia}
\end{itemize} } 
\end{dfn}

\noindent 
We call each $f_{i,\zeta}$ with $(f_{i,\zeta},k_{i,\zeta})\!\in\!\cT_{p,k}(f)$ 
a {\em nodal polynomial} of $\cT_{p,k}(f)$. 
It is in fact possible to list all the roots of $f$ in $\Z/(p^k)$ from the 
data contained in $\cT_{p,k}(f)$ \cite{krrz19,DMS19}. We 
will instead use $\cT_{p,k}(f)$, with $k$ determined by root separation/
valuation condition, to efficiently {\em count} the roots of $f$ in $\Z_p$, 
and then in $\Q_p$ by rescaling.  

\begin{ex} 
{\em $\cT_{p,k}(x^2)$ is a chain of length $\floor{\frac{k-1}{2}}$ for any 
$p,k$. \dia } 
\end{ex} 
\begin{ex} 
{\em Let $f(x)\!=\!1-x^{397}$. Then $\cT_{17,k}(f)$, for any $k\!\geq\!1$,  
consists of a single node, labelled $(1-x^{397},k)$, since $\tilde{f}$ has 
no degenerate roots in $\F_{17}$. In particular, $f$ has $1$ as 
its only root in $\Q_{17}$. 
\dia}   
\end{ex}            
\begin{ex}         
{\em Let $f(x)\!=\!1-x^{340}$. Then, when 
$k\!\in\!\{1,2\}$, the tree $\cT_{17,k}(f)$ consists 
of a single root node, labelled $(1-x^{340},k)$. 
However, when $k\!\geq\!3$, the tree $\cT_{17,k}(f)$ has depth $1$, 
and consists of the aforementioned root node {\em and} exactly $4$ child 
nodes, labelled $(f_{1,\zeta_0},k-2)$ where the 
$\tf_{1,\zeta_0}$ are, respectively, $14x$, $12x+10$, 
$5x+15$, and $3x+3$. Note that $\tilde{f}$ has 
exactly $4$ roots $\zeta_0\!\in\!\F_{17}$ ($1$, $4$, $13$, and $16$), 
each of which is degenerate, and the roots $\zeta_1\!\in\!\F_{17}$ of the 
$\tf_{1,\zeta_0}$ encode the ``next'' base-$17$ digits ($0$, $2$, 
$14$, and $16$) of the roots of $f$ 
in $\Z/(17^2)$. In particular, the roots of $f$ in $\Q_{17}$ are 
$1+0\cdot 17+\cdots$, $4+2\cdot 17+\cdots$,
$13+14\cdot 17+\cdots$, 
and $16+16\cdot 17+\cdots$ and are all {\em non}-degenerate. 
\dia}  
\end{ex} 

Nodal polynomials --- originally defined for efficient root counting over 
$\Z/(p^k)$ --- thus encode individual base-$p$ digits of roots of 
$f$ in $\Z_p$. Their degree also decays in a manner depending on root 
multiplicity.\footnote{Over any field $K$, we define the {\em multiplicity  
of a root $\zeta\!\in\!K$} of $f\!\in\!K[x]$ as the greatest $m$ with  
$(x-\zeta)^m|f$.} 
\begin{lem} {\em \label{lem:nodal} 
\cite[Lem.\ 2.2 \& 3.6]{krrz19} Suppose $f\!\in\!\Z[x]\setminus p\Z[x]$ 
has degree $d$, $f_{0,0}\!:=\!f$, $i\!\geq\!1$, 
$\mu\!:=\!\zeta_0+\cdots+p^{i-2}\zeta_{i-2}$ is a root of the mod $p^{i-1}$ 
reduction of $f$, $\zeta'\!:=\!\mu+p^{i-1}\zeta_{i-1}$, 
the pairs $(f_{i-1,\mu},k_{i-1,\mu})$ and $(f_{i,\zeta'},k_{i,\zeta'})$ both 
lie in $\cT_{p,k}(f)$, and $\zeta_{i-1}$ has multiplicity $m$ as a root of 
$\tf_{i-1,\mu}$ in $\F_p$. Then $\cT_{p,k}(f)$ has depth 
$\leq\!\floor{(k-1)/2}$ and at most $\floor{d/2}$ nodes at depth 
$i\!\geq\!1$.  Also, 
$\deg \tf_{i,\zeta'}\!\leq\!s(f_{i-1,\mu},\zeta_{i-1})\!\leq\!\min\{k_{i-1,
\mu}-1,m\}$, and 
$f_{i,\zeta'}(x)\!=\!p^{-s}f(\zeta_0+\zeta_1p+\cdots +\zeta_{i-1}p^{i-1}+p^ix)$ 
where $s\!:=\!\sum^{i-1}_{j=0}
s(f_{j,\zeta_0+\cdots+\zeta_{j-1} p^{j-1}},\zeta_j)\!\geq\!2i$.
In particular, 
$f(\zeta_0+\zeta_1p+\cdots +\zeta_{i-1}p^{i-1})\!=\!0$ mod $p^s$ and  
$f'(\zeta_0+\zeta_1p+\cdots +\zeta_{i-1}p^{i-1})\!=\!0$ mod $p^i$. 
\qed}
\end{lem} 

\noindent 
Note that the first assertion of Lemma \ref{lem:nodal} gives us an upper bound 
on the depth of $\cT_{p,k}(f)$ as a function of $k$. We will also need to 
consider lower bounds on $k$ that guarantee that $\cT_{p,k}(f)$ has enough 
depth to be useful for approximating roots in $\Z_p$.  

Let $n_p(f)$ denote the number of non-degenerate roots in 
$\F_p$ of the mod $p$ reduction of $f$. 
\begin{lem} {\em \label{lem:ulift} 
Suppose $f\!\in\!\Z[x]$, $\zeta\!=\!\sum^\infty_{j=0}\zeta_jp^j\!\in\!\Z_p$ 
is a non-degenerate root of $f$, and  
let $D$ be the maximum of $\ord_p(\zeta-\xi)$ over
all distinct non-degenerate roots $\zeta,\xi\!\in\!\Z_p$ of $f$
(if $f$ has at least $2$ non-degenerate roots in $\Z_p$) or $0$
(if $f$ has $1$ or fewer non-degenerate roots in $\Z_p$). 
Then for all $k$ sufficiently large, $\cT_{p,k}(f)$ has a nodal polynomial 
$f_{j,\zeta'}$ such that $j\!\leq\!\floor{(k-1)/2}$ and 
$\zeta'+p^j\zeta_j\!=\!\zeta$ mod $p^{j+1}$ for some {\em non}-degenerate root 
$\zeta_j$ of $\tf_{j,\zeta'}$. 
Furthermore, for $k$ sufficiently large we also have that 
$\cT_{p,k}(f)$ has depth 
$\geq\!D$, the set $\{(g,j)\!\in\!\cT_{p,k}(f)\; | \; n_p(g)\!>\!0\}$ 
remains fixed and finite, and $f$ has exactly 
$\sum\limits_{(g,j)\in \cT_{p,k}(f)} n_p(g)$ non-degenerate roots in $\Z_p$.   
}  
\end{lem}

\noindent 
{\bf Proof:} First note that  
$f(\zeta_0+\cdots+\zeta_ip^i)\!=\!0$ mod $p^{i+1}$ for all $i\!\geq\!0$. 
By Definitions \ref{dfn:crazytree} and \ref{dfn:tree}, 
$s_0\!:=\!s(f,\zeta_0)\!\in\!\{1,\ldots,m\}$, 
where $m$ is the multiplicity of $\zeta_0$ 
as a root of $\tf$ (thanks to Lemma \ref{lem:nodal}).  
Should $m\!=\!1$ then $s_0\!=\!1$, leaving $f_{0,0}\!=\!f$ as our desired nodal 
polynomial (with $\zeta_0$ a non-degenerate root of $\tf_{0,0}$) 
for all $k\!\geq\!1$. Otherwise, $s_0\!\geq\!2$ (by the definition of 
$s(\cdot,\cdot)$), in which case $k\!\geq\!1+s_0 \Longrightarrow 
\cT_{p,k}(f)$ will have $f_{1,\zeta_0}(x)\!=\!p^{-s_0}f(\zeta_0+px)$ as a 
nodal polynomial. However, we need to check if $\zeta_1$ 
is a non-degenerate root for $\tf_{1,\zeta_0}$ or not. 

Proceeding inductively, note that if $i\!\geq\!1$, 
$\zeta'\!:=\!\zeta_0+\zeta_1p +\cdots+\zeta_{i-1}p^{i-1}$, 
$s_i\!:=\!s(f_{i,\zeta'},\zeta_i)$, and $s'\!:=\!s_0+\cdots+s_i$, then 
$s_i\!\in\!\{1,\ldots,m\}$ where $m$ is now the multiplicity of $\zeta_i$ 
as a root of $\tf_{i,\zeta'}$. 
As before, $m\!=\!1$ implies that $f_{i,\zeta'}$ is our desired nodal 
polynomial  (with $\zeta_i$ a non-degenerate root of $\tf_{i,\zeta'}$) 
for all $k\!\geq\!1+s'$. Otherwise, $s_i\!\geq\!2$, in which case 
$k\!\geq\!1+s'\Longrightarrow  
\cT_{p,k}(f)$ will have $f_{i+1,\zeta'+p^i\zeta_i}(x)\!=\!p^{-s'}
f(\zeta'+p^i\zeta_i+p^{i+1}x)$ as a nodal polynomial, and then we check 
if $\zeta_{i+1}$ is a non-degenerate root for $\tf_{i+1,\zeta'+p^i\zeta_i}$ 
or not. 

Our induction must end, in finitely many steps, with our desired 
$f_{j,\zeta'}$. To see why, first observe that nodal polynomials always have 
integer coefficients and, if $d'\!:=\!\ord_pf'(\zeta)$, then  
$d'\!<\!\infty$ since $\zeta$ is a non-degenerate root and 
thus $f'(\zeta)\!=\!\alpha p^{d'}$ mod $p^{d'+1}$ for some 
$\alpha\!\in\!\Z_p\setminus p\Z_p$. So if our induction reaches $f_{i,\zeta'}$ 
with $i\!\geq\!d'$, then $\zeta'\!=\!\zeta_0+\cdots+p^{d'-1}\zeta_{d'-1}
\Longrightarrow f'_{d',\zeta'}(\zeta_{d'})\!=\!\alpha p^{2d'
-(s_0+\cdots+s_{d'-1})}$. We thus obtain $2d'\!\geq\!s_0+\cdots+s_{d'-1}$ 
and, for all $i\!\in\!\N$ {\em with $f_{i,\zeta'}$ belonging to a node of 
$\cT_{p,k}(f)$ with a child}, the definition of $s_i$ tells us that 
$s_i\!\geq\!2$. Since $\ord_p f'(\zeta_0+\cdots+p^i\zeta_i)\!=\!d'$ for all 
$i\!\geq\!d'$, we must eventually encounter a $j\!\geq\!d'$ with 
$s_j\!=\!1$, meaning no child for $f_{j,\zeta'}$. So our induction ends 
with a nodal polynomial $f_{j,\zeta'}$ with no degenerate roots. Moreover, 
we must have $\tf_{j,\zeta'}(\zeta_j)\!=\!0$ mod $p$ (by definition of 
$\zeta$ and $f_{j,\zeta'}$) and thus 
$\zeta_j$ must be a non-degenerate root of $\tf_{j,\zeta'}$. 
Also, our upper bound on $j$ is immediate from Lemma \ref{lem:nodal}. 

To prove that $\cT_{p,k}(f)$ has depth for $D$ for $k$ large enough, 
note that an $f$ with {\em no} non-degenerate 
roots in $\Z_p$ can {\em not} have a tree $\cT_{p,k}(f)$ having nodal 
polynomials with non-degenerate roots in $\F_p$. This is because of the 
equality $f_{i,\zeta'}(x)\!=\!p^{-s}f(\zeta_0+\zeta_1p + \cdots + 
\zeta_{i-1}p^{i-1}+p^i x)$ from Lemma \ref{lem:nodal}: $\tf_{i,\zeta'}$ 
having a non-degenerate root in $\F_p$ would imply by Hensel's Lemma 
that $f$ has a root $\zeta\!\in\!\Z_p$ with $\ord_p f'(\zeta)\!<\!\infty$. 
So in this case, the stated set of $(g,j)$ is empty for all $k\!\geq\!1$ 
and the stated sum is $0$. In particular, $\cT_{p,k}(f)$ always at 
least has its root node (by definition) and thus $D\!\geq\!0$.  

Similarly, an $f$ with just one non-degenerate root in $\Z_p$ can {\em not} 
have a tree $\cT_{p,k}(f)$ having two distinct nodal polynomials having 
non-degenerate roots mod $p$. (Likewise, $\cT_{p,k}(f)$ having a single 
nodal polynomial with two distinct non-degenerate roots mod $p$ is 
impossible.) So in this case, the stated set of $(g,j)$ 
has cardinality $1$ (with $n_p(g)\!=\!1$ for exactly one pair 
$(g,j)$) for all $k$ as specified in the first assertion of 
our lemma, which we've already proved. So the remaining assertions follow.  

So let us now assume $f$ has at least $2$ distinct non-degenerate roots in 
$\Z_p$. There are clearly no more than $\deg f$ such roots, so our 
first assertion implies that, for $k$ sufficiently large, {\em every}   
non-degenerate root $\zeta\!\in\!\Z_p$ of $f$ has an   
associated node in $\cT_{p,k}(f)$ encoding $\zeta$, i.e., 
$\cT_{p,k}(f)$ has depth at least $D$ for $k$ sufficiently large. 
Clearly then, the set $\{(g,j)\!\in\!\cT_{p,k}(f)\; | \; n_p(g)\!>\!0\}$ is 
finite and will not change as $k$ increases: This is because the set can not 
lose elements as $k$ increases, and any new element would introduce a 
new non-degenerate root for $f$ via Hensel's Lemma. 

So we now only need to prove that the stated sum counts roots correctly. 
Toward this end, note by construction that every non-degenerate root 
$\zeta\!\in\!\Z_p$ of $f$ is associated to a unique sequence 
of the form $(\zeta_0,\ldots, \zeta_i)\!\in\!\F^{i+1}$ with 
$\zeta_0,\ldots,\zeta_{i-1}$ all degenerate 
roots for previously defined nodal polynomials, but 
with $\zeta_i$ a {\em non}-degenerate root of $\tf_{i,\zeta'}$. 
So the number of non-degenerate roots of $f$ in $\Z_p$ is no greater than the 
stated sum. 
 
To conclude, note that Hensel's Lemma (and our 
earlier observation that nodal polynomials
are rescaled shifts of $f$) implies that each 
non-degenerate root in $\F_p$ of a nodal polynomial lifts to  
a unique root of $f$ in $\Z_p$. Furthermore, since 
the derivatives of nodal polynomials are rescaled shifts of $f'$, 
each such lifted root is a non-degenerate root. So the number of 
non-degenerate roots of $f$ in $\Z_p$ is at least as large as the stated sum, 
and we are done. \qed   

\subsection{Trees and Extracting Digits of Radicals} 
We prove the following useful lemma in Remark \ref{rem:low} of Section 
\ref{sec:trinosolqp}:  
\begin{lem} {\em \label{lem:binodepth} 
Suppose $f(x)\!=\!c_1+c_2x^d\!\in\!\Z[x]$ with $c_1c_2\!\neq\!0$ mod $p$ 
and $\ell\!:=\!\ord_p d$. Then every {\em non}-root nodal polynomial 
$f_{i,\zeta}$ of $\cT_{p,k}(f)$ satisfies $\deg \tf_{i,\zeta}\!\leq\!2$ or 
$\deg \tf_{i,\zeta}\!\leq\!1$, according as $p\!=\!2$ or $p\!\geq\!3$. 
In particular, $f(\zeta_0)\!=\!0$ mod $p$ for some 
$\zeta_0\!\in\!\{0,\ldots,p-1\} \Longrightarrow s(f,\zeta_0)\!\leq\!\ell+1$. 
\qed } 
\end{lem} 

\begin{rem} \label{rem:depth} 
{\em It is a simple exercise to prove, from Lemma \ref{lem:binodepth} and 
Definition \ref{dfn:tree}, that $\cT_{p,k}(f)$ always has depth $\leq\!1$ 
for $f\!\in\!\Z[x]$ a binomial with $f(0)\!\neq\!0$ 
mod $p$. The family of examples $x^{p^2}-1$  
(for any $k\!\geq\!4$) shows that this depth can be attained for 
any prime $p$. \dia } 
\end{rem} 

With our tree-based encoding of $p$-adic roots in place, we can now 
prove that it is easy to find approximate roots in $\Q_p$ for binomials 
when $p$ is fixed.  
\begin{thm} {\em
\label{thm:binoqp}
Suppose $f\!\in\!\Z[x]$ is a binomial
of degree $d$ with coefficients of absolute value at most $H$, $f(0)\!\neq\!0$, 
$\gamma\!=\!\gcd(d,\max\{2,p-1\})$, and $\{\zeta_1,\ldots,\zeta_\gamma\}$ is 
the set of roots of $f$ in $\Q_p$. Then in time 
$O\!\left(\left(\frac{p}{\gamma}+\gamma+\log d\right)
\log(dp)\log\log(dp)+\log(H)\log(pH)\log\log(pH)\right)$,
we can find, for each $j\!\in\!\{1,\ldots,\gamma\}$, a $z^{(j)}_0\!\in\!\Q$ of 
logarithmic height $O\!\left(\log\left(dH^{1/d}\right)\right)$ 
that is an approximate root with associated true root $\zeta_j$. } 
\end{thm}

An algorithm that proves Theorem \ref{thm:binoqp} when $p$ is odd is outlined 
below.\\  
\mbox{}\scalebox{.96}[.8]{\fbox{\mbox{}\hspace{.1cm}\vbox{
\begin{algor} {\em 
\label{algor:binoqp}
{\bf (Solving Binomial Equations Over $\pmb{\Q^*_p}$ for odd $p$)}  
\mbox{}\\
{\bf Input.} An odd prime $p$ and 
$c_1,c_2,d\!\in\!\Z\setminus\{0\}$ with $|c_i|\!\leq\!H$ for all $i$. \\ 
{\bf Output.} A true declaration that $f(x)\!:=\!c_1+c_2x^d$ has 
no roots in $\Q_p$, or $z_1,$ 
$\ldots,z_\gamma\!\in\!\Q$ with\\ 
\mbox{}\hspace{1.8cm}\scalebox{.98}[1]{logarithmic 
height $O\!\left(\log\left(dH^{1/d}\right)\right)$ such that $\gamma\!=\!
\gcd(d,p-1)$, $z_j$ is an approximate}\\ 
\mbox{}\hspace{1.8cm}\scalebox{.96}[1]{root with associated true root 
$\zeta_j\!\in\!\Q_p$ for all $j$, and the $\zeta_j$ are pair-wise distinct.}\\
{\bf Description.} \\ 
1: \scalebox{.95}[1]{If $\ord_p c_1\!\neq\!\ord_p c_2$ mod $d$ then say 
{\tt ``No roots in $\Q_p$!''} and {\tt STOP}.} \\
2: Let $\ell\!:=\!\ord_p d$ and replace $f$ with $f(x)\!:=\!c'_1+c'_2x^d$ 
where $c'_i\!:=\!\frac{c_i}{p^{\ord_p c_i}}$ for all $i$. \\
3: \scalebox{.93}[1]{If $\left(-\frac{c'_1}{c'_2}
\right)^{p^{\ell}(p-1)/\gamma}\!\!\!\!\!\!\!\!\!\neq\!1$ 
mod $p^{2\ell+1}$ then say {\tt ``No roots in $\Q_p$!''} and {\tt STOP}.}\\ 
4: \scalebox{.96}[1]{Let $\delta\!:=\!1$. If $d\!\leq\!-1$ then 
set $\delta\!:=\!-1$ and respectively replace $d$ by $|d|$ and $f(x)$ by 
$x^df(1/x)$.}\\ 
5: \scalebox{.895}[1]{Let $g$ be any generator for $\F^*_p$, 
$r\!:=\!(d/\gamma)^{-1}$ mod $p-1$, 
$c'\!:=\!(-c'_1/c'_2)^r$ mod $p$, and $\tilde{h}(x)\!:=\!x^{\gamma}-c'$.}\\ 
6: Find a root $x_1\!\in\!\left\{g^0,\ldots,g^{\frac{p-1}{\gamma}-1}\right\}$ 
of $\tilde{h}$ via brute-force search.\\ 
7: For all $j\!\in\!\{2,\ldots,\gamma\}$ let 
$x_j\!:=\!x_{j-1}g^{(p-1)/\gamma}$ mod $p$. \\ 
8: \scalebox{.87}[1]{If $\ell\!\geq\!1$ then, for each $j\!\in\!\{1,\ldots,
\gamma\}$, replace $x_j$ by $x_j -\frac{f(x_j)/p^\ell}
{f'(x_j)/p^\ell}\!\in\!\Z/(p^2)$.}\\ 
9: {\tt Output} $\left\{(x_1 p^{\ord_p(c_1/c_2)/d})^\delta,\ldots,(x_\gamma 
p^{\ord_p(c_1/c_2)/d})^\delta\right\}$. } 
\end{algor}} 
}}
\begin{rem} {\em
Step 6 above is designed for simplicity rather than practicality, and can be 
sped up considerably if one one avails to more sophisticated algorithms with 
complexity linear in $\gcd(d,p-1)$ and quasi-linear in $\log(pd)$: See, e.g., 
\cite{amm,caofan,chokwonlee}. \dia}
\end{rem}  

The following algorithm proves the $p\!=\!2$ case of Theorem 
\ref{thm:binoqp}.\\
\noindent
\mbox{}\scalebox{.96}[.8]{\fbox{\mbox{}\hspace{.1cm}\vbox{
\begin{algor} {\em
\label{algor:binoq2}
{\bf (Solving Binomial Equations Over $\pmb{\Q^*_2}$)}
\mbox{}\\
{\bf Input.}\hspace{.2cm}$c_1,c_2,d\!\in\!\Z\setminus\{0\}$ with
$|c_i|\!\leq\!H$ for all $i$. \\
{\bf Output.} A true declaration that $f(x)\!:=\!c_1+c_2x^d$ has
no roots in $\Q_2$, or $z_1,\ldots,z_\gamma\!\in\!\Q$ with\\
\mbox{}\hspace{1.8cm}\scalebox{.98}[1]{logarithmic
height $O\!\left(\log\left(dH^{1/d}\right)\right)$ such that $\gamma\!=\!
\gcd(d,2)$, $z_j$ is an approximate}\\
\mbox{}\hspace{1.8cm}\scalebox{.96}[1]{root of $f$ with associated true root
$\zeta_j\!\in\!\Q_p$ for all $j$, and the $\zeta_j$ are pair-wise distinct.}\\
{\bf Description.} \\
1: \scalebox{.95}[1]{If $\ord_2 c_1\!\neq\!\ord_2 c_2$ mod $d$ then say
{\tt ``No roots in $\Q_p$!''} and {\tt STOP}.} \\
2: Let $\ell\!:=\!\ord_2 d$ and replace $f$ with $f(x)\!:=\!c'_1+c'_2x^d$ where
$c'_i\!:=\!\frac{c_i}{2^{\ord_2 c_i}}$ for all $i$. \\
3: \scalebox{.93}[1]{If $c'_1\!\neq\!-c'_2$ mod $8$ or $\left(-\frac{c'_1}{c'_2}
\right)^{2^{\ell-1}}\!\!\!\!\!\!\!\!\!\neq\!1$
mod $2^{2\ell+1}$ then say {\tt ``No roots in $\Q_2$!''} and {\tt STOP}.}\\
4: \scalebox{.95}[1]{Let $\delta\!:=\!1$. If $d\!\leq\!-1$ then
set $\delta\!:=\!-1$ and respectively replace $d$ by $|d|$ and $f(x)$ by
$x^df(1/x)$.}\\
5: Let $x_1\!:=\!1$. If $\gamma\!=\!1$ then {\tt GOTO} Step 7.\\
6: Let $x_2\!:=\!3$. \\
7: {\tt Output} $\left\{x_1 2^{\ord_2(c_1/c_2)/d},\ldots,x_\gamma 
2^{\ord_2(c_1/c_2)/d}\right\}$. }
\end{algor}}
}}

\begin{rem} {\em 
Our correctness proof below shows that, for {\em binomials}, knowing the 
\underline{$2$} most significant base-$p$ digits of a root in $\Q_p$ 
is enough to yield an approximate root in the sense of Smale, 
{\em independent of $d$ and $H$}. Note, however, that each subsequent 
application of Newton's method to refine an approximation has complexity 
depending on $\log(dH)$ as well as $\log p$. \dia} 
\end{rem} 
\begin{rem} {\em 
We point out that the approximate roots output by our two algorithms 
above require the use of Newton iteration {\em applied to $f_{1,\zeta_0}$ 
(instead of $f$) when $p|d$}. This is clarified in our correctness proof 
below. \dia }  
\end{rem} 

\noindent
{\bf Proof of Theorem \ref{thm:binoqp}:} 
It clearly suffices to prove the correctness of Algorithms \ref{algor:binoqp} 
and \ref{algor:binoq2}, and then analyze their complexity. 

\medskip 
\noindent 
{\bf Correctness: (Case of odd $p$)}  
Theorem \ref{thm:newt} implies that Step 1 merely 
checks whether the valuations of the roots of $f$ in $\C^*_p$ in fact lie in 
$\Z$, which is necessary for $f$ to have roots in $\Q^*_p$. 
Steps 2 and 4 allow us to reduce our search for approximate roots to 
$(\Z/(p^{2\ell+1}))^*$ and assume positive degree $d$. 

\scalebox{.96}[1]{Lemma \ref{lem:binoqp} implies that 
Step 3 simply check that the coset of roots of $f$ in $\C^*_p$ 
intersects $\Z^*_p$.}  

Step 5 is the application of an automorphism of $\F^*_p$  
so we can reduce the degree of our binomial to $\gamma$, which is possibly 
much smaller than both $p-1$ and $d$.  

Steps 6--7 then clearly find the correct coset of $\F^*_p$ that makes 
$f$ vanish mod $p$. In particular, by Hensel's Lemma, Step 9 clearly gives 
the correct output if $\ell\!=\!0$. (Recall that we have replaced 
each coefficient $c_i$ of $f$ with $c'_i$.) 

If $\ell\!\geq\!1$ then let $\zeta_0$ be any $x_j$ from Step 8. We then 
have $\deg \tilde{f}_{1,\zeta_0}\!\leq\!1$ thanks 
to Lemma \ref{lem:binodepth}. 
Furthermore, Definition \ref{dfn:crazytree} tells us that the unique root 
$\zeta_1\!\in\!\F_p$ of $\tilde{f}_{1,\zeta_0}$ is exactly the next 
base-$p$ digit of a unique root $\zeta\!\in\!\Z_p$ of $f$ with 
$\zeta\!=\!\zeta_0$. 
Also, $\deg \tilde{f}_{1,\zeta_0}$ must be $1$ (for otherwise $\tf$ would 
not vanish on its coset of roots in $\F^*_p$) and $s(f,\zeta_0)\!\geq\!2$ 
since $\ell\!\geq\!1$ forces $\zeta_0$ to be a degenerate root of $\tf$. 
Lemma \ref{lem:nodal} then tells us that Hensel's Lemma 
--- applied to $f_{1,\zeta_0}(x)\!=\!p^{-s(f,\zeta_0)}f(\zeta_0+px)$ and 
start point $\zeta_1\!\in\!\Z/(p)$ --- implies that $\zeta_0+\zeta_1p$ 
yields Newton iterates rapidly converging to a true root 
$\zeta\!\in\!\Z_p$. So Step 8 in fact refines $x_1$ to the mod $p^2$ 
quantity $\zeta_0+\zeta_1 p$, and thus Steps 7--9 indeed give us suitable 
approximants in $\Q$ to all the roots of $f$ in $\Q_p$. So our algorithm 
is correct.  

Note also that the outputs, being integers in $\{0,\ldots,p^2-1\}$ 
rescaled by a factor of $p^{\ord_p(c_1/c_2)/d}$ (or possibly the reciprocals 
of such quantities), clearly each have bit-length $O\!\left(\log(p)
+\frac{|\log(c_1/c_2)|}{d\log p}\log p\right)\!=\!
O\!\left(\log(p)+\frac{\log H}{d}\right)\!=\!O\!\left(\log\!\left(
pH^{1/d}\right)\right)$. \qed  

\smallskip 
\noindent 
{\bf (Case of $p\!=\!2$)} The proof is almost the same                   
as the Correctness proof for odd $p$, save that we respectively
replace Lemma \ref{lem:binoqp} and Algorithm \ref{algor:binoqp} by
Lemma \ref{lem:binoq2} and Algorithm \ref{algor:binoq2}.
In particular, Steps 5--8 of Algorithm \ref{algor:binoqp} collapse into
Steps 5--6 of Algorithm \ref{algor:binoq2}.

So we must explain Steps 5--6 here: These steps give us the mod
$4$ reductions of the $\gamma$ many roots of $f$ in $\Z_2$, since
Steps 5 and 6 are executed only after Steps 1 and 3 certify that
$f$ indeed has roots in $\Z_2$. (Remember
that $\gamma\!\in\!\{1,2\}$ for $p\!=\!2$.) Furthermore, Hensel's Lemma
implies that the root $1$ of $\tf$ lifts to the sole root of
$f$ in $\Z_2$ when $\ell\!=\!0$. So the case $\ell\!=\!0$ is done.

If $\ell\!\geq\!1$ then there is one more complication: The nodal
polynomial $\tf_{1,1}$ is now quadratic. This is because
Lemma \ref{lem:binodepth} tells us that $\deg \tf_{1,1}\!\leq\!2$.
Furthermore, $\ell\!\geq\!1$ implies that $\gamma\!=\!2$ (assuming there 
are roots in $\Z_2$ and the algorithm hasn't terminated already) and thus $f$ 
must have exactly $2$ roots in $\Z_2$. Lemma \ref{lem:ulift} then tells
us that $\deg \tf_{1,1}\!\leq\!1$ would imply $f$ has $\leq\!1$
root in $\Z_2$. Therefore, $\tf_{1,1}$ must be quadratic.

Furthermore, $\tf_{1,1}$ must also have $2$ distinct roots: This is
because $\tf_{1,1}$ equal to $x^2$ or $1+x^2\!=\!(1+x)^2$ mod $2$ would
imply that no nodal polynomial $\tf_{i,\zeta}$, for $i\!\geq\!1$, has a
non-degenerate root. So, again by Lemma \ref{lem:ulift}, we would
not attain $2$ roots in $\Z_2$. (Similarly, it is impossible for
$\tf_{1,1}$ to be irreducible.) Therefore, the mod $4$ reductions
of the two roots of $f$ in $\Z_2$ must be $1$ and $3$. So Steps 5--6
are indeed correct.

Lemma \ref{lem:nodal} then tells us that Hensel's Lemma
--- applied to $f_{1,1}(x)\!=\!2^{-s(f,1)}f(1+2x)$ and {\em either}
start point $0$ or $1$ in $\Z/(2)$ --- implies that $1+0$ and $1+1\cdot 2$ 
yield sequences of iterates rapidly converging to true roots in $\Z_2$. So
Steps 5--7 indeed give us suitable approximants in $\Q$ to all the roots
of $f$ in $\Q_2$, and our algorithm is correct.

Note also that the outputs, being integers in $\{1,3\}$ rescaled by a
factor of $2^{\ord_2(c_1/c_2)/d}$ (or possibly the reciprocals
of such quantities), clearly each have bit-length\\
\mbox{}\hfill $O\!\left(\frac{|\log(c_1/c_2)|}{d\log 2}\log 2\right)\!=\!
O\!\left(\frac{\log H}{d}\right)\!=\!O\!\left(\log\!\left(
H^{1/d}\right)\right)$. \hfill \qed

\medskip 
\noindent 
{\bf Complexity Analysis: (Case of odd $p$)} Via Corollary \ref{cor:binomod}, 
\cite{generator}, and Theorem \ref{thm:cxity},  
it is easily checked that Steps 1--5 of  
Algorithm \ref{algor:binoqp} have respective complexity: \\  
\mbox{}\hspace{1cm}$O(\log(H)\log(pH)\log\log(pH))+O(\log(d)\log\log d)$; 
\ \  $O(\log(d)\log(dp)\log\log(dp))$;\\ 
\mbox{}\hspace{1cm}$O(\log^2(d)\log\log d)$; \ \ 
(time neglible compared to the preceding quantities); \ and\\  
\mbox{}\hspace{1cm}$O(p^{1/4}\log(p)\log\log(p))+O(\log^2(p)\log\log(p))$.\\  
This adds up to time no worse than\\ 
\mbox{}\hspace{1cm}$O(p^{1/4}\log(p)\log\log(p)
+\log(H)\log(pH)\log\log(pH) +\log(d)\log(dp)\log\log(dp))$\\ 
so far. 
Steps 6--7 (whose complexity dominates the 
complexity of Steps 6--9), involve  $\frac{p-1}{\gamma}-1$ multiplications in 
$\F_p$ and $\gamma-1$ multiplications in $\Z/(p^{2\ell+1})$. 
Since $\ell\log p\!\leq\!\log d$, this takes time no worse than 
$O(\frac{p}{\gamma}\log(p)\log\log(p)+\gamma \log(d)\log\log d)$, which 
is bounded from above by $O\!\left(\left(\frac{p}{\gamma}+\gamma\right)
\log(dp)\log\log(dp)\right)$. Note also that $\frac{p}{\gamma}+\gamma\!\geq
\!2\sqrt{p}$ by the Arithmetic-Geometric Inequality.  
So our final complexity bound is bounded from above by\\ 
\mbox{}\hspace{1.8cm}
$O\!\left(\left(\frac{p}{\gamma}+\gamma+\log d\right)
\log(dp)\log\log(dp)+\log(H)\log(pH)\log\log(pH)\right)$. \ \ \  \qed 

\smallskip
\noindent 
{\bf (Case of $p\!=\!2$)} We simply use the same techniques
as for Algorithm \ref{algor:binoqp}, save for Steps 5--8 there
being collapsed into Steps 5--6 in Algorithm \ref{algor:binoq2}. \qed 

\section{Proving Theorem \ref{thm:tri}: Trinomial Roots Never Get to Close}  
\label{sec:trisepqp} 
Let us first recall the following version of {\em Yu's Theorem}:
\begin{thm} \label{thm:yu}
\cite[Pg.\ 190]{yu07} {\em
Suppose $p$ is any prime, $n\!\geq\!2$, $\alpha_1,\ldots,\alpha_n\!\in\!\Q$
with $\alpha_i = r_i/s_i$ a reduced fraction for
each $i$, and $b_1,\ldots,b_n\!\in\!\Z$ are not all zero. Then
$\alpha_1^{b_1}\cdots\alpha_n^{b_n} \neq 1$ implies that
$\alpha_1^{b_1}\cdots\alpha_n^{b_n} - 1$ has $p$-adic valuation 
strictly less than 

\smallskip 
\noindent 
\mbox{}\hspace{1.75cm}$
\log(2)\log_p(2n)n^{5/2}(256e^2)^{n+1} p \log_p(B)
\prod^n_{i=1}\max\left\{\log|r_i|,\log|s_i|,\frac{1}{16e^2}\right\}$, 

\smallskip 
\noindent 
\scalebox{.9}[1]{where $B\!:=\!\max\{|b_1|,\ldots,|b_n|,3\}$. 
In particular, $\log(2)256e^2\!<\!1312$, $256e^2\!<\!1892$, and 
$\frac{1}{16e^2}\!<\!0.0085$. \qed}} 
\end{thm}

We will prove the square-free case of Theorem \ref{thm:tri} here, postponing 
the proof of the non-square-free case to Section \ref{sub:degen2}.  
To prove that two distinct roots $\zeta_1,\zeta_2\!\in\!\C_p$ of a 
square-free trinomial $f$ can not be too close, we will prove that $f'$ has a 
root $\tau\!\in\!\C_p$ 
with three special properties: (i) $|f(\tau)|_p$ is not too small, 
(ii) $|\zeta_1-\zeta_2|_p\!\geq\!p^{-1/(p-1)}|\zeta_1-\tau|_p$, and (iii) 
$|\zeta_1-\tau|_p$ is not too small. Step (i) is where we avail to Yu's 
Theorem, so let us now quantify our approach.  
\begin{prop} {\em \label{prop:1}
Suppose $f(x)\!=\!c_1+c_2x^{a_2} +c_3x^{a_3}\!\in\!\Z[x]$ is a trinomial of 
degree $d\!=\!a_3\!>\!a_2\!\geq\!1$, with all its coefficients having 
absolute value at most $H$, and $\tau\!\in\!\C_p$ is a root 
of $f'$. Then $\tau^{a_3-a_2}\!=\!-\frac{a_2 c_2}{a_3 c_3}$ and 
$f(\tau)\!=\!c_1+c_2\tau^{a_2}\left(1-\frac{a_2}{a_3}\right)$. \qed } 
\end{prop}
\begin{lem}{\em 
\label{lem:1} 
Following the notation above, assume further that 
$f$ is square-free. Then\\  
$\norm{f(\tau)}\!\geq\!\exp\!\left[-O(p\log_p(d)\log^2(dH))\right]$. }  
\end{lem}

\noindent 
{\bf Proof:} First note that if $f$ is square-free then $f$ has no 
repeated factors, and thus no degenerate roots in $\C_p$. So 
$f(\tau)\!\neq\!0$. Proposition \ref{prop:1} we then obtain that 
$\ord_p f(\tau)$ is 
\begin{eqnarray} 
\label{eqn:yuyu} 
\text{\scalebox{.8}[1]{
$\ord_p(c_1+c_2\tau^{a_2}(1-a_2/a_3))
=\ord_p(c_1) + \ord_p(-1) + 
\ord_p\left(\frac{-(a_3-a_2)c_2}{a_3c_1}
\left(-\frac{a_2 c_2}{a_3 c_3}\right)^{a_2/(a_3-a_2)}-1\right)$.}}   
\end{eqnarray} 

Clearly, $\ord_p c_1 \leq \frac{\log H}{\log p}$ and $\ord_p(-1)\!=\!0$. 
To bound the third summand on the right-hand side of Equality (\ref{eqn:yuyu}) 
above, let $T\!:=\!\frac{-(a_3-a_2)c_2}{a_3c_1}
\left(-\frac{a_2 c_2}{a_3 c_3}\right)^{a_2/(a_3-a_2)}$ 
and observe that\linebreak 
$T^{a_3-a_2}-1\!=\!\prod^{a_3-a_2}_{j=1}(T-\omega^j)$  
for $\omega\!\in\!\C_p$ a primitive $(a_3-a_2)$-th root 
of unity. In particular, $T^{a_3-a_2}\!\neq\!1$ since 
$f(\tau\omega^j)\!\neq\!0$ for all $j\!\in\!\{1,\ldots,a_3-a_2\}$, thanks to 
Proposition \ref{prop:1} and $f$ not having any degenerate roots. So then 
$M\!:=\!\ord_p(T^{a_3-a_2}-1) = \sum_{j=1}^{a_3-a_2}\ord_p(T-\omega^j)
<\infty$, with the $(a_3-a_2)$-th term of the sum exactly 
$\ord_p(T-\omega^{a_3-a_2})\!=\!\ord_p(T-1)$, 
i.e., the third summand from Equality (\ref{eqn:yuyu}).  

Suppose $\ord_p T\!<\!0$. Then for each $i\!\in\!\{1,\ldots,a_3-a_2\}$ we 
have $\ord_p(T-\omega^j)\!=\!\ord_pT\!<\!0$, since roots 
of unity always have $p$-adic valuation $0$. We must then have 
$\ord_p f(\tau)\!=\!\ord_p(c_1) + \ord_p(T-\omega^{a_3-a_2}) 
\!<\!\frac{\log_p(dH)}{1}$ (by Theorem \ref{thm:newt}) and we obtain our lemma.

On the other hand, should $\ord_pT \geq 0$, we get
$\ord_p(T-\omega^j)\geq j\ord_p(\omega)\!=\!0$, for each $j$. 
So $M\!\geq\!\ord_p(T-1)$ and we'll be done if we find a 
sufficiently good upper bound on $M$.

By luck, $M$ is boundable directly from Yu's Theorem (Theorem \ref{thm:yu} 
here) upon setting $n\!=\!2$, $\alpha_1\!=\!-\frac{(a_3-a_2)c_2}{a_3c_1}$, 
$\alpha_2\!=\!-\frac{a_2c_2}{a_3 c_3}$, $b_1\!=\!a_3-a_2$, and $b_2\!=\!a_2$. 
In particular, we can assume $|r_i|,|s_i|\!\leq\!dH$ for $i\!\in\!\{1,2\}$ and 
$B\!=\!\max\{d,3\}$, and move the $\log p$ factors in the denominator 
so that 
$M\!<\!\log(2)256e^2\log(4)2^{5/2}(256e^2)^2p\log\max\{d,3\}
\left(\max\left\{\log_p(dH),
\frac{1}{16e^2\log p}\right\}\right)^2$. For $d\!=\!2$ we 
get $f(\tau)\!=\!\frac{c_1}{4c_3}(4c_1c_3-c^2_2)$, which is a rational 
number that this an integer of absolute value at most $H^2+4H$ divided by an 
integer of absolute value at most $4H$. Such a rational number clearly has 
valuation no greater than $\log_p(H^2+4H)\!=\!O(\log_p H)$ and thus 
$|f(\tau)|_p\!\geq\!e^{-O(\log H)}$ when $d\!=\!2$. 
Since $d\!\geq\!2$ for an arbitrary trinomial, and $H\!\geq\!1$,
we then obtain $M\!<\!36791093348 p\log(d)\log^2_p(dH)\!=\!O(p\log(d)
\log^2_p(dH))$. In other words, the third summand from (1) is bounded from 
above by the last $O$-bound, and thus $\ord_p f(\tau)\!=\!O(M)$ since 
$\frac{\log H}{\log p}\!=\!O(M)$. Since 
$|f(\tau)|_p\!=\!e^{-\log(p)\ord_p f(\tau)}$, 
we are done. \qed  

\medskip 
The Ultrametric Inequality directly yields the following:  
\begin{prop}{\em \label{prop:dev_small} 
If $f\!\in\!\Z[x]$ and $t\!\in\!\C_p$ then 
$\norm{t}\leq 1\Longrightarrow \norm{f'(t)}\leq 1$. \qed } 
\end{prop}

\medskip 
\scalebox{.96}[1]{Below is a rescaled {\em $p$-adic} version of 
{\em Rolle's Theorem}, based on \cite[Sec.\ 2.4, Thm., Pg.\ 316]{Rob00}.}   
\begin{thm}\label{thm:rolle} {\em Let $f\in \mathbb{C}_p[x]$ have two distinct 
roots $\zeta_1, \zeta_2\in \mathbb{C}_p$ with 
$\norm{\zeta_1-\zeta_2}\!=\!cp^{1/(p-1)}$ for some $c\!>\!0$. 
Then $f'$ has a root $\tau\!\in\!\C_p$ with $|\zeta_1-\tau|_p,
|\zeta_2-\tau|_p\!\leq\!c$. \qed } 
\end{thm}

We can now prove part of one of our main results.  

\medskip 
\noindent 
{\bf Proof of the Square-Free Case of Theorem \ref{thm:tri}:} 
Note that $\zeta_i\!\neq\!0 \Longrightarrow |\ord_p \zeta_i|\!\leq\!\log_p H$ 
thanks to Theorem \ref{thm:newt}. So then 
$\ord_p(\zeta_1-\zeta_2)\!\geq\!-\log_p H$ for any pair of 
distinct roots $\zeta_1,\zeta_2\!\in\!\C_p$ of $f$ and, if 
$\zeta_1\zeta_2\!=\!0$, we also have  
$\ord_p(\zeta_1-\zeta_2)\!\leq\!\log_p H$. 
So $\log H\!\geq\!\log|\zeta_1-\zeta_2|_p$ and, if $\zeta_1\zeta_2\!=\!0$ 
then we also have $\log|\zeta_1-\zeta|_p\!\geq\!-\log H$. 
So we may assume $\zeta_1\zeta_2\!\neq\!0\!\neq\!f(0)$.  

\scalebox{.95}[1]{For convenience, let us abbreviate the first (larger) 
$O$-bound stated in our theorem by $O(M)$.}  

\medskip 
\noindent 
{\bf Case 1: (Both roots are small: $\pmb{\norm{\zeta_1}, 
\norm{\zeta_2}\leq 1}$.)}  \\
Suppose $\norm{\zeta_1-\zeta_2}\!>\!p^{-2/(p-1)}$ ($=\!
e^{-2\log(p)/(p-1)}$). Since $2\log(p)/(p-1)=O(M)$ we are done.   

Now assume that $\norm{\zeta_1-\zeta_2}\leq p^{-2/(p-1)}$. Then by Theorem  
\ref{thm:rolle} $f'$ has a root $\tau\!\in\!\mathbb{C}_p$ with  
$\norm{\zeta_i-\tau} \leq p^{1/(p-1)}\norm{\zeta_1-\zeta_2} 
\leq p^{-1/(p-1)}$ for all $i\!\in\!\{1,2\}$. 
Since $f$ is square-free, Lemma \ref{lem:1} implies that 
$\norm{f(\tau)} \geq e^{-O(M)}$. Applying Theorem \ref{thm:rolle} to\\ 
\mbox{}\hfill $g(x)\!:=\!f(x)-\frac{f(\tau)-f(\zeta_1)}{\tau-\zeta_1}x
-\frac{\tau f(\zeta_1) -\zeta_1 f(\tau)}{\tau-\zeta_1}$\hfill\mbox{}\\  
(which vanishes at $\tau$ and $\zeta_1$), we then see that there is a 
$\mu\!\in\!\C_p$ with $\norm{\mu-\zeta_1}\leq 1$ 
such that $g'(\mu)\!=\!0$, i.e., 
$f(\tau) = f(\tau) - f(\zeta_1) = f'(\mu)(\tau-\zeta_1)$. 
Note that $|\mu|_p\!\leq\!1$ since $|\mu|_p\!>\!1$ would 
imply that $|\mu|_p\!>\!|\zeta_1|_p$ and thus 
$|\mu-\zeta_1|_p\!=\!|\mu|_p\!>\!1$, giving us a contradiction. As 
$f(\tau)\!\neq\!0$ we get $f'(\mu)\neq 0$ and $\tau\!\neq\!\zeta_1$. 
From Proposition \ref{prop:dev_small} 
we have $\norm{f'(\mu)}\!\leq\!1$, so then $\norm{\tau-\zeta_1}\!=\!  
\frac{\norm{f(\tau)}}{\norm{f'(\mu)}} \geq e^{-O(M)}$. 
We thus get  $\norm{\zeta_1-\zeta_2} \geq p^{-1/(p-1)}\norm{\tau-\zeta_1}\!
\geq\!e^{-O(M)-\frac{\log p}{p-1}}\!=\!e^{-O(M)}$. \qed 

\medskip 
\noindent
{\bf Case 2: (Both roots are large: $\pmb{\norm{\zeta_1}, 
\norm{\zeta_2} > 1}$.)} Simply observe that $1/\zeta_1$ and $1/\zeta_2$ are 
roots of the {\em reciprocal polynomial} $f^*(x)\!:=\!x^{\deg f}f(\frac{1}{x})$.
In particular, we can apply Case 1 to the trinomial $f^*$ since
$\norm{\frac{1}{\zeta_1}}, \norm{\frac{1}{\zeta_2}}<1$.
We then obtain $\norm{\frac{1}{\zeta_1} - \frac{1}{\zeta_2}} 
\geq e^{-O(M)}$. Hence $\norm{\zeta_1-\zeta_2} = \norm{\zeta_1}
\norm{\zeta_2} \norm{\frac{1}{\zeta_1} - \frac{1}{\zeta_2}}
\geq \norm{\frac{1}{\zeta_1} - \frac{1}{\zeta_2}}\geq 
e^{-O(M)}$. \qed
 
\medskip  
\noindent
{\bf Case 3: (Only one root has norm $\pmb{>1}$.)}  \\ 
Without loss of generality, we may assume that $|\zeta_1|_p\!\leq\!1\!<\!
|\zeta_2|_p$. We then simply note that, as $\norm{\zeta_1}\neq \norm{\zeta_2}$, 
we have $\norm{\zeta_1-\zeta_2} = \max\left\{\norm{\zeta_1}, \norm{\zeta_2}
\right\}\!>\!1$  and we are done. \qed 

\section{Proving Theorem \ref{thm:tetra}: Tetranomial Roots Can Get Too Close} 
\label{sec:tetra}

\subsection{The Case of Prime $\pmb{p}$} 
Let $g(x) = p^{2h}f(x+p^{h-1}) = p^{2h}(x+p^{h-1})^d 
- p^{2h}\left(\frac{x+p^{h-1}}{p^h} 
- \frac{1}{p}\right)^2$  
$=p^{2h}(x+p^{h-1})^d - x^2$. 
Then $g$ has the same roots as $f_{d,p}$, save for a ``small'' shift by 
$p^{h-1}$. Rescaling, we get \scalebox{.95}[1]
{$G(x):=\frac{g(p^{(h-1)d/2+h}x)}{p^{(h-1)d+2h}}
=p^{-(h-1)d-2h}\left[ p^{2h}(p^{(h-1)d/2+h}x+p^{h-1})^d - 
p^{(h-1)d+2h}x^2\right]$}\linebreak 
$=\sum_{i=0}^{d}{d\choose i}p^{(h-1)(di/2-i)+ih}x^i - x^2=
1-x^2 \mod p^{d(h-1)/2+1}$, 
which is square-free for odd prime $p$. So if $p$ is odd, then 
Hensel's Lemma implies that there are 
roots $\zeta_1,\zeta_2 \in \Z_p$ of $G$ such that 
$\zeta_1 \equiv 1 \mod p^{d(h-1)/2+1}$ and $\zeta_2 \equiv -1 \mod 
p^{d(d-1)/2+1}$. 

On the other hand, if $p\!=\!2$, then, as $h\!>\!2$, we have $p^{d(h-1)/2+1} 
\!\geq\!8$. So, since $G(x)\!=\!1-x^2\!=\!(3-x)(5-x) \mod 2^3$, we obtain that 
$G$ is square-free in $\Z_2[x]$. Hensel's Lemma then implies that there are 
roots $\zeta_1,\zeta_2 \in \Z_p$ of $G$ such that 
$\zeta_1\!=\!3$ mod $p^{d(h-1)/2+1}$ and $\zeta_2\!=\!5$ mod 
$p^{d(h-1)/2+1}$. 

So, whether $p$ is odd or even, we obtain two roots 
$x_1,x_2\!\in\!\Z_p$ of $G$ with $\norm{x_1}\!=\!\norm{x_2} 
\!=\!1$. For each $i\!\in\!\{1,2\}$, $y_i\!=\!p^{(h-1)d/2+h}x_i$ is then 
the corresponding root of $g$. So $\zeta_1\!:=\!y_1 + p^{h-1}$ and 
$\zeta_2\!:=\!y_2 + p^{h-1}$ are two roots of $f$ in $\Z_p$ such that
$\norm{\zeta_1-\zeta_2}\!=\!\norm{(y_1+p^{h-1})-(y_2+p^{h-1})}\!=\!
\norm{y_1-y_2} \leq \max\left\{\norm{y_1}, \norm{y_2}\right\}\!=\!
p^{-(h-1)d/2-h} = p^{-\Omega(dh)}$. \qed  

\begin{rem} \label{rem:tetra} {\em 
From our proof, we see that $f_{d,p}$ has two roots of the form\\ 
\mbox{}\hfill $\zeta_i\!=\!p^{h-1}+\eps_ip^{(h-1)d/2}+O(p^{1+(h-1)d/2})$ 
\hfill\mbox{}\\ 
with $i\!\in\!\{1,2\}$ and $\{\eps_1,\eps_2\}$ equal to $\{\pm 1\}$ 
or $\{3,5\}$, according as $p$ is odd or even. In particular, by 
direct evaluation, it is easily checked that 
$\ord_p f'_{d,p}(\zeta_i)\!=\!\ord_p(d)+(h-1)(d-1)$. In other 
words, we can need as many as $\Omega(d\log H)$ of the most significant 
base-$p$ digits of a root of a tetranomial in order to use it as a 
start point for Newton iteration. We will see in  
Section \ref{sec:central} that $O_p(\log^3(\max\{d,H\})\log(d))$ 
base-$p$ digits suffice for trinomials. \dia } 
\end{rem} 

\subsection{The Case $\pmb{p\!=\!\infty}$} 
Shifting by $\frac{1}{2^{h-1}}$, we get $g(x):=f_{d,\frac{1}{2}}(x+2^{1-h}) 
= (x+2^{1-h})^d - 2^{2h}x^2$\\ 
$=2^{d(1-h)} + d2^{(d-1)(1-h)}x+\left({d\choose 2} 2^{(d-2)(1-h)} 
-2^{2h}\right)x^2 + {d\choose 3} 2^{(d-3)(1-h)} x^3 + \cdots+ x^d$.
We will see momentarily that, unlike $\anewt(f)$ (which has $3$ lower
edges), $\anewt(g)$ will have just $2$ lower edges.
(See the right-hand illustration in Example \ref{ex:newts}.)
This will force (via Theorem \ref{thm:newt}) the existence of two distinct
roots of small norm for $g$, thus yielding two nearby roots of $f$ after
undoing our earlier shift.

Toward this end, note that the three lowest order terms of $g$ contribute
the points\linebreak
\scalebox{.91}[1]{$p_0:= (0,d(h-1)\log2)$, $p_1:=(1,(d-1)(h-1)\log2-\log d)$,
and $p_2 = \left(2,-\log\left(4^h - \frac{{d\choose 2}}
{2^{(d-2)(h-1)}} \right)\right)$}\linebreak
as potential vertices of $\anewt(g)$.
Observe that $\frac{{d\choose 2}}{2^{(d-2)(h-1)}}\!<\!0.059$ for all
$h\!\geq\!3$ and $d\!\geq\!4$, and thus $p_2$ is the only point of $\anewt(f)$
with negative $y$-coordinate. So $p_2$ is a vertex of $\anewt(f)$, and
all edges with vertices to the right of $p_2$ have positive slope.
Furthermore, the slopes of the line segments $\overline{p_0p_1}$
and $\overline{p_0p_2}$ are respectively $-(h-1)\log(2)-\log d$ and
a number less than $-\frac{1}{2}\log(4^h-0.059)-\frac{1}{2}d(h-1)\log 2$.

Since $2^{h-1}\!<\!\sqrt{4^h-0.059}$ and $\log d\!<\!\frac{1}{2}d(h-1)\log 2$
for all $d\!\geq\!4$ and $h\!\geq\!3$, we thus see that the slope
of $\overline{p_0p_2}$ is more negative. So the leftmost lower edge of
$\anewt(g)$ has vertices $p_0$ and $p_2$. It is easily checked that
the slope of this edge is less than $-10.3$, which is in turn clearly
$<\!-2\log 3$. So by Theorem \ref{thm:newt}, there are
two roots $z_1,z_2$ of $g$ such that
\begin{align*}
\log|z_i| \leq \frac{1}{2}\left[-\log\left(2^{2h} - {d\choose 2}
2^{(d-2)(1-h)}\right)  - d(h-1)\log2\right].
\end{align*}
These two roots thus satisfy $|z_i| = 2^{-\Omega(dh)}$. Now, for
$i\!\in\!\{1,2\}$, $\zeta_i = z_i+2^{1-h}$ yields roots of $f_{d,\frac{1}{2}}$
with $|\zeta_1-\zeta_2| = |z_1+2^{1-h}-(z_2+2^{1-h})| \leq |z_1|+|z_2| < 
2^{-\Omega(dh)}$. \qed

\section{Valuation Bounds from Discriminants and Repulsion From Degeneracy}  
\label{sec:central}  
While we we were able to prove a special case of our bound for the 
minimal root spacing of trinomials, we will need to examine the 
roots in $\C^*_p$ more carefully for trinomials that have degenerate 
roots in $\C^*_p$. We will see that the roots appear to repel more 
strongly in the degenerate case, and a key tool to prove this  
is the {\em trinomial discriminant}. 
\begin{dfn} \label{dfn:disc} {\em 
\cite{gkz94} Suppose $f(x)\!=\!c_1+c_2x^{a_2}+c_3x^{a_3}\!\in\!\Z[x]$ is a 
trinomial with $a_3\!>\!a_2\!\geq\!1$, $r\!:=\!\gcd(a_2,a_3)$, and 
$\ba_i\!:=\!\frac{a_i}{r}$ for all 
$i$. We then define the {\em trinomial discriminant} to be\\ 
\mbox{}\hfill  
$\Delta_{\mathrm{tri}}(f)\!:=\!\ba^{\ba_3}_3c^{\ba_3-\ba_2}_1 c^{\ba_2}_3 
- \ba^{\ba_2}_2 (\ba_3-\ba_2)^{\ba_3-\ba_2} (-c_2)^{\ba_3}$. \hfill \dia}
\end{dfn} 

\noindent 
Up to a sign factor, our definition agrees with the  
definition of the {\em $\{0,a_2,a_3\}$-discriminant} from \cite[Ch.\ 9, pp.\ 
274--275, Prop.\ 1.8]{gkz94} when $\gcd(a_2,a_3)\!=\!1$.   
We will also need to recall the following facts: 
\begin{lem} {\em 
\label{lemma:tri} 
\cite[Lemma 40]{airr}
Following the notation of Definition \ref{dfn:disc}:\\  
(1) If $c_1c_3\!\neq\!0$ then 
$\Delta_{\mathrm{tri}}(f)\!\neq\!0\Longleftrightarrow f$ has no degenerate
roots in $\C_p$. Furthermore,\linebreak 
\mbox{}\hspace{.8cm}$p\nmid c_1c_3\gcd(a_2,a_3)$ also implies the equivalence  
$\Delta_{\mathrm{tri}}\!\left(\tf\right)\!\neq\!0$ mod $p \Longleftrightarrow 
\tf$ has no\linebreak 
\mbox{}\hspace{.8cm}degenerate roots in $\bF_p$. \\ 
(2) If $\Delta_{\mathrm{tri}}(f)\!\neq\!0$ then 
\scalebox{.93}[1]{$\Delta_{\mathrm{tri}}(f)=\left(\frac{c_3}{c_1}
\right)^{\ba_2-1}\!\!\!\!\!\!\!\!\!\!\!\!\prod\limits_{\xi\in\C_p \; : \; 
\barf(\xi)=0}\barf'
(\xi)=(-1)^{\ba_3(\ba_3-\ba_2)}\!\!\!\!\!\!\!\!\!\!\prod\limits_{\xi\in\C_p 
\; : \; \barf(\xi)=0} \left(\ba_2c_2+\ba_3c_3\xi^{\ba_3-\ba_2}\right)$}
\linebreak 
\mbox{}\hspace{.8cm}where $\barf\!\in\!\Z[x]$ is the unique polynomial 
satisfying $f(x)\!=\!\barf(x^r)$ identically. \qed}   
\end{lem}
\begin{rem}
{\em The second sentence of Assertion (1) appears not to be well-known but 
does follow easily from the development of \cite[Ch.\ 9]{gkz94}, upon observing 
that $p\nmid \gcd(a_2,a_3) \Longrightarrow$ the matrix \scalebox{.7}[.7]
{$\begin{bmatrix} 1 & 1 & 1\\ 0 & a_2 & a_3\end{bmatrix}$} has rank $2$. 
Should $p|\gcd(a_2,a_3)$ then it is easily checked that every root in 
$\F^*_p$ of the trinomial $\tf$ above is degenerate. \dia}  
\end{rem} 

Recall that the {\em classical degree $d$ discriminant} of a polynomial 
$g(x)\!=\!c_0+\cdots+c_dx^d\!\in\!\C_p[x]$ is $\Delta_d(g)\!:=\!
\frac{\res_{d,d-1}(f,f')}{c_d}$ where $\res_{d_1,d_2}(g_1,g_2)$ denotes the 
well-known {\em resultant} of two univariate polynomials, $g_1$ and $g_2$, 
having respective degrees $\leq\!d_1$ and $\leq\!d_2$ (see, e.g., 
\cite[Ch.\ 12]{gkz94}). We will also need some deeper facts 
about the discriminants of trinomials, and prove repulsion from 
degenerate roots along the way: 
\begin{lem} {\em 
\label{lemma:tri2} 
Suppose $f(x)\!=\!c_1+c_2x^{a_2}+c_3x^{a_3}\!\in\!\Z[x]$ has 
degree $d\!=\!a_3\!>\!a_2\!\geq\!1$, $c_1c_2c_3\!\neq\!0$, and  
$|c_i|\!\leq\!H$ for all $i$. 
Assume further that $f$ has a degenerate root $\tau\!\in\!\C_p$, 
$r\!:=\!\gcd(a_2,a_3)$, and $\ba_i\!:=\!\frac{a_i}{r}$ for all $i$.  
Finally, let\\ 
\vbox{$Q(x)\!:=\!(\ba_3-\ba_2)\left(1+2x+3x^2+\cdots+(\ba_2-1)x^{\ba_2-2}
\right)$\\ 
\mbox{}\hfill $+\ba_2\left((\ba_3-\ba_2)x^{\ba_2-1}+(\ba_3-\ba_2-1)x^{\ba_2}
+\cdots+1\cdot x^{\ba_3-2} \right)$}\\ 
and $q(x)\!:=\!(\ba_3-\ba_2) -\ba_3x^{\ba_2}+\ba_2x^{\ba_3}$. Then:\\   
(1) Any degenerate root $\tau\!\in\!\C_p$ of $f$ satisfies $\tau^r\!\in\!\Q^*$ 
and $(\tau^{a_2},\tau^{a_3})\!=\!\frac{c_1}{a_3-a_2} 
\left(-\frac{a_3}{c_2},\frac{a_2}{c_3}\right)$.\linebreak 
\mbox{}\hspace{.7cm}Furthermore, if $p\nmid (a_3-a_2)c_1$, then 
any degenerate root $\tilde{\tau}\!\in\!\bF_p$ of $\tf$ 
satisfies\\ 
\mbox{}\hspace{.7cm}$(c_2\tilde{\tau}^{a_2},
c_3\tilde{\tau}^{a_3})\!=\!\frac{c_1}
{a_3-a_2}(-a_3,a_2)$ and, if $p\nmid c_2c_3$ in addition, then 
$\tilde{\tau}^r\!\in\!\F^*_p$.\\  
(2) \scalebox{.92}[1]{The polynomial $q$ has $1$ as its unique degenerate root 
in $\C_p$ and $q(x)\!=\!Q(x)(x-1)^2$ identically.}\\ 
(3) We have $Q(1)\!=\!\ba_2\ba_3(\ba_3-\ba_2)/2$ and,  
for $\ba_3\!\geq\!4$, $\Delta_{\ba_3-2}\!\left(Q\right)=
\ba_3 (\ba_2\ba_3(\ba_3-\ba_2))^{\ba_3-4}J$,\\ 
\mbox{}\hspace{.7cm}where 
$J\!=\!O(\ba^2_2\ba^3_3(\ba_3-\ba_2)^2)$ is a nonzero integer.\\ 
(4) For $\ba_3\!\geq\!4$ we have 
$\Delta_{\ba_3-2}(Q)\!=\!\ba^{\ba_3-4}_2 \!\!\!\!\!\! 
\prod\limits_{\mu\in\C_p \; : \; Q(\mu)=0} Q'(\mu)$.\\
(5) \scalebox{.97}[1]{$|\ord_p(\zeta-\tau)|\!\leq\!\log_p\frac{(d-r)d^3
H}{8r^4}\!<\!4\log_p\frac{dH^{1/4}}{r}$ for any non-degenerate 
root $\zeta\!\in\!\C_p$ of $f$.} 
}  
\end{lem}  

\noindent 
{\bf Proof of Lemma \ref{lemma:tri2}:} Assertions (1)--(3) are immediate upon 
applying \cite[Lemma 40]{airr} to the polynomial $\bar{f}$ from Lemma 
\ref{lemma:tri} (which satisfies $f(x)\!=\!\barf(x^r)$). Assertion (4) follows 
similarly from \cite[Product Formula, Pg.\ 398]{gkz94}, which is a product 
formula for resultants. Assertion (5) will follow routinely upon proving that 
the roots of $Q$ can't be too close to $1$, and that the same holds for the 
$1/r$-th powers of the roots of $Q$ as well. In particular, we'll\linebreak 
\scalebox{.92}[1]{soon see 
that the $r$th powers of the non-degenerate roots of $f$ are mild rescalings 
of the roots of $Q$.}  

\smallskip 
\noindent 
{\bf Assertion (5):} To simplify matters, we will first reduce to the case 
$r\!=\!1$. Since the polynomial $\barf$ from Lemma \ref{lemma:tri} 
is an instance of the case $r\!=\!1$, and the roots of $\barf$ are 
the $r$th powers of the roots of $f$, we can perform our reduction  
by showing that a sufficiently good upper bound on $|\ord_p(\zeta^r-\tau^r)|$ 
implies our desired upper bound on $|\ord_p(\zeta-\tau)|$. So first note that 
if $\ord_p\zeta\!\neq\!\ord_p\tau$ then $\ord_p(\zeta-\tau)
\!=\!\min\{\ord_p\zeta,\ord_p\tau\}$. In particular, 
since $a_3\!=\!r\ba_3$, and $a_2$ and $a_3-a_2$ are positive multiples of $r$, 
Theorem \ref{thm:newt} implies:
\begin{eqnarray}
\label{fact:nicenorm}
\mbox{}\hspace{.8cm} 
\text{\scalebox{.9}[1]{{\em {\em Any} root of $f$ in $\C_p$ must have 
valuation in the closed interval $\left[\frac{\ord_p(c_2/c_3)}{r},
\frac{\ord_p (c_1/c_2)}{r}\right]$}}}\\ 
\text{{\em or have valuation exactly 
$\frac{\ord_p(c_1/c_3)}{r\ba_3}$, according as $\ord_p\frac{c^2_2}{c_1c_3}\!
\leq\!0$ or not.}} 
\hspace{1.65cm}\mbox{} \nonumber   
\end{eqnarray} 
So $|\ord_p (\zeta-\tau)|\!\leq\!\frac{\log_p H}{r}\!<\!
\log_p\frac{(d-r)d^3H}{8r^4}$, and the last inequality clearly 
holds when $\frac{d}{r}\!\geq\!2$. 
We may thus \fbox{assume $\ord_p\zeta\!=\!\ord_p\tau$}. 

Now, if $r\!>\!1$, then we can observe that 
\begin{eqnarray} 
\label{eqn:trick} 
\ord_p(\zeta^r-\tau^r) & = & r\ord_p(\zeta)+
\ord_p\left(1-\left(\frac{\tau}{\zeta}\right)^r\right).
\end{eqnarray} 
Letting $\omega\!\in\!\C_p$ be any primitive $r$th root of unity, we then 
obtain $\ord_p\left(1-\left(\frac{\tau}{\zeta}\right)^r\right)
\!=\!\sum^{r-1}_{j=0}\ord_p\left(1-\frac{\tau \omega^j}{\zeta}\right)$. 
Since each term in the preceding sum is clearly nonnegative we must then 
have $\ord_p\left(1-\frac{\tau}{\zeta}\right)\!\leq\!
\ord_p\left(1-\left(\frac{\tau}{\zeta}\right)^r\right)$. 
So if we have $\ord_p(\zeta^r-\tau^r)\!\leq\!M$ for some 
$M\!\geq\!r\ord_p\zeta$ then Equality (\ref{eqn:trick}) implies 
$\left|\ord_p\left(1-\frac{\tau}{\zeta}\right)\right|\!\leq\!M-r\ord_p\zeta$. 
Fact (\ref{fact:nicenorm}) then implies\\ 
\mbox{}\hfill  
$|\ord_p(\zeta-\tau)|\!=\!\left|\ord_p(\zeta)+\ord_p\left(1
-\frac{\tau}{\zeta}\right)\right|\!\leq\!M-(r-1)\ord_p \zeta\!\leq\!M
+\frac{r-1}{r}\log_p H$.\hfill\mbox{}\\ 
Since $\frac{1}{r}+\frac{r-1}{r}\!=\!1$, we will clearly establish 
Assertion (5) if we can prove $\ord_p\left(\zeta^r-\tau^r\right)\!\leq\!
\log_p\frac{(d-r)d^3H^{1/r}}{8r^4}$. Since every root of $\barf$ is the $r$th 
power of a root of $f$ (and vice-versa), and since 
$\deg \barf\!=\!\frac{d}{r}$ and 
$\gcd(\ba_2,\ba_3)\!=\!1$, Fact (\ref{fact:nicenorm}) implies that it suffices 
to prove the following half of the $r\!=\!1$ case of Assertion (5): 
$\ord_p (\zeta-\tau)\!\leq\!
\log_p\frac{(d-1)d^3H}{8}$. (Our stated bound is implied by the preceding 
bound since $\ord_p\zeta\!=\!\ord_p\tau\Longrightarrow 
\ord_p(\zeta-\tau)\!\geq\!0$.)  
We will thus \fbox{assume $\gcd(a_2,a_3)\!=\!1$ henceforth}. 

\smallskip 
\noindent 
{\em (The Case $d\!\in\!\{2,3\}$)} Note that $d\!\geq\!2$ because
$f$ is a trinomial. The case $d\!=\!2$ is then vacuously
true since a quadratic with a degenerate root has no non-degenerate roots.  

For $d\!=\!3$, Assertion (2) of our lemma tells us that there is 
only one non-degenerate root $\zeta$ and it is rational. So, evaluating 
the factorization of $f$ at $0$, we must have 
$\tau^2\zeta\!=\!-\frac{c_1}{c_3}$. Assertion (1) of our 
lemma tells us that $\tau^3\!=\!\frac{c_1a_2}{(3-a_2)c_3}$ and thus  
$\frac{\zeta}{\tau}\!=\!-\frac{3-a_2}{a_2}$. So we obtain 
$\ord_p(\tau-\zeta)\!=\!\ord_p(\tau)+\ord_p(1-\frac{\zeta}{\tau})\!=\!
\frac{\ord_p((c_2a_2)/(3c_3))}{3-a_2}
+\ord_p\left(\frac{3-a_2}{a_2}\right)$, where the last equality 
follows from Theorem \ref{thm:newt} applied to $f'$.  
Since $|c_2a_2|\!\leq\!2H$ and $3-a_2\!\leq\!2$, it easily follows that 
$\ord_p(\tau-\zeta)\!\leq\!\log_p(4H)\!<\!\log_p \frac{(d-1)d^3H}{8}$. 
Our assertion thus holds when $d\!\leq\!3$. \qed  

\smallskip 
\noindent 
{\em (The Case $d\!\geq\!4$)} We will first prove an upper bound on 
$\ord_p(1-\mu)$ for all roots $\mu\!\in\!\C_p\setminus\{1\}$ of $q$ . 
Observe that Assertion (2) and the classical theory of discriminants 
\cite[Ch.\ 12]{gkz94} imply that $Q$ has exactly $a_3-2$ distinct roots in 
$\C^*_p$ and $\Delta_{a_3-2}(Q)\!\neq\!0$. The first half of Assertion (3) 
then tells us that $\prod\limits_{\mu\in\C_p \; : \; Q(\mu)=0} 
(1-\mu)\!=\!\frac{Q(1)}{a_2}\!=\!\frac{a_3 (a_3-a_2)}{2}$,  
since the leading coefficient of $Q$ is $a_2$. So then  
\begin{eqnarray} 
\label{eqn:ordsum}  
\sum\limits_{\mu\in\C_p \; : \; Q(\mu)=0}\ord_p(1-\mu) & =  & \ord_p\left(
\frac{a_3(a_3-a_2)}{2}\right)\!\leq\!
\log_p\left(\frac{a_3(a_3-a_2)}{2}\right)\!\leq\!\log_p\binom{d}{2}. 
\end{eqnarray}  

Thanks to Theorem \ref{thm:newt}, $\ord_p a_2\!=\!0$ (i.e., the 
leading coefficient of $Q$ not being divisible by $p$) implies that all the 
roots $\mu\!\in\!\C_p$ of $Q$ have nonnegative valuation. So then 
$\ord_p(1-\mu)\!\geq\!0$ and, thanks to Bound (\ref{eqn:ordsum}), we obtain 
$\ord_p(1-\mu)\!\leq\!\log_p\binom{d}{2}\!<\!\log_p\frac{(d-1)d^3\cdot d}{8}$. 
(Note that the coefficients of $q$ have absolute value at most $d\!=\!a_3$.) 
So we may assume $\sigma\!:=\!\ord_p a_2\!>\!0$ henceforth. 

Since $\gcd(a_2,a_3)\!=\!1$ we must have 
$\ord_p a_3\!=\!0\!=\!\ord_p(a_3-a_2)$. Theorem 
\ref{thm:newt} applied to $q$ then tells us that $Q$ has exactly $a_3-a_2$ 
roots in $\C_p$ of $p$-adic valuation $-\frac{\sigma}{a_3-a_2}$, and 
exactly $a_2-2$ roots $\mu\!\in\!\C_p$ of $p$-adic valuation $0$, since 
$q(x)\!=\!Q(x)(x-1)^2$. In particular, 
$\ord_p(1-\mu)\!=\!-\frac{\sigma}{a_3-a_2}\!\geq\!-\log_p(d-1)$ 
on the set of roots with negative valuation, and $\ord_p(1-\mu)\!\geq\!0$ at 
the roots $\mu\!\in\!\C_p$ with $\ord_p\mu\!=\!0$. 

Equality (\ref{eqn:ordsum}) then implies that each of the $a_3-2$ roots 
$\mu\!\in\!\C_p$ of $Q$ with $\ord_p \mu\!=\!0$ must satisfy 
$\ord_p(1-\mu)\!=\!(a_3-a_2)\frac{\sigma}{a_3-a_2}
+\ord_p\left(\frac{a_3(a_3-a_2)}{2}\right)
\!=\!\ord_p\left(\frac{a_2a_3(a_3-a_2)}{2}\right)\!\leq\!\log_p\left(
\frac{a_3a_2(a_3-a_2)}{2}\right)$. 
By the Arithmetic Geometric Inequality, $a_2(a_3-a_2)\!\leq\!a^2_3/4$,   
so we arrive at $\ord_p(1-\mu)\!\leq\!\log_p(d^3/8)\!<\!
\log_p((d-1)d^3\cdot d/8)$ and we have proved  
Assertion (5) in the special case $f(x)\!=\!q(x)$. 

A direct computation via Assertion (1) of our lemma then yields 
$f(x)\!=\!\frac{c_1}{(a_3-a_2)\tau^2}q(x/\tau)$ identically. So 
the roots of $f$ are simply scalings of the roots of $q$ by a factor 
$\tau$. Since $f'(\tau)\!=\!0$, Theorem \ref{thm:newt} implies that 
$\ord_p \tau\!=\!\frac{\ord_p(a_2c_2) - \ord_p(a_3c_3)}{a_3-a_2}$, 
which clearly lies in the closed interval $[-\log_p(dH),\log_p((d-1)H)]$.  
So then $\ord_p(\tau-\zeta)\!=\!\ord_p\tau+\ord_p(1-\mu)$ for\linebreak 
\scalebox{.95}[1]{some root 
$\mu\!\in\!\C_p$ of $Q$. In other words, 
$\ord(\tau-\zeta)\!\leq\!\log_p((d-1)Hd^3/8)\!=\!\log_p((d-1)d^3H/8)$. 
\qed}

\medskip 
Assertion (1) of Lemma \ref{lemma:tri2} tells us that degenerate roots in 
$\C^*_p$ of trinomials satisfy binomial equations with well-bounded 
coefficients. Our earlier Algorithms \ref{algor:binoqp} and \ref{algor:binoq2} 
thus imply that degenerate roots of trinomials are easy to approximate. Our 
final step in proving Theorem \ref{thm:tri} will be estimating 
the spacing of {\em non}-degenerate roots in $\C_p$ for trinomials having 
degenerate roots in $\C_p$.

\subsection{Completing the Proof of Theorem \ref{thm:tri}: Degenerate Root 
Spacing} 
\label{sub:degen2} 
First note that we may assume $\zeta_1\zeta_2\!\neq\!0\!\neq\!
f(0)$, since this initial reduction to nonzero roots (from the proof of the 
square-free case in Section \ref{sec:trisepqp}) 
does not require $f$ to be square-free.  
Note also that Proposition \ref{prop:bi} and 
Assertion (5) of Lemma \ref{lemma:tri2} 
tells us that our sharper lower bound holds if at least one 
$\zeta_i$ is a degenerate root. So we may assume that $\zeta_1$ and 
$\zeta_2$ are both non-degenerate roots. Furthermore, letting 
$r\!:=\!\gcd(a_2,a_3)$, we can reduce to special case $r\!=\!1$ 
via the same argument as from the proof of Assertion (5) of 
Lemma \ref{lemma:tri2}. So we will also assume $\gcd(a_2,a_3)\!=\!1$. 

Our proof then follows almost exactly the format of the square-free case, with 
just two small changes: (a) We replace $f$ by the polynomial 
$F(x)\!:=\!\frac{f(x)}{(x-\tau)^2}$, where $\tau\!\in\!\Q$ is the unique 
degenerate root of $f$. (That $f$ has exactly one degenerate root, and 
it has multiplicity $2$, follows from Assertions (1) and (2) of Lemma 
\ref{lemma:tri2}.) (b) We replace Lemma \ref{lem:1} by a direct proof that 
$\norm{F(\tau)}\!\geq\!e^{-O(\log(dH))}$.  

To prove the last bound, observe that 
$F(\tau)\!=\!\frac{c_1}{(a_3-a_2)\tau^2}Q(1)$. Since 
$\ord_p\tau\!=\!\frac{\ord_p(a_2c_2/(a_3c_3))}{a_3-a_2}$, Assertion (3) 
of Lemma \ref{lemma:tri2} then tells us that\\ 
\mbox{}\hspace{2cm}$\ord_p F(\tau)\!\leq\!
\log_p(H)+\log_p(dH)+\log_pO(a^2_2a^3_3(a_3-a_2)^2)\!=\!O(\log_p(dH))$. \qed 

\section{Solving Trinomials over $\Q_p$} \label{sec:trinosolqp}
Unlike the binomial case (see Remark \ref{rem:depth}), the tree 
$\cT_{p,k}(f)$ can have depth $\Omega(\log_p(dH))$ or greater for a 
trinomial $f\!\in\!\Z[x]$ with $p\nmid f(0)$ and $k$ sufficiently large 
\cite{fgr}. However, Lemma \ref{lem:low} below 
will show that the structure of $\cT_{p,k}(f)$ is still simple: 
{\em No} path in $\cT_{p,k}(f)$ has more than $2$ vertices of out-degree 
more than $2$. Corollary \ref{cor:trinodepth} below  
will establish how large $k$ must be so that $\cT_{p,k}(f)$ is deep enough to 
encode (via Lemma \ref{lem:ulift}) all the non-degenerate roots of $f$ in 
$\Z_p$, {\em and} do so with sufficient accuracy for Newton iteration to 
converge quickly. Our estimates on $k$ will enable us to 
approximate all the roots of $f$ in $\Q_p$ in time 
$(p\log_p d)^3\log^{4+o(1)}(dH)$. Mild assumptions on the exponents of 
$f$ can also guarantee that the root node of $\cT_{p,k}(f)$ has $O(\sqrt{p})$ 
or even fewer children, and the presence of degenerate roots in $\Q^*_p$ for 
$f$ enables even tighter estimates for $k$. 
Each of these restrictions leads to speed-ups we will describe.

\subsection{Trees and Trinomials} 
\begin{lem} {\em \label{lem:low} 
Suppose $f(x)\!=\!c_1+c_2x^{a_2}+c_3x^{a_3}\!\in\!\Z[x]$ is a trinomial of
degree\linebreak
$d\!=\!a_3\!>\!a_2\!\geq\!1$, with all its coefficients having absolute value
at most $H$. 
Then every {\em non}-root  
nodal polynomial $f_{i,\zeta}$ of $\cT_{p,k}(f)$ with $\zeta\!\neq\!0$ 
mod $p$ satisfies $\deg \tf_{i,\zeta}\!\leq\!4$, $\deg 
\tf_{i,\zeta}\!\leq\!3$, or $\deg \tf_{i,\zeta}\!\leq\!2$, according as 
$p\!=\!2$, $p\!=\!3$, or $p\!\geq\!5$.}  
\end{lem} 
\begin{ex} 
{\em One can check that for $f(x)\!:=\!x^{10}+11x^2-12$, the tree 
$\cT_{2,8}(f)$ is isomorphic to \raisebox{-.15cm}
{\epsfig{file=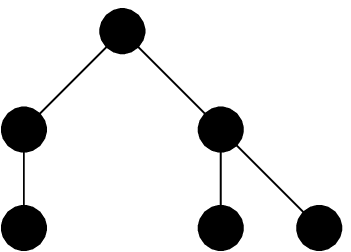,height=.5cm,clip=}}. 
In particular, this $f$ has exactly $6$ roots in $\Q^*_2$: 
$\tf_{2,2}\!=\!\tf_{2,1}\!=\!\tf_{2,3}\!=\!x^2+x$ and  
each of these (terminal) nodal polynomials has exactly $2$ 
non-degenerate roots in $\F_2$. Remembering the earlier digits encoded 
in $\cT_{2,8}(f)$, these $6$ roots then each lift to a unique root of $f$ in 
$\Z_2$. Note that $\tf_{1,1}(x)\!=\!x^4+x^2$ has degree $4$. \dia} 
\end{ex} 
\begin{ex}
{\em Composing Example \ref{ex:tri} with $x^2$, let us take
$f(x)\!:=\!x^{20}-10x^2+738$. One then sees that the tree
$\cT_{3,7}(f)$ is isomorphic to \raisebox{-.15cm}
{\epsfig{file=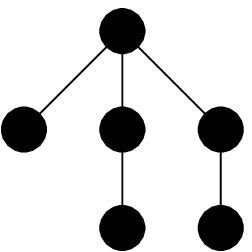,height=.5cm,clip=}}.
In particular, this $f$ has exactly $8$ roots in $\Q^*_3$, each arising
as a Hensel lift of a non-degenerate root in $\F_3$ of some nodal polynomial:
$\tf_{1,0}$, $\tf_{1,1}$, $\tf_{2,1}$, $\tf_{1,2}$, and $\tf_{2,8}$
respectively contribute $2$, $1$, $2$, $1$, and $2$ roots.
Note that $\tf_{1,2}(x)\!=\!x^3+2x^2+x$ has degree $3$. \dia}
\end{ex}

To prove Lemma \ref{lem:low} we will need a powerful result of Lenstra  
\cite{len99} on the Newton polygons of shifted sparse polynomials. First, let 
us define $d_m(r)$ to be the least common multiple of all integers that 
can be written as the product of at most $m$ pairwise distinct positive 
integers that are at most $r$, and set $d_m(r)\!:=\!1$ if $mr\!=\!0$. 
\begin{thm} \cite[Sec.\ 3]{len99}
\label{thm:lenstra} 
{\em Suppose $f\!\in\!\Q[x]$ is a $t$-nomial, $g(x)\!=\!f(1+px)$, and 
$r$ is the largest nonnegative integer such that
$r-\ord_p d_{t-1}(r)\!\leq\!\underset{0\leq j\leq t-1}{\max}\{j-\ord_p(j!)\}$. 
Then any lower edge of 
$\newt_p(g)$ with inner normal $(v,1)$ with 
$v\!\geq\!1$ lies in the strip $[0,r]\times \R$. \qed}   
\end{thm} 

\noindent 
We point out that the vector of parameters $(t,r,v)$ from our statement 
above would be $(k+1,m,\nu(x-1))$ in the notation of \cite{len99}, and the 
parameter $r$ there is set to $1$ in our application here. 

\medskip 
\noindent 
{\bf Proof of Lemma \ref{lem:low}:} First note that replacing 
$x$ by $cx$, for any $c\!\in\!\{1,\ldots,p-1\}$, preserves the number of roots 
of $f$ in $\Z_p$ and (up to relabelling the $\zeta$ in the subscripts of the  
$f_{i,\zeta}$) the tree $\cT_{p,k}(f)$. So 
to study $\tf_{1,\zeta_0}$ with $\zeta_0\!\in\!\{1,\ldots,p-1\}$,  
it suffices to study $\tf_{1,1}$. 

Note that the lower hull of any Newton polygon can be identified with a 
piecewise linear convex function on an interval. In particular, 
$f_{1,1}(x)\!=\!p^{-s(f,1)}f(1+px)$ and thus the lower hull of 
$\newt_p(f_{1,1})$ can be identified with the sum of the lower hull of 
$\newt_p(f(1+x))$ and the function $x-s(f,1)$. Note also that by the 
definition of $\newt_p$, the minimal $y$-coordinate of a point 
of $\newt_p(f(1+px))$ is exactly $s(f,1)$. 

Theorem \ref{thm:lenstra} then 
tells us that all lower edges of $\newt_p(f_{1,1})$ of non-positive slope 
lie in the strip $[0,r]\times \R$, where $r$ is the largest 
nonnegative integer such that 

\smallskip 
\noindent 
\mbox{}\hspace{.5cm} ($\star$) \hspace{4.5cm} 
$r-\ord_p d_2(r)\!\leq\!\eps_p$,

\smallskip 
\noindent 
where $\eps_2\!=\!1$ and $\eps_p\!=\!2$ for all $p\!\geq\!3$. In particular, 
the definition of $\newt_p(f_{1,1})$ tells us that $p$ divides the coefficient 
of $x^j$ in $f_{1,1}$ for all $j\!\geq\!r+1$ and thus 
$\deg \tf_{1,1}\!\leq\!r$. 

By Lemma \ref{lem:nodal}, all other non-root nodal 
polynomials $f_{i,\zeta}$ with $\zeta\!\neq\!0$ mod $p$ satisfy 
$\deg \tf_{i,\zeta}\!\leq\!\deg \tf_{1,1}$. So it 
suffices to prove that $r$ satisfies the stated bounds of our lemma. 
This is easily verified by first observing that $d_2(0)\!=\!d_2(1)\!=\!1$ 
and $d_2(2)\!=\!2$. So Inequality ($\star$) certainly holds for 
$r\!\in\!\{0,1,2\}$, regardless of $p$. Observing that $d_2(3)\!=\!6$ 
and $d_2(4)\!=\!24$, we then see that Inequality ($\star$) 
holds at $r\!=\!4$ (resp.\ $r\!=\!3$) when $p\!=\!2$ (resp.\ 
$p\!=\!3$). 

So it is enough to show that: 
(i) $r-\ord_2 d_2(r)\!\geq\!2$ for 
$r\!\geq\!5$, (ii) $r-\ord_3 d_2(r)\!\geq\!3$ for $r\!\geq\!4$, and 
(iii) $r-\ord_p d_2(r)\!\geq\!3$ for $r\!\geq\!3$ and 
$p\!\geq\!5$. 
From \cite[Prop.\ 2.4]{len99}, we have $\ord_p d_2(r)\!\leq\!\frac{2\log r}
{\log p}$. Note that, for any fixed $p$, the quantity $r-\frac{2\log r}
{\log p}$ is an increasing function of $r$ for $r\!\geq\!\frac{2}{\log p}$. 
Furthermore, $\ceil{7-\frac{2\log 7}{\log p}}\!\geq\!2$ 
for all $p\!\geq\!2$ and $\ceil{5-\frac{2\log 5}{\log p}}\!\geq\!3$
for all $p\!\geq\!3$. Noting that $d_2(5)\!=\!120$ and $d_2(6)\!=\!360$, it is 
then easily checked that (i)--(iii) all hold. \qed 
\begin{rem} 
\label{rem:low} 
{\em The proof of Lemma \ref{lem:binodepth} is simply the variation of the 
proof above where we replace Inequality ($\star$) by 
$r-\ord_p d_1(r)\!\leq\!1$, replace $d_2(r)$ with $d_1(r)$, and 
let $\eps_p\!=\!1$ for {\em all} $p$. In\linebreak 
\scalebox{.97}[1]{particular, the definition of $s(f,\zeta_0)$ tells us 
that $s(f,\zeta_0)\!\leq\!1+\ord_p f'(\zeta_0)\!=\!1+\ord_p d\!=\!1+\ell$.  
\dia}} 
\end{rem}   

It seems harder to get an upper bound on $s(f,\zeta_0)$ for trinomials 
than binomials. Nevertheless, we can derive a bound quadratic in $\log d$ and 
linear in $\log H$, and thereby estimate how large $k$ must be for our 
tree $\cT_{p,k}(f)$ to be deep enough for trinomial root approximation.   
\begin{cor} 
\label{cor:trinodepth} {\em 
Suppose $f(x)\!=\!c_1+c_2x^{a_2}+c_3x^{a_3}\!\in\!\Z[x]$ has degree 
$d$, $0\!<\!a_2\!<\!a_3$, $p\nmid c_1$, $c_2c_3\!\neq\!0$, and 
$|c_i|\!\leq\!H$ for all $i$. Let $r\!:=\!\gcd(a_2,a_3)$, 
define $S_0$ to be the maximum of $s(f,\zeta_0)$ (see 
Definition \ref{dfn:crazytree}) for any $\zeta_0\!\in\!\{1,\ldots,p-1\}$ 
satisfying $f(\zeta_0)\!=\!f'(\zeta_0)\!=\!0$ 
mod $p$, and set $S_0\!:=\!0$ should there be no such $\zeta_0$. Also let 
$D$ be the maximum of $\ord_p(\zeta-\xi)$ over 
all distinct non-degenerate roots $\zeta,\xi\!\in\!\Z_p$ of $f$ 
(if $f$ has at least $2$ non-degenerate roots in $\Z_p$) or $0$ 
(if $f$ has $1$ or fewer non-degenerate roots in $\Z_p$); and 
define $M_p$ to be $4$, $3$, or $2$, according as $p$ is $2$, $3$, 
or $\geq\!5$. Then: 

\smallskip 
\noindent 
1.\ $k\!\geq\!1+S_0\min\{1,D\}+M_p\max\{D-1,0\}  
\Longrightarrow$ the depth of $\cT_{p,k}(f)$ is at least $D$.\\  
2.\ $a_2\!=\!1\Longrightarrow S_0\leq 2+\ord_p(d(d-1)c_3/2)<2
+2\log_p(dH)$.\\   
3.\ $d\!\geq\!3\Longrightarrow$\\ 
\mbox{}\hspace{.6cm}\scalebox{.95}[1]{$S_0<2+2\ord_p(r)+ 
\log_p\left(\frac{d}{r} \left(\frac{d}{r}-1\right)H\right)
+\log(2)\log(4)2^{57/2}e^6 p\log\left(\frac{d}{r}-1\right)\log_p
\left(\frac{d}{r}\left(\frac{d}{r}-1\right)H\right)$}\\ 
\mbox{}\hspace{1cm}$<2+\log_p\left(\frac{d}{r}\left(\frac{d}{r}-1\right)H
\right)+147164373392p\log\left(\frac{d}{r}-1\right)
\log_p\left(\frac{d}{r}\left(\frac{d}{r}-1\right)H\right)$.\\  
4.\ $f$ has a degenerate root in $\C_p \Longrightarrow 
S_0\!\leq\!2+2\log_p(r)+\log_p(d/r)$.\\  
5.\ The lower bound for $k$ from Assertion (1) can be attained for 
$k\!=\!O(p\log^2_p(dH)\log d)$ or\\ 
\mbox{}\hspace{.5cm}$k\!=\!O(\log_p(dH))$, according 
as $f$ has no degenerate roots in $\C_p$, or at least one such  
root. } 
\end{cor} 
\begin{rem}
\label{rem:extremal} {\em 
Note that $d\!\geq\!2$ for any trinomial, and $d\!=\!2$ implies 
$a_2\!=\!1$ above. \dia } 
\end{rem} 

\noindent 
{\bf Proof of Corollary \ref{cor:trinodepth}:}\\  
{\bf Assertion (1):} $\cT_{p,k}(f)$ 
always includes a root node by definition, so the case $D\!=\!0$ is 
trivial and we assume $D\!\geq\!1$. 

Our lower bound on $k$ then follows easily from Lemma \ref{lem:low}:  
Since $f$ has distinct non-degenerate roots $\zeta,\xi\!\in\!\Z_p$ 
with $\ord(\zeta-\xi)\!\geq\!1$ by assumption, this means that 
$\zeta\!=\!\xi$ mod $p$ and thus $\tf$ must have a degenerate 
root $\zeta'_0\!\in\!\{1,\ldots,p-1\}$ (since $p\nmid c_1$). 
Having $k\!\geq\!1+S_0$ then simply allows the 
root node to have maximally many child nodes (and thus depth $\!\geq\!1$), 
thanks to Definition \ref{dfn:crazytree}. Furthermore, thanks 
to Lemma \ref{lem:low}, the summand $M_p\max\{D-1,0\}$ simply guarantees that 
$\cT_{p,k}$ has depth $D$ and that $\cT_{p,k}(f)$ has maximally many nodes 
at depth $\leq\!D$. (Note that for any nodal polynomial $f_{i,\zeta'}$ with 
$i\!\geq\!1$, we have that $s(f_{i,\zeta'},\zeta_i)$ is bounded from 
above by $4$, $3$, or $2$, according as $p$ is $2$, $3$, or $\geq\!5$, 
thanks to Lemma \ref{lem:nodal}.) In particular, we\linebreak 
\scalebox{.9}[1]{see that any $k$ satisfying 
our lower bound yields a $k$ satisfying all the assumptions of 
Lemma \ref{lem:ulift}. \qed}  

\smallskip 
\noindent 
{\bf Assertion (2):} Immediate from $s(f,\zeta_0)
\!\leq\!2+\ord_p\frac{f''(\zeta_0)}{2}$ (thanks to the definition of 
$s(\cdot,\cdot)$ as a minimum), $f''(\zeta_0)\!=\!d(d-1)c_3\zeta^{d-2}_0$, 
and $\ord_p \zeta_0\!=\!0$. \qed 

\medskip 
\noindent 
{\bf Note.} {\em \fbox{We now temporarily assume that $\gcd(a_2,a_3)\!=\!1$}, 
to simplify the proofs of Assertions\linebreak 
\scalebox{.97}[1]{(3) and (4), and show later how to reduce the 
case $\gcd(a_2,a_3)\!>\!1$ to the case $\gcd(a_2,a_3)\!=\!1$. \dia}}  

\medskip 
\noindent 
{\bf Assertion (3):} \scalebox{.95}[1]{First note that we must have 
$p\nmid c_2$ or $p\nmid c_3$ in order for $\tf$ to have a root in $\F^*_p$.}  

Since $f'(\zeta_0)\!=\!a_2c_2\zeta^{a_2-1}_0+a_3c_3\zeta^{a_3-1}_0\!=\!0$ 
mod $p$, and $\gcd(a_2,a_3)\!=\!1$, 
we see that $p|a_2 \Longrightarrow$\linebreak 
\scalebox{.95}[1]{$\ord_p c_3\!=\!\ord_p a_2\!>\!0$ and 
$p\nmid a_3c_2$. In which case, $\ord_p f'(\zeta_0)\!=\!\ord_p(a_2c_2)
+\ord_p\left(1-\frac{-a_3c_3}{a_2c_2}\zeta^{a_3-a_2}_0\right)$,}\linebreak 
and then we can bound $\ord_p f'(\zeta_0)$ from above by the $n\!=\!2$ 
case of Yu's Theorem if the second valuation is {\em not} $\infty$. Should 
this valuation be $\infty$, then we can instead apply the $n\!=\!2$ 
case of Yu's Theorem to $\ord_p f''(\zeta_0)\!=\!\ord_p(a_2(a_2-1)c_2)
+\ord_p\left(1-\frac{-a_3(a_3-1)c_3}{a_2(a_2-1)c_2}\zeta^{a_3-a_2}_0\right)$, 
since $\frac{a_3-1}{a_2-1}\!\neq\!1$. So we obtain our stated bound 
directly from Theorem \ref{thm:yu}, and the fact that $s(f,\zeta_0)\!\leq\!
\min\{1+\ord_p f'(\zeta_0),2+\ord_p f''(\zeta_0)\}$.   

Similarly, $p|a_3 \Longrightarrow \ord_p c_2\!=\!\ord_p a_3\!>\!0$ and 
$p\nmid a_2c_3$. In which case, $\ord_p f'(\zeta_0)\!=\!\ord_p(a_3c_3)
+\ord_p\left(1-\frac{-a_2c_2}{a_3c_3}\zeta^{a_2-a_3}_0\right)$, and 
we proceed in the same way as the last paragraph to obtain our stated bound. 

So let us now assume $p\nmid a_2a_3$. Then $f'(\zeta_0)\!=\!0$ mod $p 
\Longrightarrow p\nmid c_2c_3$, since $\ord_p\zeta_0\!=\!0$ and 
$p$ can not divide both $c_2$ and $c_3$. So then we again attain 
our bound just like in the last paragraph. \qed  

\smallskip 
\noindent 
{\bf Assertion (4):} 
Note that $p\nmid c_1$ implies that any 
degenerate root $\tau\!\in\!\C_p$ of $f$ must be nonzero.  
Lemma \ref{lemma:tri2} then tells us that $\tau$ is the only degenerate 
root of $f$ in $\C_p$ and $\tau\!\in\!\Q^*_p$. Moreover, from the proof of 
Lemma \ref{lemma:tri2}, we have 
$f(\tau x)\!=\!\frac{c_1}{(a_3-a_2)\tau^2}q(x)$ identically and 
$\ord_p \tau\!=\!\frac{\ord_p(a_2c_2) - \ord_p(a_3c_3)}{a_3-a_2}$. 
(Recall that $q(x)\!=\!(a_3-a_2)-a_3x^{a_2}+a_2x^{a_3}$ has $1$ as its 
unique degenerate root in $\C_p$.) 

Now, we must have $p\nmid c_2$ or $p\nmid c_3$ in order for there to be any 
roots at all for $\tf$. 

\smallskip
\noindent 
{\bf Sub-Case $p\nmid c_2$.} If $\tau$ has negative valuation, then we must 
have $p|c_3$ by Theorem \ref{thm:newt}. 
Also, $f'(\zeta_0)\!=\!\zeta^{a_2-1}_0(c_2a_2+c_3a_3
\zeta^{a_3-a_2}_0)\!=\!0$ mod $p \Longrightarrow p|a_2$ since $p\nmid c_2$. 
Since $\ord_p \tau\!=\!\frac{\ord_p(a_2c_2) - \ord_p(a_3c_3)}
{a_3-a_2}\!<\!0$ by assumption, we must have 
$\ord_p (a_3c_3)\!>\!\ord_p(a_2c_2)$ and thus 
\linebreak 
$\ord_p f'(\zeta_0)\!=\!\ord_p(c_2a_2)\!=\!\ord_p(a_2)$.  
In other words, $\ord_p \tau\!<\!0 \Longrightarrow S_0\!\leq\!1+\ord_p(a_2)$. 

So let us now assume $\ord_p\tau\!=\!0$. Then by our identity 
$f(\tau x)\!=\!\frac{c_1}{(a_3-a_2)\tau^2}q(x)$, and the fact that 
$\tau\!\in\!\Q^*$ (via Assertion (1) of Lemma \ref{lemma:tri2}), 
the vector of coefficient valuations for $f$ and the 
vector of coefficient valuations for $q$ differ by a multiple of $(1,1,1)$. 
So our assumptions that $p\nmid c_1c_2$ and $\gcd(a_2,a_3)\!=\!1$ 
imply that $p\nmid (a_3-a_2)a_3$. So then, $a_3-a_2$ is invertible mod $p$ 
and, by the rescaling between $f$ and $q$, we have that $\tf$ and $\tq$ share 
the same value of $S_0$ (as well as the same number of degenerate roots 
in $\{1,\ldots,p-1\}$). So let us now work with $q$ instead, and assume 
for the remainder of this sub-case that $\zeta_0$ is a degenerate root of 
$\tq$ mod $p$.   

If $p|a_2$ then 
$\ord_p q'(\zeta_0)\!=\!\ord_p(a_2)+\ord_p\left(-1+\zeta^{a_3-a_2}_0\right)$ 
(since $p\nmid a_3$). Also, $\ord_p q''(\zeta_0)\!=\!\ord_p(a_2)
+\ord_p(-a_2+a_3\zeta^{a_3-a_2}_0
-(-1+\zeta^{a_3-a_2}_0))$. Since $p|a_2$ and $p\nmid a_3$, we see 
that\linebreak
$\ord_p (-1+\zeta^{a_3-a_2})\!>\!0$ implies that 
$\ord_p q''(\zeta_0)\!=\!\ord_p a_2$. 
On the other hand, if\linebreak 
$\ord_p (-1+\zeta^{a_3-a_2})\!=\!0$, then 
$\ord_p q'(\zeta_0)\!=\!\ord_p a_2$ from our earlier formula for 
$\ord_p q'(\zeta_0)$. So by the definition of $s(\cdot,\cdot)$, we 
obtain $S_0\!\leq\!2+\ord_p a_2$. 

To conclude, $p\nmid a_2$, combined with our earlier conclusion 
that $p\nmid (a_3-a_2)a_3$, implies that $\zeta_0\!=\!1$, thanks to 
Assertion (1) of Lemma \ref{lemma:tri2}. In which case, $q'(1)\!=\!0$ 
but $q''(1)\!=\!a_2a_3((a_3-1)-(a_2-1))\!=\!a_2a_3(a_3-a_2)$ and thus 
$S_0\!\leq\!2$.  

\smallskip
\noindent 
{\bf Sub-Case $p\nmid c_3$.} Here, we must have $\ord_p \tau\!=\!0$ and 
thus $\ord_p(a_2c_2)\!=\!\ord_p(a_3c_3)$ by our earlier formula 
for $\ord_p\tau$. In particular, we must have $\ord_p(a_2c_2)\!=\!\ord_p a_3$ 
since $p\nmid c_3$. Note also that $p|a_2$ thus implies $p|a_3$, which 
would contradict $\gcd(a_2,a_3)\!=\!1$. So we must also have $p\nmid a_2$ 
and thus $\ord_p c_2\!=\!\ord_p a_3$.  
Since we already proved the Sub-Case $p\nmid c_2$, let us now assume 
$p|c_2$ (and thus $p|a_3$).  

By our identity $f(\tau x)\!=\!\frac{c_1}{(a_3-a_2)\tau^2}q(x)$, and the fact 
that $\tau\!\in\!\Q^*$ (via Assertion (1) of Lemma \ref{lemma:tri2}), 
the vector of coefficient valuations for $f$ and the 
vector of coefficient valuations for $q$ differ by a multiple of $(1,1,1)$. 
So our assumptions that $p\nmid c_1c_3$ and $\gcd(a_2,a_3)\!=\!1$ 
imply that $p\nmid (a_3-a_2)a_2$. So then, $a_3-a_2$ is invertible mod $p$ 
and, by the rescaling between $f$ and $q$, we have that $\tf$ and $\tq$ share 
the same value of $S_0$ (as well as the same number of degenerate roots 
in $\{1,\ldots,p-1\}$). So let us now work with $q$ instead, and assume 
now that $\zeta_0$ is a degenerate root of $\tq$ mod $p$.   

Observe then that $\ord_p q'(\zeta_0)\!=\!\ord_p(a_3)
+\ord_p\left(-1+\zeta^{a_3-a_2}_0\right)$
(since $p\nmid a_2$). Also,\linebreak 
$\ord_p q''(\zeta_0)\!=\!\ord_p(a_3)
+\ord_p(-a_2+a_3\zeta^{a_3-a_2}_0
-(-1+\zeta^{a_3-a_2}_0))$. Since $p|a_3$ and $p\nmid a_2$, we see that $\ord_p 
(-1+\zeta^{a_3-a_2})\!>\!0$ implies that $\ord_p q''(\zeta_0)\!=\!\ord_p a_3$.
On the other hand, if $\ord_p (-1+\zeta^{a_3-a_2})\!=\!0$, then 
$\ord_p q'(\zeta_0)\!=\!\ord_p a_3$ from our earlier formula for
$\ord_p q'(\zeta_0)$. So by the definition of $s(\cdot,\cdot)$, we
obtain $S_0\!\leq\!2+\ord_p a_3$. \qed 

\smallskip 
\noindent 
{\bf Extending to $\pmb{\gcd(a_2,a_3)\!>\!1}$.} 
To complete our proofs of Assertions (3) and (4) let us assume 
$r\!:=\!\gcd(a_2,a_3)\!>\!1$ and recall that $\barf$ is the unique 
polynomial in $\Z[x]$ satisfying $f(x)\!=\!\barf(x^r)$ identically. 
Clearly then, $\deg \barf\!=\!\frac{\deg f}{r}$ and 
any root $\tau\!\in\!\C_p$ of $f$ induces a root  
$\tau^r$ of $\barf$. Furthermore, $\tf$ having a degenerate root 
$\zeta_0\!\in\!\{1,\ldots,p-1\}$ clearly implies that the mod $p$ 
reduction of $\barf$ has $\mu_0$ as a degenerate root, where 
$\mu_0\!\in\!\{1,\ldots,p-1\}$ 
is the mod $p$ reduction of $\zeta^r_0$.
The Chain Rule then implies $\ord_p f'(\zeta_0)\!=\!
\ord_p(r)+\ord_p \barf'(\mu_0)\!\leq\!\log_p(d)+\ord_p\barf'(\mu_0)$. 

Should $f'(\zeta_0)$ vanish identically, then Assertion (1) 
of Lemma \ref{lemma:tri2} easily implies that all the degenerate roots 
of $f$ have multiplicity $2$ and thus $f''(\zeta_0)$ can not vanish. 
In which case, via the Chain Rule again, 
$\ord_p f''(\zeta_0)\!=\! 2\ord_p(r) + \ord_p \barf''(\mu_0)
\!\leq\!2\log_p(d)+\ord_p \barf''(\mu_0)$. So our general formula 
follows immediately from the case $r\!=\!1$, which we've already proved. \qed 

\smallskip 
\noindent 
{\bf Assertion (5):} Immediate from Assertions (3) and (4), and Theorem 
\ref{thm:tri}. \qed   

\subsection{Building Trees Efficiently} 
It is easy to see that the only degenerate root the quadratic 
trinomial $c_1+c_2x+x^2\!\in\!\Z[x]$ 
can have mod $p$ is exactly $-c_2/2$ when $p\!\geq\!3$. 
(For $p\!=\!2$ it is clear that the only monic degenerate quadratics are 
$x^2+1$ and $x^2$, with respective degenerate roots $1$ and $0$.) 
It will be useful to have a similar statement for trinomials with  
$(p,d)\!\in\!\{2,3\}\times\{3,4\}$. 
\begin{prop} 
\label{prop:degen} {\em 
Suppose $f(x)\!=\!c_0+c_1x+c_2x^2+c_3x^3+c_4x^4\!\in\!\Z[x]$ has degree 
$d\!\geq\!2$, and $|c_i|\!\leq\!H$ for all $i$. Then:\\  
0. The discriminant of $f$ can be evaluated in time $O(\log(\max\{p,H\})
\log\log\max\{p,H\})$.\\  
1. When $p\!\leq\!3$ we can find all the degenerate roots of $f$ in $\F_p$ (or 
correctly declare there\linebreak 
\mbox{}\hspace{.6cm}none) in time $O(\log H)$. In particular, 
$f$ has at most $1$ (resp.\ $2$) degenerate root(s) in $\F_p$,\linebreak
\mbox{}\hspace{.6cm}according as $d\!\leq\!3$ or $d\!=\!4$.\\ 
2. For any prime $p$ we can find all the non-degenerate roots of $f$ (or 
correctly declare there\linebreak 
\mbox{}\hspace{.6cm}are none) in deterministic time $O(p^{1/2}\log^2p)$.}   
\end{prop} 

\noindent 
{\bf Proof:} Assertion (0) follows from the definitions of the quartic, cubic, 
and quadratic discriminants (see, e.g., \cite[Ch.\ 12]{gkz94}), 
Theorem \ref{thm:cxity}, and the fact that evaluating 
$\Delta_d(f)$ reduces to evaluating  
a $7\times 7$, $5\times 5$, or $3\times 3$ determinant in the coefficients of 
$f$ (followed by division by the leading coefficient of $f$), after 
reducing the coefficients mod $p$. 

For Assertion (1), first note that $p\!\leq\!3$ implies that we can reduce 
the coefficients of $f$ and $f'$ mod $p$ in time $O(\log H)$ thanks to Theorem 
\ref{thm:cxity}. We can then simply use brute-force (over a search space 
with at most $3$ elements!) to find all the 
degenerate roots of $f$ in time $O(1)$. 
In particular, since any degenerate root must have multiplicity $\geq\!2$, 
the only way $f$ can have more than $1$ degenerate root is for $d\!=\!4$, 
in which case there can be no more than $2$ degenerate roots. For instance, 
$x^4+x^2+1$ (resp.\ $x^4+x^2$) has degenerate roots 
$\{\pm 1\}\!\in\!\F_3$ (resp.\ $\{0,1\}\!\in\!\F_2$). 

Assertion (2) follows immediately from Shoup's deterministic algorithm 
for factoring arbitrary univariate polynomials over a finite field 
\cite{shoup}, upon specializing to degree $\leq\!4$. \qed 

\begin{lem} {\em \label{lem:cxity}
For any trinomial $f(x)\!=\!c_1+c_2x^{a_2}+c_3x^{a_3}\!\in\!\Z[x]$ of degree 
$d$, with $p\nmid c_1$, $0\!<\!a_2\!<\!a_3$, and $|c_i|\leq\!H$ for all $i$, 
let $\nu$ denote the number of degenerate roots of $\tf$ in $\F^*_p$ and 
let $\cD$ denote the depth of $\cT_{p,k}(f)$. Then 
$\cT_{p,k}(f)$ has $\leq\!1+\left(2\cD-1\right)\nu$ 
nodes; and we can 
compute the mod $p$ reductions of all the nodal polynomials 
$f_{i,\zeta}$ of $\cT_{p,k}(f)$, as well as all the values of the 
$s(f_{i-1,\mu},\zeta_{i-1})$, in deterministic time\\  
$O\!\left((p+\log d)\log(dp)\log(\log(dp))+p\log^2(p)\log\log(p)\right.$\\ 
\mbox{}\hfill$+\left. \nu \cD [k\log(p)\log(k\log p)\log(d) 
+\log H]+\log(H)\log(dpH)\log\log(dpH) \right)$. }  
\end{lem}

\noindent 
{\bf Proof:} By Lemma \ref{lem:low}, all non-root nodal 
polynomials have mod $p$ reduction of degree no greater than $4$. Thus, the 
root node of $\cT_{p,k}(f)$ has $\leq\!\nu$ ($\leq\!p-1$) children, 
and any node at depth $\geq\!1$ has no more than $2$ children (since a 
polynomial of degree $\leq\!4$ has $\leq\!2$ degenerate roots). Lemma 
\ref{lem:nodal} also tells us that $\deg \tf_{i,\mu+\zeta_{i-1}p^{i-1}}$ is at 
most the multiplicity of $\zeta_{i-1}\!\in\!\F^*_p$ as a root of 
$\tf_{i-1,\mu}$. So any node $v$ that has an ancestor at level $\geq\!1$ with 
$2$ children can have no more than $1$ child. Thus, there can be no more than 
$2\nu$ nodes at depth $i\!\geq\!2$. It is then clear that $\cT_{p,k}(f)$ has 
at most $1+\left(2\cD -1\right)\nu$ nodes. 

We now check whether $\tf$ has any degenerate roots in $\F_p$: 
By assumption, they must lie in $\F^*_p$. Also, should $p|c_3$, $\tf$ would 
be a binomial and thus have degenerate roots in $\F^*_p$ only if $p|a_2$;  
in which case any root of $\tf$ in $\F^*_p$ is degenerate.   
We can then decide if there are degenerate roots simply by checking whether 
$(-c_1/c_2)^{(p-1)/\gcd(a_2,p-1)}\!=\!1$ mod $p$, which can 
be done in time $O(\log(dH)\log(\log(dH))+\log^2(p)\log\log p)$ via 
Theorem \ref{thm:cxity}. Should there be any degenerate roots, there 
will then be exactly $\gcd(a_2,p-1)$ many, 
and we can then find them in time no worse than 
$O((p+\log d)\log(dp)\log(\log(dp))+\log(H)\log(pH)\log\log(pH))$ 
via brute-force (much like our earlier complexity analysis of Steps 
5--7 of Algorithm \ref{algor:binoqp}). 

So let us assume $p\nmid c_3$. Note that $p|\gcd(a_2,a_3) \Longrightarrow$ 
every root of $\tf$ in $\F^*_p$ is degenerate, in which case we 
can simply find all these roots first by reducing the coefficients 
(resp.\ exponents) of $\tf$ mod $p$ (resp.\ mod $p-1$) in 
time\\ 
\mbox{}\hspace{1.5cm}$O(\log(\max\{d,p\})\log(\log\max\{d,p\})+
\log(\max\{p,H\})\log\log\max\{p,H\})$\\ 
and then applying brute-force search in time $O(p\log^2(p)\log\log p)$. So let 
us assume\linebreak 
$p\nmid \gcd(a_2,a_3)$. Observe then that  
$\tf$ has degenerate roots in $\F^*_p \Longleftrightarrow 
\Delta_{\mathrm{tri}}(\tf)\!=\!0$ mod $p$, thanks to Assertion (1) of Lemma 
\ref{lemma:tri}. In particular, by Theorem \ref{thm:cxity}, 
$\Delta_{\mathrm{tri}}(\tf)$ can be computed mod $p$ in time 
$O(\log(\max\{d,p\})\log(\log\max\{d,p\})+
\log(\max\{H,p\})\log\log\max\{H,p\})$ 
(to reduce the exponents of $\Delta_{\mathrm{tri}}(\tf)$ mod 
$p-1$ and the power bases mod $p$) plus $O(\log^2(p)\log\log p)$ 
to compute the monomials of $\Delta_{\mathrm{tri}}(\tf)$.
If $\Delta_{\mathrm{tri}}(\tf)\!\neq\!0$ mod $p$ then we know 
$\tf$ has no degenerate roots and then $\cT_{p,k}(f)$ is simply a single 
root node. Otherwise, let $r'\!:=\!\gcd(a_2,a_3,p-1)$ and apply the 
Extended Euclidean Algorithm (in time $O(\log(p)\log^2 \log p)$ via 
Theorem \ref{thm:cxity}) to $a_2$ mod $p-1$ and $a_3$ mod $p-1$ to find 
$\alpha,\beta\!\in\!\Z$ with logarithmic height 
$O(\log p)$ such that 
$\alpha(a_2 \text{ mod } p-1)+\beta(a_3 \text{ mod } 
p-1)\!=\!r'$. 
Assertion (1) of Lemma \ref{lemma:tri2} then tells us that the 
degenerate roots of $\tf$ in $\F^*_p$ are exactly the roots of 
$g(x)\!:=\!x^{r'}-(-1)^\alpha\left(\frac{c_1}{a_3-a_2}\right)^{\alpha+\beta}
\left(\frac{a_3}{c_2}\right)^\alpha
\left(\frac{a_2}{c_3}\right)^\beta$ in $\F^*_p$. 
Lemmata \ref{lem:binoqp} and \ref{lem:binoq2} and Theorem \ref{thm:cxity} then 
easily imply that deciding whether $g$ has any roots in 
$\F^*_p$ takes time $O(\log^2(p)\log\log p)$, and there are exactly 
$r'$ many degenerate roots in $\F^*_p$ if so. Just as in the last 
paragraph, we can then apply brute-force to $g$ in time\\  
\mbox{}\hspace{2cm}$O((p+\log d)\log(dp)\log(\log(dp))
+\log(H)\log(pH)\log\log(pH))$\\  
to find all the degenerate roots of $\tf$ in $\F^*_p$. 

Assuming $\tf$ has degenerate roots in $\F^*_p$, let us now see how to compute 
the child nodes of the root node in $\cT_{p,k}(f)$: 
First note that the coefficient of $x^i$ in the monomial term expansion of 
$c(\mu+px)^a$ mod $p^j$ (for $i\!\leq\!j$) is simply $c\binom{a}{i}
\mu^{a-i}p^i$ mod $p^j$. Also, Lemma \ref{lem:nodal} tells us that 
$f_{i,\zeta}(x)\!=\!p^{-s}f(\mu+p^i x)$ mod $p^j$ for suitable $(s,\mu,j)$. 
Putting this together, this means we can compute $s(f,\zeta_0)$ and 
$\tf_{1,\zeta_0}$ (for all degenerate roots $\zeta_0\!\in\!\F^*_p$ of $\tf$) 
by evaluating $\zeta^{a_2}_0$ and $\zeta^{a_3}_0$ mod $p^k$, 
$\binom{a_2}{i}$ and $\binom{a_3}{i}$ for $i\!\in\!\{0,1,2\}$ if $p\!\geq\!5$,  
and $O(1)$ additional ring operations in $\Z/(p^k)$. (We instead take 
$i\!\in\!\{0,1,2,3\}$ or $\{0,1,2,3,4\}$ according as $p$ is $3$ or $2$.) 
Via Recursive Squaring (a.k.a.\ the Binary Method \cite[pp.\ 102--103]{bs}), 
Theorem \ref{thm:cxity} tells us that we can compute the $a_2$nd and $a_3$rd 
powers of all the degenerate roots $\zeta_0\!\in\!\F^*_p$ in time 
$O(v\cdot \log(d) \cdot k\log(p)\log(k\log p))$, and 
the remaining operations are negligible in comparison. In particular, 
each $s(f,\zeta_0)$ can be computed by bisection and the resulting complexity 
is also 
negligible compared to the preceding $O$-estimate.  

So in summary, all computations necessary to find all child nodes of 
the root node take time no greater than 

\smallskip 
\noindent 
$\displaystyle{O((p\log(p)+\log d)\log(dp) \log\log(dp) 
 +p\log^2(p)\log(\log p)}$\\  
\mbox{}\hspace{4cm}
$\displaystyle{+\log(H)\log(dpH)\log\log(dpH) + \nu k \log(d)\log(p)
 \log(k\log p))}$.  

\smallskip 
Having computed all the mod $p$ reductions of the nodal polynomials 
$\tf_{1,\zeta_0}$ at depth $1$, we then proceed inductively, performing 
almost the same calculations as in the last two paragraphs. The only 
difference, assuming $p\!\geq\!5$, is then applying applying the quadratic 
discriminant (instead of the trinomial discriminant) to detect and find the 
{\em sole} degenerate root of $f_{i-1,\mu}$ (for $i\!\in\!\{2,
\ldots,k-1\}$), should there be one. (Should $p\!\in\!\{2,3\}$ then we 
simply apply Proposition \ref{prop:degen} instead, and possibly have 
two degenerate roots in the worst case when $p\!=\!2$.) This eliminates the 
need for brute-force search, and gives us an improved complexity bound 
of $O(k\log(p)\log(k\log p) \log(d) +\log H)$ 
to compute the children (no more than two) of each $f_{i-1,\mu}$. 

Summing all the resulting complexity estimates over all $O(\nu \cD)$ children, 
we are done. \qed 

\medskip 
\begin{cor} {\em 
\label{cor:average} 
Following the notation of Lemma \ref{lem:cxity}, we have the following 
improved\linebreak 
complexity bounds for computing the mod $p$ reductions of all the 
nodal polynomials of $\cT_{p,k}(f)$ and their respective $s(\cdot,\cdot)$ 
values: 

\smallskip 
\noindent 
1.\ If we only wish to construct the sub-tree of $\cT_{p,k}(f)$ corresponding 
to $\zeta_0\!=\!1$, and correctly\\ 
\mbox{}\hspace{.5cm}declare whether $1$ is a degenerate 
root of $f$:\\  
\mbox{}\hspace{.5cm}Deterministic time\\  
\mbox{}\hspace{.5cm}$O\!\left(\cD[k\log(p)\log(k\log p)\log(d)
+\log H]+\log^2(p)\log\log(p)\right.$\\  
\mbox{}\hspace{2.8cm}\scalebox{1}[1]{$\left.
+\log(\max\{d,p\})\log(\log\max\{d,p\})
+ \log(\max\{p,H\})\log(\log\max\{p,H\}) \right)$.} 

\smallskip 
\noindent 
2.\ If the exponents are $\{0,a_2,a_3\}$ with 
$\gcd(a_2a_3(a_3-a_2),(p-1)p)\!\leq\!2$:\\ 
\mbox{}\hspace{.5cm}Deterministic time\\ 
\mbox{}\hspace{.5cm}$O\!\left(p^{1/2}\log^2(p)+\cD[k\log(p)\log(k\log p)\log(d)
+\log H]\right.$\\
\mbox{}\hspace{2.8cm}$\left.
+\log(\max\{d,p\})\log(\log\max\{d,p\})
+ \log(\max\{p,H\})\log(\log\max\{p,H\}) \right)$,\\
\mbox{}\hspace{.5cm}or Las Vegas randomized time\\ 
\mbox{}\hspace{.5cm}$O(\cD[k\log(p)\log(k\log p)\log(d)
+\log H]+\log^{2+o(1)}(p)$\\
\mbox{}\hspace{1.5cm}$+\log(\max\{d,p\})\log(\log\max\{d,p\})
+ \log(\max\{p,H\})\log(\log\max\{p,H\}) 
 \vphantom{\cD^{\cD^\cD}}
)$.} 
\end{cor} 

\medskip 
\begin{rem} 
\label{rem:bottle} 
{\em While we state a randomized speed-up in Assertion (2) above, any 
asymptotic gains are unfortunately overwhelmed once we insert our upper bounds 
on $k$ and $\cD$ for the non-degenerate case from Corollary 
\ref{cor:trinodepth} and Theorem \ref{thm:tri}. Nevertheless, we state our 
bounds in a refined way above, should better bounds on $k$ and $\cD$ become 
available in the future. \dia } 
\end{rem} 

\noindent 
{\bf Proof of Corollary \ref{cor:average}:} In what follows, we keep in mind 
the template of the proof of Lemma \ref{lem:cxity}, and simply point out the 
key changes resulting in speed-ups. 

\smallskip 
\noindent 
{\bf Assertion (1):} Here there is no need to search for roots of 
$\tf$: We merely evaluate $\tf$ and $\tf'$ at $1$ to see if $1$ 
is a degenerate root. This amounts to time

\smallskip 
\noindent 
\mbox{}\hspace{1cm}$O(\log(\max\{d,p\})\log(\log\max\{d,p\})
+ \log(\max\{p,H\})\log(\log\max\{p,H\}))$

\smallskip 
\noindent 
to reduce exponents mod 
$p-1$ and coefficients  mod $p$, and then time    
$O(\log^2(p)\log\log(p))$ for the evaluation. At this point, we 
also know if $1$ fails to be a degenerate root of $\tf$. 

We then need time $O(\max\{k\log p,\log H\}\log\max\{k\log p,\log H\})$ 
to reduce the coefficients of $\tf$ mod $p^k$, and then time 
$O(\log(d) k\log(p)\log(k\log p))$ to compute $s(f,1)$ and the child 
node of the root node. For the remaining descendants, Lemma 
\ref{lem:low} tells us that there are at most $2$ children, and 
any subsequent siblings can have no further offspring with more than 
one child. Also, as observed earlier, we can find the degenerate roots 
of the mod $p$ reduction of any non-root nodal polynomial in time 
$O(\log H)$. So the remaining child nodes take time  
$\cD-1$ times $O(k\log(d)\log(p)\log(k\log p)+\log(H))$ to compute. \qed 

\smallskip 
\noindent 
{\bf Assertion (2):} 
The gcd assumption on the exponents implies there can be 
at most $2$ degenerate roots for $\tf$ in $\F_p$ (and they are 
nonzero since we originally assumed $p\nmid c_1$ in Lemma \ref{lem:cxity}): 
This follows 
from basic group theory if $p|c_3$ and via Lemma \ref{lemma:tri2} if 
$p\nmid c_3$. 

If $p|c_3$ then we can decide whether $\tf$ has a degenerate root in $\F^*_p$ 
by computing $g_1\!:=\!\gcd(\tf,x^{p-1}-1)$ and checking whether 
$\deg g_1\!\geq\!1$ or not: If $\deg g_1\!=\!1$ then we can easily find the 
unique root of $g_1$ using one arithmetic operation in $\F_p$.  
If $\deg g_1\!=\!2$ then we can find the roots either in 
deterministic time $O(p^{1/2}\log^2 p)$ via Shoup's fast deterministic 
factoring algorithm \cite{shoup},  
or Las Vegas time $\log^{2+o(1)}p$ via the fast randomized factorization 
algorithm of Kedlaya-Umans \cite{ku08}. Furthermore, $g_1$ can be computed 
efficiently by first computing $x^{a_2}$ mod $x^{p-1}-1$ via Recursive Squaring 
(a.k.a.\ the Binary Method \cite[pp.\ 102--103]{bs}), and then 
computing the rest of $\tf$ mod $x^{p-1}-1$. This entails $O(\log d)$ 
reductions (of exponents) mod $p-1$, along with $3$ arithmetic operations 
in $\F_p$, meaning additional (deterministic) time 
$O(\log(d)\log(\max\{d,p\}) \log\log\max\{d,p\})$ 
via Theorem \ref{thm:cxity}. 

\scalebox{.94}[1]{If $p\nmid c_3$ then we can decide whether $\tf$ has a 
degenerate root in $\F^*_p$ by first checking 
$\Delta_{\mathrm{tri}}(\tf)\!\stackrel{?}{=}\!0$}\linebreak   
\scalebox{.93}[1]{mod $p$, which takes time $O(\log(\max\{d,p\}) 
\log(\log\max\{d,p\}) +\log(\max\{H,p\})\log\log\max\{H,p\})$}\linebreak 
(as already observed in our last 
proof). If this discriminant indeed vanishes mod $p$ then we compute  
$g_2\!:=\!\gcd(\tf,\tf')\!=\!\gcd(\tf,\tf'/x^{a_2-1})$. Like $g_1$, the 
polynomial $g_2$ has degree $\leq\!2$, and it can be computed efficiently, 
along with its roots (if any) in deterministic time 
\[O(p^{1/2}\log^2(p)+\log(d)\log(\max\{d,p\})\log\log\max\{d,p\}),\]   
or Las Vegas time
\[ O(\log^{2+o(1)}(p)+\log(d)\log(\max\{d,p\})\log\log\max\{d,p\}).\]  

\scalebox{.98}[1]{We then proceed as in the proof of Assertion (1), with at 
worst twice as many children. \qed}  

\subsection{The Algorithm that Proves Theorem \ref{thm:big}}\mbox{}\\  
\label{sub:mainalgor} 
Recall that a {\em terminal} node of a tree is a node with no children. 

\smallskip 
\noindent  
\scalebox{.95}[1]{\fbox{\mbox{}\hspace{.3cm}\vbox{
\begin{algor} {\em
\label{algor:trinosolqp}
{\bf (Solving Trinomial Equations Over $\pmb{\Q^*_p}$)}
\mbox{}\\
{\bf Input.} A prime $p$ and
$c_1,c_2,c_3,a_2,a_3\!\in\!\Z\setminus\{0\}$ with $|c_i|\!\leq\!H$ for all 
$i$ and $1\!\leq\!a_2\!<\!a_3\!=:\!d$. \\
{\bf Output.} \scalebox{.95}[1]{A true declaration that 
$f(x)\!:=\!c_1+c_2x^{a_2}+c_3x^{a_3}$ 
has no roots in $\Q_p$, or $z_1,\ldots,z_m\!\in\!\Q$}\linebreak 
\mbox{}\hspace{1.8cm}\scalebox{.92}[1]{with logarithmic
height $O\!\left(p^2 \log^4(dH)\right)$ such that $m$ is the number of 
roots of $f$ in $\Q_p$, $z_j$}\linebreak  
\mbox{}\hspace{1.8cm}\scalebox{.92}[1]{is an approximate root of 
$f$ with associated true root $\zeta_j\!\in\!\Q_p$ for all $j$, 
and $\#\{\zeta_j\}\!=\!m$.}\\
{\bf Description.} \\
1: If [$\ord_p\frac{c^2_2}{c_1c_3}\!
\geq\!0$ and $\ord_p c_1\!\neq\!\ord_p c_3$ mod $a_3$] or\\ 
\mbox{}\hspace{.5cm}[$\ord_p\frac{c^2_2}{c_1c_3}\!
<\!0$ and $\ord_p c_1\!\neq\!\ord_p c_2$ mod $a_2$ and $\ord_p c_2\!\neq\!
\ord_p c_3$ mod $a_3-a_2$]\\ 
\mbox{}\hspace{.5cm}then say {\tt ``No roots in $\Q_p$!''} and {\tt STOP}. \\
2: Rescale and invert roots if necessary, so that we may assume  
$p\nmid c_1 c_2$ and $\ord_p c_3\!\geq\!0$. \\  
3: Decide, via gcd-free bases, $\Delta_{\mathrm{tri}}(f)\!\stackrel{?}{=}\!0$. 
If so, set $\delta\!:=\!1$. Otherwise, set $\delta\!:=\!0$.\\ 
4: If $\delta\!=\!1$ then, via Algorithm \ref{algor:binoqp} (or its $p\!=\!2$ 
version, Algorithm \ref{algor:binoq2}), {\tt output} the $2$ most\\ 
\mbox{}\hspace{.5cm}significant base-$p$ digits of each degenerate root of 
$f$ in $\Z_p$ with valuation $0$.\\ 
5: Set $k$ to be the lower bound from Corollary \ref{cor:trinodepth} 
(employing the stated upper bound on $S_0$,\\
\mbox{}\hspace{.5cm}\scalebox{.98}[1]{and the upper bound 
on $D$ from Theorem \ref{thm:tri}, should $S_0$ or $D$ not be known), and  
compute}\\
\mbox{}\hspace{.5cm}the mod $p$ reductions $\tf_{i,\zeta}$ of all the 
nodal polynomials of $\cT_{p,k}(f)$.\\ 
6: By computing $\deg \gcd(\tf_{i,\zeta},x^p-x)$ for the non-root nodal 
polynomials of $\cT_{p,k}(f)$, and brute-\\
\mbox{}\hspace{.5cm}force search over $\F^*_p$ for $\tf$, determine 
which nodal polynomials have non-degenerate roots. \\  
7: {\tt Output} every non-degenerate root $\zeta_0\!\in\!\F_p$ of $\tf$. Also 
{\tt output}, for each non-root nodal\\ 
\mbox{}\hspace{.5cm}polynomial $f_{i,\zeta}$ found in Step 6, the 
set $\left\{\zeta+p^i\zeta_i\; | \; \zeta_i\!\in\!\F_p \text{ and } 
\tf_{i,\zeta}(\zeta_i)\!=\!0\!\neq\!\tf'_{i,\zeta}(\zeta_i)\right\}$.\\   
8: If $p|c_3$ then rescale and invert roots to compute approximants for  
the remaining roots of $f$\\ 
\mbox{}\hspace{.5cm}\scalebox{.95}[1]{in $\Q_p$, by computing roots of 
valuation $0$ for a rescaling of the reciprocal polynomial $f^*$.}} 
\end{algor}}
}} 
\begin{rem} {\em
We point out that some of the approximate roots output by our algorithm 
above require the use of Newton iteration {\em applied to $f_{i,\zeta}$
(instead of $f$)}. This is clarified in our correctness proof below. \dia }
\end{rem}

\smallskip 
\noindent 
{\bf Proof of Theorem \ref{thm:big}:} 
First note that the root $0$ is trivially detected by checking whether 
the constant term $c_1$ is $0$. So we may assume $c_1\!\neq\!0$ 
and focus on roots in $\Q^*_p$. Note also that the rescalings from 
Steps 2 and 8 (which are simply replacements of $f$ with 
$p^{j_1}f(p^{j_2}x)$ for suitable $j_1,j_2\!\in\!\Z$) 
result in a possible increase in the bit-sizes our outputs, 
but this increase is $O(\log H)$ thanks to Theorem \ref{thm:newt}. So 
we focus on roots in $\Z_p$ of valuation $0$, and assume $p\nmid c_1$ 
and $\ord_p(c_2)\ord_p(c_3)\!=\!0$. 

Condition (1) (the logarithmic height bound for our approximate roots) then 
clearly holds thanks to Step 5 of our algorithm, the definition 
of $\cT_{p,k}(f)$, Lemma \ref{lem:ulift}, Theorem \ref{thm:tri}, and 
Corollary \ref{cor:trinodepth}.  

Condition (2) (on the convergence of the Newton iterates) follows easily from 
the definition of $f_{i,\mu}$. In particular, Lemma \ref{lem:nodal} tells us 
that $f_{i,\mu}(x)\!=\!p^{-s}f(\mu+p^i x)$ mod $p^j$ for suitable $(s,\mu,j)$, 
and thus a non-degenerate root $\zeta_i\!\in\!\F_p$ of $\tf_{i,\zeta}$ 
yields a root $\mu+p^i\zeta_i$ of $f$ mod $p^{i+1}$. 
Moreover, by Hensel's Lemma, $z_0\!:=\!\zeta_i$ is an approximate root of 
$f_{i,\mu}$, meaning that the sequence $(\mu+p^i z_n)_{n\in\N}$ 
derived from the iterates $(z_n)_{n\in\N}$ coming from applying Newton 
iteration to $(f_{i,\mu},z_0)$ satisfies 
$|\xi-(\mu+p^i z_n)|_p\!\leq\!\left(\frac{1}{p}\right)^{2^{n-1}}
|\xi-(\mu+p^iz_0)|_p$, where $\xi\!\in\!\Z_p$ is some true (non-degenerate) 
root of $f$. From Lemma \ref{lem:ulift} (and our choice of $k$ via 
Corollary \ref{cor:trinodepth}) we know that {\em all} the 
non-degenerate roots of $f$ can be recovered this way, and uniquely so. 

Condition (3) on correctly counting the roots of $f$ in 
$\Q_p$ follows immediately from Steps 3--8. In particular, 
Step 4 correctly counts the degenerate roots in $\Q_p$ thanks to our 
earlier work on Algorithms \ref{algor:binoqp} and \ref{algor:binoq2}.  
Also, Corollary \ref{cor:trinodepth} and Lemma \ref{lem:ulift}  
tell us that the outputs from Step 7 are a collection of 
approximate roots that, en masse, converge to the set of 
non-degenerate roots of $f$ in $\Z_p$ of valuation $0$, with no overlap. 
Step 8 then accounts for the remaining degenerate and non-degenerate roots 
in $\Q_p$. 

The time complexity estimates from our theorem will follow from our complexity 
analysis of Algorithm \ref{algor:trinosolqp} below. 
First, however, let us prove correctness for our algorithm. 

\smallskip 
\noindent 
{\bf Correctness:} Via Theorem \ref{thm:newt}, Step 1  
guarantees that $f$ has roots of integral valuation, which is a 
necessary condition for their to be roots in $\Q_p$. 
Steps 2 and 8 involves substitutions that only negligibly affect 
the heights of the coefficients, similar to the binomial case (where 
the underlying rescalings are stated in finer detail). 

Step 3 correctly detects degenerate roots in $\C^*_p$ thanks to Lemma 
\ref{lemma:tri}. 
As observed above, Steps 4--7 correctly count the number of non-degenerate 
roots of $f$ in $\Z_p$ of valuation $0$. In particular, 
Step 4 is accomplished via Lemmata \ref{lemma:tri} and \ref{lemma:tri2}, 
and the characterization of degenerate roots from the latter lemma implies 
that we can use the Extended Euclidean Algorithm to find a binomial 
efficiently encoding the degenerate roots of $f$ in $\Q_p$ (as already 
detailed in the third paragraph of the proof of Lemma \ref{lem:cxity}). 
\qed 

\smallskip 
\noindent 
{\bf Complexity Analysis:} Steps 1, 2, and 8 involve basic field 
arithmetic that will be dominated by Steps 3--7. So we will focus on 
Steps 3--7 only. 

Step 3 can be accomplished in time $O(\log^2(dH))$ via \cite[Thm.\ 39]{airr}. 
Note in particular that detecting vanishing for $\Delta_{\mathrm{tri}}(f)$ 
is much easier than computing its valuation. 

Step 4 takes time $O((p+\log(dH))\log(dpH)\log\log(dpH))$ thanks to 
Theorem \ref{thm:binoqp}. 

Letting $\nu$ and $\cD$ respectively denote the number of degenerate roots of 
$\tf$ in $\F^*_p$ and the depth of $\cT_{p,k}(f)$, Step 5 takes time 
$O\!\left(\nu p^2\log^4(dH) 
\log^3_p(d)\log\left(p\log(dH)\right)\right)$ or\\ 
$O\!\left((p+\log d)\log(dp)\log\log(dp)+p\log^2(p)\log\log(p)\right.$\\ 
\mbox{}\hfill $\left. +\nu \log^2(dH)\log(d)\log_p\log(dH)
+\log(H)\log(dpH)\log\log(dpH) \right)$,\\ 
according as $\delta\!=\!0$ or 
$\delta\!=\!1$. This follows immediately from an elementary calculation, upon 
substituting the corresponding value of $k$ from Corollary 
\ref{cor:trinodepth} into Lemma \ref{lem:cxity}, 
using the fact that the depth $\cD$ is bounded from above by one of our 
two bounds from Theorem \ref{thm:tri}. 

The brute-force portion of Step 6 clearly takes time 
$O(p\log^2(p)\log\log p)$ via Theorem \ref{thm:cxity}. 
Lemma \ref{lem:cxity} tells us that $\cT_{p,k}(f)$ has $O(\nu \cD)$ nodes, 
and Lemma \ref{lem:low} tells us that each  
non-root nodal polynomials has mod $p$ reduction with degree $\leq\!4$. 
So the remaining multi-node gcd computation takes time $O(\nu 
\cD\cdot \log(p)\log \log p)$ via Theorem \ref{thm:cxity}. 
So the overall time for Step 6 is 
$O(p\left[\nu \log^2(dH)\log_p(d)+\log^2 p\right]\log\log p)$ or 
$O([p\log^2(p)+\nu \log(dH)]\log\log p)$, 
according as $\delta$ is $0$ or $1$, thanks to Theorem \ref{thm:tri}.   

As for Step 7, we already know the non-degenerate roots in $\F_p$ 
of $\tf$ from Step 6. For the remaining nodes, observe that Lemma 
\ref{lem:low} tells 
us that the mod $p$ reductions of the non-root nodal polynomials have 
degree at most $4$. Also, the root has $\nu$ children, each yielding 
a tree that is a chain with (at worst) one bifurcation. Furthermore, note that 
the presence of a non-degenerate root in $\F_p$ for $\tf_{i,\zeta}$ implies 
that $\tf_{i,\zeta}$ can have at most $1$ degenerate root in $\F_p$, meaning 
that its child will have degree at most $2$ by Lemma \ref{lem:nodal}. Finally, 
note that once a quadratic $\tf_{i,\zeta}$ has a non-degenerate root 
in $\F_p$, it can no longer have any children. In other words, we have shown 
that there can be at most $O(\nu)$ nodes having $\tf_{i,\zeta}$ 
possessing a non-degenerate root. 
Applying Shoup's deterministic factoring algorithm \cite{shoup} 
to the non-root nodal polynomials, we then see that finding 
the non-degenerate roots for our entire tree takes time 
$O(\nu\cdot p^{1/2}\log^2 p)$. 

In summary, we see that Step 5 dominates our overall 
complexity when $\delta\!=\!0$, yielding a bound of 
\fbox{$O\!\left(\nu p^2\log^4(dH)\log^3_p(d)\log(p\log(dH))\right)$}. When 
$\delta\!=\!1$, Steps 4, 5, and 7 dominate together, yielding an overall 
complexity bound of\\  
\fbox{\vbox{
$O\!\left((p+\log(dH))\log(dpH)\log\log(dpH)+p\log^2(p)\log\log(p)\right.$\\
\mbox{}\hfill $\left. +\nu[p^{1/2}\log^2(p)+\log^2(dH)\log(d)\log_p\log(dH)]
\right)$.}} 

\noindent 
Noting that $\nu\!\leq\!p-1$, we are done after an elementary calculation. \qed

\begin{rem} 
\label{rem:main} 
{\em A consequence of our proof is that it also contains a proof of 
the deterministic complexity bound of Corollary \ref{cor:degen}, since we 
included above the case where $f$ has a degenerate root. To get the 
Las Vegas randomized bound, we simply replace the fast deterministic
\linebreak 
\scalebox{.95}[1]{factoring 
algorithm from \cite{shoup} in Step 7 with the fast randomized factoring 
algorithm from \cite{ku08}. \dia}}   
\end{rem} 

\subsection{``Typical'' Exponents, Las Vegas, and a Combined Speed-Up}  
\label{sub:faster} 
For our final speed-ups we will make use of the fact that 
trinomials can only vanish on a small number of cosets in $\F^*_q$: Building 
on earlier results from \cite{cfklls00,bcr13,kel16}, Kelley and Owen 
proved \cite[Thm.\ 1.2]{kelleyowen} that $c_1+c_2x^{a_2}+c_3x^{a_3} 
\!\in\!\F_q[x]$, with $q$ a prime power, vanishes at no more than
$\floor{\frac{1}{2}+\sqrt{\frac{q-1}{r'}}}$ cosets of the size $r'$
subgroup of $\F^*_q$ (and nowhere else), where $r'\!=\!\gcd(a_2,a_3,q-1)$. In
particular, this bound is optimal for $\F_q$ an even degree extension of a
prime field. For $q$ {\em prime}, there is even computational evidence (for
all $q\!\leq\!292837$) that the number of such cosets might in fact no greater
than $2\log q$ \cite{cgrw}. 

It is easy to see that, for any fixed prime $p$, 
$\gcd(a_2a_3(a_3-a_2),(p-1)p)\!\leq\!2$ for a positive 
density subset of $(a_2,a_3)\!\in\!\N^2$. (Simply pick $a_2$ and 
$a_3$ to avoid certain arithmetic progressions depending on $p$ and 
the divisors of $p-1$.) So one can argue that a 
large fraction of trinomials over $\Z$ have $O(\sqrt{p})$ roots 
in $\F_p$ and, via Lemma \ref{lem:ulift}, $O(\sqrt{p})$ roots in $\Q_p$. 
A propos of this paucity of roots for ``most'' 
exponents, let us recall a useful 
trick that will allow us to significantly reduce the degree of a large 
fraction of trinomials over $\F_p$: Via a fast algorithm\linebreak 
\scalebox{.95}[1]{for the {\em Shortest Lattice Vector Problem} in $\Z^2$ 
(see, e.g., \cite{dpv11}), one can prove the following result:}   
\begin{lem} \label{lem:lattice} 
{\em \cite[Special Case of Lemma 1.11]{bcr13} 
Given any prime $p$, and $a_2,a_3\!\in\!\N$ with 
$0\!<\!a_2\!<\!a_3\!<\!p-1$ and $r'\!:=\!\gcd(a_2,a_3,(p-1)p)$, one can 
find within $\log^{O(1)}p$ bit operations an integer $e$ such 
that for all $i\in \{2,3\}$, $ea_i\!=\!m_i$ mod $p-1$ and $|m_i|\! 
\leq\!r'\sqrt{2(p-1)}$. \qed } 
\end{lem}

\medskip 
\noindent 
{\bf Proof of Corollary \ref{cor:faster}:} We follow the template of the 
proof of Theorem \ref{thm:big}, save for some key differences. 
The first main difference is that, 
under our assumptions, we can compute the tree $\cT_{p,k}(f)$ faster 
via Corollary \ref{cor:average} instead of Lemma \ref{lem:cxity}. We 
then need to compute the non-degenerate roots of all the nodal polynomials, 
so the next key difference is that we can use degree reduction to speed up 
this up at the root node. (The remaining nodes receive no further speed-up 
unless randomization is used.) 

So we merely need to recompute our complexity bounds. Recall 
that $\cD$ denotes the depth of the tree $\cT_{p,k}(f)$, and $\nu$ 
is the number of children of the root node (which for $k$ sufficiently 
large, is the number of degenerate roots of $\tf$). We note the 
changes to the complexity of Algorithm \ref{algor:trinosolqp} below, in both 
the {\em restricted root} case  (where we only seek root of the form 
$p^j+O(p^{j+1})$) and the {\em small gcd} case (where we assume 
$\gcd(a_2a_2(a_3-a_2),(p-1)p)\!\leq\!2$):  

\begin{quote} 
A.\ Step 4 can be sped up to deterministic time\\ 
\mbox{}\hspace{.5cm}$O(\log^2(p)\log\log(p) 
+ \log(\max\{d,p\})\log(\log\max\{d,p\})$\\ 
\mbox{}\hfill $+ \log(\max\{p,H\})\log(\log\max\{p,H\}))$\\ 
\mbox{}\hspace{.5cm}in the restricted root case; or deterministic time\\  
\mbox{}\hspace{.5cm}$O(p^{1/2}\log^2(p) 
+ \log(\max\{d,p\})\log(\log\max\{d,p\})$\\ 
\mbox{}\hfill $+ \log(\max\{p,H\})\log(\log\max\{p,H\}))$,\\ 
\mbox{}\hspace{.5cm}or Las Vegas randomized time\\  
\mbox{}\hspace{.5cm}$O\!\left(\log^{2+o(1)}(p) + 
\log(\max\{d,p\}) \log\log \max\{d,p\})\right.$\\ 
\mbox{}\hfill $\left.+\log(\max\{H,p\})\log\log\max\{H,p\}
\vphantom{\cD^{\cD^\cD}}\right)$.\\  
\mbox{}\hspace{.5cm}in the small gcd case. 

\smallskip 
\noindent 
B.\ Step 5 can be sped up to deterministic time\\  
\mbox{}\hspace{.5cm}$O\!\left((p+\log d)\log(dp)\log(\log(dp))
   +p\log^2(p)\log\log(p)\right.$\\
\mbox{}\hspace{1.2cm}$+\left. \cD [k\log(p)\log(k\log p)\log(d) 
+\log H]+\log(H)\log(dpH)\log\log(dpH) \right)$,\\
\mbox{}\hspace{.5cm}\scalebox{.97}[1]{in both cases. If $f$ has a degenerate 
root in $\C^*_p$ then we can further speed up both}\\ 
\mbox{}\hspace{.5cm}cases to Las Vegas randomized time\\  
\mbox{}\hspace{1.2cm}
$O(\log^2_p(dH)\log(\log(dH))+\log^2(p)\log(\log p)+\log(dpH)\log\log(dpH))$.

\smallskip 
\noindent 
C.\ We replace Step 6 of Algorithm \ref{algor:trinosolqp} with the following:\\ 
\mbox{}\hspace{.6cm}\scalebox{.87}[1]{\fbox{\vbox{\noindent {\em 
6': \scalebox{1}[1]{By computing $\deg \gcd(\tf_{i,\zeta},x^p-x)$ for the 
non-root nodal polynomials of $\cT_{p,k}(f)$, and}\\ 
\mbox{}\hspace{.6cm}\scalebox{.98}[1]{factoring a degree-reduced 
version of $\tf$ (if needed), determine which nodal polynomials}\\
\mbox{}\hspace{.6cm}have non-degenerate roots in $\F_p$.}}}}\\ 
\mbox{}\hspace{.6cm}\scalebox{.95}[1]{This modified step takes deterministic 
time $O(\cD\log(p)\log\log p)$ in the restricted}\linebreak 
\mbox{}\hspace{.6cm}\scalebox{1}[1]{root case; or deterministic time 
$O(p\log^2(p)+\cD\log(p)\log\log p)$ or Las Vegas}\\ 
\mbox{}\hspace{.6cm}\scalebox{1}[1]{randomized time 
$O(p^{3/4}\log^{1+o(1)}(p)+\cD\log(p)\log\log p)$ in the small gcd case.}   

\smallskip 
\noindent 
D.\ \scalebox{1}[1]{Step 7 can be sped up to deterministic time 
$O(p^{1/2}\log^{2+o(1)}(p))$ or Las Vegas}\\
\mbox{}\hspace{.6cm}randomized time $O(\log^{2+o(1)} p)$, in both cases. 
\end{quote}
We now explain Changes A--D. 

\smallskip 
\noindent  
{\bf A.} In the restricted root case, we merely need to evaluate 
$f$ and $f'$ at $1$, so our first bound is clear. 

In the small gcd case, the number of degenerate 
roots is at most $2$ thanks to our gcd assumption and Lemma \ref{lemma:tri2}. 
So instead of employing Algorithms \ref{algor:binoqp} or \ref{algor:binoq2}, 
we simply find the degenerate roots by factoring, using either the 
fast deterministic algorithm from \cite{shoup} or the fast 
Las Vegas randomized factorization algorithm from \cite{ku08}.  

\smallskip 
\noindent 
{\bf B.} The complexity bounds follows by applying Corollary \ref{cor:average} 
instead of Lemma \ref{lem:cxity}, ultimately yielding 
$O(p^2\log^4(dH)\log^3_p(d)\log(p\log(dH)))$ via Corollary 
\ref{cor:trinodepth} and Theorem \ref{thm:tri}. As noted in Remark 
\ref{rem:bottle}, our current bounds 
for $k$ and $\cD$ obstruct any Las Vegas speed-up for Step 5 (in the 
non-degenerate case). 

\smallskip 
\noindent 
{\bf C.} The deterministic speed-ups follow from the complexity 
analysis of Algorithm \ref{algor:trinosolqp}, in the proof of 
Theorem \ref{thm:big}, simply by setting $\nu\!=\!2$ in the bound there. 
Note also that in the restricted root case, there is no need to 
search for any roots of $\tf$ since we only care about most significant 
digit $1$: We merely need to evaluate $\tf$ and $\tf'$ at $1$. 

To get our Las Vegas speed-up, we replace the brute-force search for 
degenerate roots of $\tf$ with a targeted factorization: 
First build a degree-reduced version of 
$\tf$ via Lemma \ref{lem:lattice} to apply the automorphism of 
$\F^*_p$ defined by $x\mapsto x^e$ to replace $\tf$ by 
$\tg(x)\!:=\!\tf(x^e)$, and compute $e'\!:=\!1/e$ mod
$p-1$, in deterministic time $\log^{O(1)}p$. This reduces $\deg \tf$ to 
$\deg \tg\!\leq\!2\sqrt{2(p-1)}$. 
To find the roots of $\tf$ in $\F^*_p$ we can then find the roots of 
$\tg$ in $\F^*_p$ by using the Kedlaya-Umans factorization algorithm 
\cite{ku08}, take the $e'$th powers mod $p$ of 
these roots, and then identify which of these roots of 
$\tf$ is a degenerate root found earlier. This takes time 
$O((2\sqrt{2(p-1)})^{1.5}\log^{1+o(1)}(p)
 +\log^{O(1)} p)\!=\!O(p^{3/4}\log^{1+o(1)}p)$. 

Since $\nu\!\leq\!2$ in both cases, the remaining multinodal gcd 
computation takes additional deterministic time $O(\cD\log(p)\log\log p)$. 

\smallskip 
\noindent 
{\bf D.} Since we already found the non-degenerate roots of 
$\tf$ in $\F_p$ in Step 6', we merely need to speed up finding the 
non-degenerate roots in $\F_p$ of the remaining nodal polynomials: 
We already observed in the proof of Theorem \ref{thm:big} that there are 
$O(\nu)$ nodes having a $\tf_{i,\zeta}$ possessing a non-degenerate root. 
But $\nu\!\leq\!2$ in both cases, so we only need to worry about $O(1)$ 
nodes. So our proof of Theorem \ref{thm:big} already implies a 
deterministic speed-up to $O(p^{1/2}\log^2p)$ (for $O(1)$ applications 
of Shoup's deterministic factoring algorithm \cite{shoup}), in both cases. 

However, if we replace Shoup's algorithm with the fast randomized 
factorization algorithm from \cite{ku08}, then we can speed Step 7 up to  
Las Vegas randomized time $O(\log^{2+o(1)} p)$ in both cases.

To conclude, we see that Step 5 dominates the deterministic complexity in both 
cases (restricted root and small gcd), and wipes out any Las Vegas speed-up 
unless better bounds\linebreak 
\scalebox{.95}[1]{for $k$ and $\cD$ are available. 
Summing our complexity estimates, we obtain our desired bounds. \qed}  

\medskip
An immediate consequence of our last proof --- if we can apply the sharper 
bounds for $\cD$ and $k$ from the degenerate cases of Theorem \ref{thm:tri} 
and Corollary \ref{cor:trinodepth} --- is the following combined speed-up: 
\begin{cor} 
\label{cor:final} {\em 
Following the notation of Corollary \ref{cor:degen}, we can 
speed up the Las Vegas complexity bound to 
$O\!\left(\log^2(p)\log(\log p)+\log^2(dH)\log_p(d)\log(\log(dH))\right)$ 
(in the restricted root case) or 
$O\!\left(p^{3/4}\log^{1+o(1)}(p)
+\log^2(dH)\log_p(d)\log(\log(dH))\right)$ 
(in the small gcd case). \qed}  
\end{cor} 

\section*{Acknowledgements}   
We thank Erich Bach and Bjorn Poonen for informative discussions
on Hensel's Lemma. Special thanks to Elliott Fairchild, Josh Goldstein, and 
David Zureick-Brown for inspirational conversations during the pandemic.   
We also thank the anonymous referees for helpful suggestions that improved our 
paper. 

\bibliographystyle{plain}
\bibliography{cha}

\end{document}